\documentclass[preprint]{elsarticle}

\usepackage{etex}               

\usepackage[bookmarks=true,colorlinks=true,linkcolor=blue]{hyperref}

\usepackage[margin=1.in]{geometry}

\usepackage{xcolor}
\definecolor{jr@red}{RGB}{228,26,28}
\definecolor{jr@blue}{RGB}{55,126,184}
\definecolor{jr@green}{RGB}{77,175,74}
\definecolor{jr@purple}{RGB}{152,78,163}
\definecolor{jr@orange}{RGB}{255,127,0}
\definecolor{jr@yellow}{RGB}{255,255,51}
\definecolor{jr@brown}{RGB}{166,86,40}
\definecolor{jr@pink}{RGB}{247,129,191}
\definecolor{jr@gray}{RGB}{153,153,153}

\usepackage{amsmath}          
\usepackage{amssymb}          
\usepackage{mathrsfs}         
\usepackage{mathtools}        

\usepackage{soul}             
\usepackage[normalem]{ulem}

\usepackage{caption}          
\usepackage{makecell,multirow,multicol,booktabs} 

\usepackage[group-separator={,}]{siunitx}

\usepackage{tikz}
  \usetikzlibrary{calc,positioning}
  \usetikzlibrary{arrows.meta}
\usepackage{tikz-3dplot}

\usepackage{subcaption}

\usepackage{pgfplots}
  \pgfplotsset{compat=newest}
\usepackage{graphicx}
\usepackage[all]{xy}
\usepackage{lineno}

\usepackage{hyperref}

\usepackage{amsthm}
\theoremstyle{plain}
\newtheorem{lemma}{Lemma}[section]

\newtheorem{proposition}[lemma]{Proposition}
\theoremstyle{definition}

\theoremstyle{remark}
\newtheorem{remark}[lemma]{Remark}
\usepackage{bm}

\usepackage[section]{algorithm} 
\usepackage{algpseudocode}  

\makeatletter%
\newcounter{algorithmicH} 
\let\oldalgorithmic\algorithmic
\renewcommand{\algorithmic}{%
  \stepcounter{algorithmicH} 
  \oldalgorithmic} 
\renewcommand{\theHALG@line}{ALG@line.\thealgorithmicH.\arabic{ALG@line}}
\makeatother%

\newcommand{\boldvec}[1]{\ensuremath{\mathbf{#1}}}

\ifdefined\vec
  \renewcommand{\vec}[1]{\boldvec{#1}}
\else
  \newcommand{\vec}[1]{\boldvec{#1}}
\fi

\makeatletter
\newcommand{\editHighlighting}[3]{%
  \if@display%
  \textcolor{red!40!white}{\text{\sout{#1}} \textcolor{#2}{#3}}%
  \else%
  \textcolor{red!40!white}{\sout{#1} \textcolor{#2}{#3}}%
  \fi%
}
\makeatother

\newif\ifshowstatus
\showstatustrue 

\usepackage{graphicx}
\graphicspath{{figures/}}

\usepackage{placeins}

\begin{document}

\begin{frontmatter}

\title{LGNO: A Local--Global Neural Operator for Hyperbolic Conservation Laws}

\author[gt]{Hao Wang}
\ead{wanghaomathe@gmail.com}

\author[brown]{Chi-Wang Shu\fnref{brownThanks}}
\ead{chi-wang_shu@brown.edu}

\author[gt]{Qi Tang\corref{cor1}\fnref{ascrThanks}}
\ead{qtang@gatech.edu}

\address[gt]{School of Computational Science and Engineering, Georgia Institute of Technology, Atlanta, GA 30332, USA.}
\address[brown]{Division of Applied Mathematics, Brown University, Providence, Rhode Island 02912, USA.}

\cortext[cor1]{Corresponding authors}
\fntext[brownThanks]{Research was partially supported by NSF grant DMS-2309249.}
\fntext[ascrThanks]{Research was partially supported by the U.S.\ Department of Energy Office of Science, Early Career Research Program under Award Number DE-SC0026277 and and Office of Naval Research under award number N00014-23-1-2545 (Untangling Computation).}
  
\begin{abstract}
Solutions of hyperbolic conservation laws exhibit both smooth structures across large scales and sharp localized features such as shocks and contact discontinuities, making them difficult to approximate accurately with existing neural operators. The Fourier Neural Operator (FNO) captures long-range interactions well but tends to smear localized structures through excessive numerical dissipation. To address this, we propose a Local--Global Neural Operator (LGNO) that learns a one-step discrete flow map by combining a global FNO branch for representing smooth dynamics at large scales with a local multiresolution branch for enhancing localized discontinuities and nonsmooth features. The model is trained with a one-step loss that combines a physical space prediction term and a spectral penalty on high frequencies to suppress spurious oscillations near steep fronts. On a large collection of benchmarks in one and two dimensions, LGNO consistently outperforms FNO baselines with matched parameter counts, reducing one-step errors by factors of 2--5 and remaining significantly more accurate over long autoregressive rollouts. Most strikingly, although it is trained only on short-time data from a high-order WENO-Z scheme, the long-time rollout of LGNO on a coarse \(256^2\) grid exhibits lower numerical dissipation than the same WENO-Z scheme run on a finer \(512^2\) grid, while being orders of magnitude cheaper to evaluate. These results suggest that, with an appropriate architecture and training objective, learned operators can effectively learn discrete flow maps. They further suggest that such learned operators have the potential to control long-time numerical dissipation better than the conventional shock-capturing schemes that generate the training data.
\end{abstract}


\end{frontmatter}


\section{Introduction}
\label{sec:intro}

Conservation laws arise in a broad range of applications, including gas dynamics~\cite{cauchy_conservationlaw,fluid_conservationlaw} and plasma physics such as magnetohydrodynamics~\cite{plasma_conservationlaw,christlieb2014finite}. Many of these applications are governed by hyperbolic systems of partial differential equations (PDEs)~\cite{Leveque_ConservationLaws,Leveque_ConservationLaws_fv}. A defining feature of such systems is that even smooth initial data can develop sharp gradients and discontinuities in finite time~\cite{dafermos2005hyberbolic}. This makes their numerical approximation fundamentally challenging: a successful method must accurately propagate nonlinear waves, sharply resolve localized structures such as shocks and contact discontinuities, and remain robust over long-time integration.

A large body of conventional numerical methods has been developed for conservation laws, including finite volume methods~\cite{Leveque_ConservationLaws_fv}, discontinuous Galerkin methods~\cite{Nodal_DG,OEDG}, and high order nonoscillatory reconstruction techniques such as the WENO scheme~\cite{WENO_init} and its variants such as WENO-Z~\cite{WENO_improve} and TENO~\cite{WENO_Hyperbolic_LINFU}. While these methods differ in formulation, their effectiveness is closely tied to the structure of the underlying equations~\cite{Morton2007FV,Shu2016WENODG}. In particular, shock-dominated dynamics is governed by strongly localized interactions, and classical high-resolution schemes are carefully designed around local reconstruction, numerical fluxes, and nonlinear stabilization mechanisms. These ingredients allow them to balance accuracy in smooth regions with robustness near discontinuities~\cite{SHU2020ENOWENO}. However the robustness comes at the price of numerical dissipation: reducing this dissipation while preserving fine-scale structures typically requires higher-order reconstructions, less dissipative stabilization mechanisms, or finer meshes, each of which substantially increases the computational cost.

In recent years, machine learning methods for PDEs have made rapid progress, especially in neural surrogate models for solution prediction~\cite{model_surrogate} and in operator learning approaches such as the Fourier Neural Operator (FNO)~\cite{li2021fourier} and the Deep Operator Network (DeepONet)~\cite{deeponet}. By learning the solution operator or discrete flow map directly, these approaches offer a promising route for accelerated simulation in repeated-query settings~\cite{flow_map_matching,Neural_operator_speed}. Studies in approximation theory indicate that the complexity of operator learning depends strongly on the regularity of the target map, with smoother operators generally being easier to approximate efficiently~\cite{parametric_complexity}. For nonlinear conservation laws, however, shocks and other discontinuities naturally introduce nonsmooth solution structures, which makes accurate operator learning substantially more difficult~\cite{neural_operator_struggle}.

While neural operators with nonlinear reconstruction, such as FNO, can in principle represent discontinuous solutions~\cite{Nonlinear_reconstruction}, representational capacity alone does not guarantee practical performance. Existing studies have shown that FNO can handle relatively simple discontinuous examples, but their accuracy deteriorates as the discontinuity structures become more complex~\cite{neural_operator_struggle}. Related neural-operator approaches for discontinuous Riemann problems in the compressible Euler equations have also shown promising results, but their accuracy near shocks and contact discontinuities remains limited~\cite{RiemannOnets}. This gap has motivated a growing line of structure-aware learning methods that incorporate ideas from classical shock-capturing schemes, for example by learning local numerical fluxes or neural finite-volume structures~\cite{FLUXFNO,UNFV}, indicating that locality and conservation structure are useful inductive biases for learning nonsmooth dynamics.

Motivated by this perspective, this work introduces a Local--Global Neural Operator (LGNO) for hyperbolic conservation laws that learns a one-step discrete flow map by coupling two complementary components: a global branch, built on an FNO, that captures coherent large-scale behavior and long-range interactions, and a local multiresolution convolutional branch that sharpens localized nonsmooth structures. The model is trained with a one-step loss combining a physical-space prediction term and a high-frequency spectral penalty that suppresses spurious oscillations. This design aligns the architecture with the coexistence of smooth global waves and localized nonsmooth features that characterizes conservation law solutions.

\begin{figure}[htb]
    \centering
    \begin{subfigure}[b]{0.24\textwidth}
        \centering
        \includegraphics[width=\textwidth]{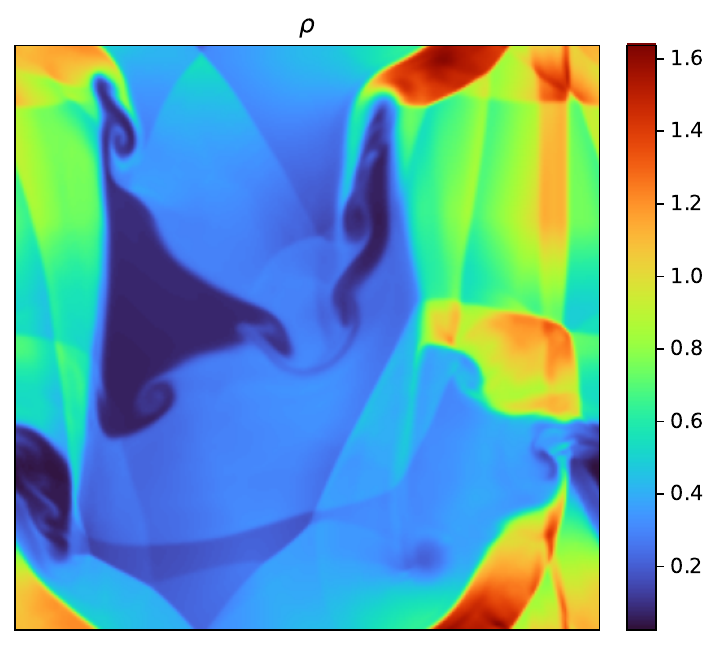}
        \caption{WENO-Z \(256^2\)}
    \end{subfigure}
    \begin{subfigure}[b]{0.24\textwidth}
        \centering
        \includegraphics[width=\textwidth]{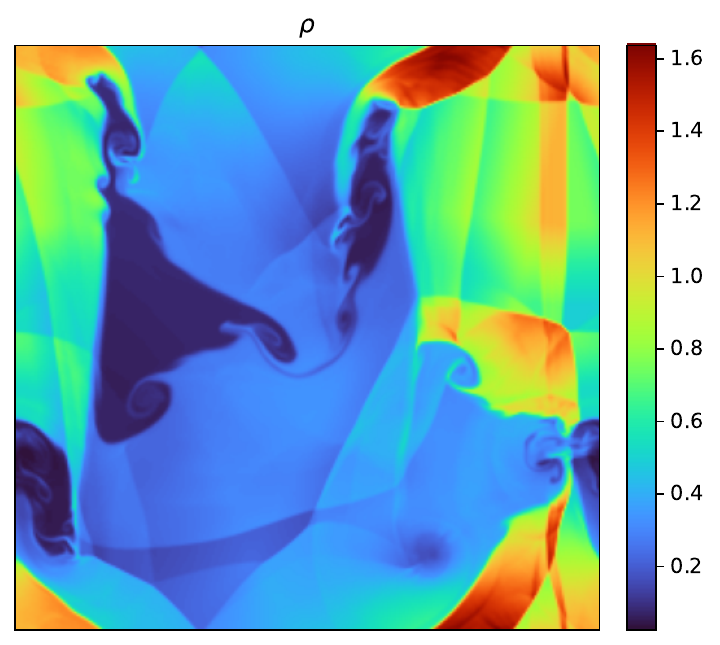}
        \caption{WENO-Z \(512^2\)}
    \end{subfigure}
    \begin{subfigure}[b]{0.24\textwidth}
        \centering
        \includegraphics[width=\textwidth]{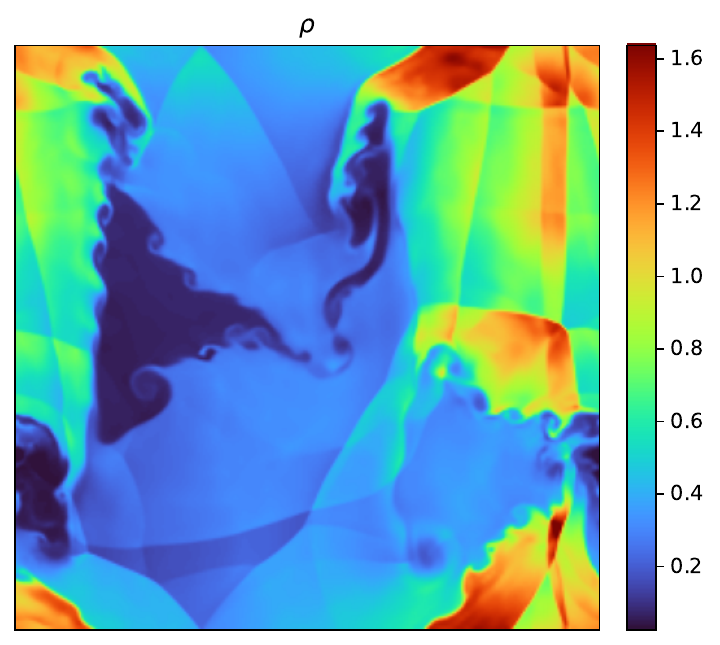}
        \caption{LGNO \(256^2\)}
    \end{subfigure}
    \begin{subfigure}[b]{0.24\textwidth}
        \centering
        \includegraphics[width=\textwidth]{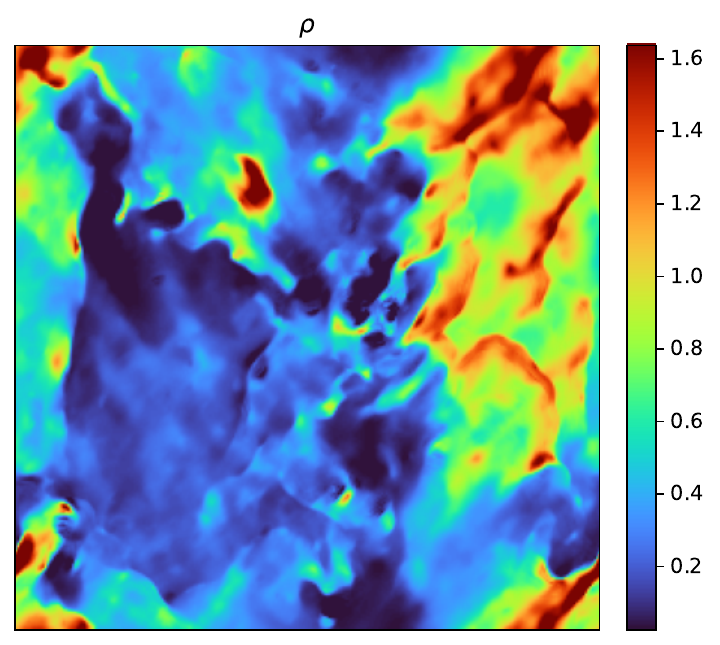}
        \caption{FNO \(256^2\)}
    \end{subfigure}
    \caption{\textbf{LGNO is less dissipative than WENO-Z at a higher resolution.} Density field for the two-dimensional Euler equations after \(50\) autoregressive rollout steps with \(\Delta t=0.01\). Trained only on short-time WENO-Z data at \(512^2\) resolution, the learned LGNO on the coarse \(256^2\) grid (c) preserves more vortical structures along the contact discontinuities and shear layers than the WENO-Z reference at both \(256^2\) (a) and \(512^2\) (b), reflecting lower long-time numerical dissipation, whereas the FNO baseline at \(256^2\) (d) is severely degraded over the long rollout and loses the coherent structure.}
    \label{fig:intro}
\end{figure}

The main contributions of this work are threefold. First, we develop a local multiresolution branch that complements the global FNO branch, together with a simple learned boundary padding for outflow boundary conditions. Second, through an extensive study on standard one- and two-dimensional hyperbolic benchmarks, we show that LGNO captures shock and vortical structures significantly more accurately than a parameter-matched FNO. Third, and most notably, we demonstrate that, although it is trained only on short-time data from the WENO-Z scheme, LGNO on a coarse \(256^2\) grid produces long-time rollouts that are less dissipative than the same scheme run on a finer \(512^2\) grid, while being orders of magnitude cheaper to evaluate, as illustrated in Figure~\ref{fig:intro}.

This may appear counterintuitive, since the training data are produced by the shock-capturing scheme itself; the reason is that such schemes are accurate over short times but accumulate excessive numerical dissipation when their dissipative reconstruction acts at every step over a long integration, whereas the learned one-step flow map of LGNO inherits only the lower dissipation of the short-time training data as it is composed over long rollouts. This indicates that, with an appropriate architecture and training procedure, learned operators have the potential to control long-time numerical dissipation better than the shock-capturing schemes that generate their training data.

The remainder of this paper is organized as follows. Section~\ref{sec:background} reviews hyperbolic conservation laws and operator learning. Section~\ref{sec:methods} presents the proposed LGNO and its data-driven training strategy. Section~\ref{sec:results} evaluates LGNO on several benchmark problems against parameter-matched FNO baseline and the \mbox{WENO-Z} scheme, and Section~\ref{sec:conc} concludes the paper. Implementation details, computational cost comparisons, additional numerical results, and benchmark setups are provided in the appendices.


\section{Background and Related Work}\label{sec:background}

\subsection{Hyperbolic Conservation Laws and Conventional Numerical Schemes}\label{subsec:background_conservation}

We consider systems of hyperbolic conservation laws posed on a spatial domain \( \Omega \subset \mathbb{R}^d \):
\begin{equation}
\label{eq:conslaw_general}
\frac{\partial \mathbf{u}}{\partial t}
+\sum_{\alpha=1}^{d}\frac{\partial \mathbf{f}_{\alpha}(\mathbf{u})}{\partial x_{\alpha}}
=0, \quad \mathbf{x}\in \Omega, \quad t > 0.
\end{equation}
Here \( \mathbf{u}=(u_1,\dots,u_m)^{\top}\in \mathbb{R}^m \) denotes the vector of conserved variables, and \( \mathbf{f}_{\alpha}:\mathbb{R}^m\to\mathbb{R}^m \), \(\alpha=1,\dots,d\), are the Cartesian components of the flux. The system is supplemented with an initial condition
\[
\mathbf{u}(\mathbf{x},0)=\mathbf{u}_0(\mathbf{x}), \quad \mathbf{x}\in\Omega,
\]
together with suitable boundary conditions on \(\partial\Omega\).

A defining feature of nonlinear hyperbolic conservation laws is that even smooth initial data may develop steep gradients and discontinuities in finite time~\cite{dafermos2005hyberbolic}, so that classical solutions generally fail to exist globally and one must work with weak solutions. In integral form, Eq.~\eqref{eq:conslaw_general} implies the conservation relation
\begin{equation}
\frac{\mathrm{d}}{\mathrm{d}t}\int_{K}\mathbf{u}(\mathbf{x},t)\,\mathrm{d}\mathbf{x}
+
\int_{\partial K}\mathbf{F}(\mathbf{u})\cdot \mathbf{n}\,\mathrm{d}s
=0,
\label{eq:integral_conservation}
\end{equation}
for every control volume \(K\subset\Omega\), where
\(
\mathbf{F}(\mathbf{u})
=
(\mathbf{f}_1(\mathbf{u}),\dots,\mathbf{f}_d(\mathbf{u}))
\)
and \(\mathbf{n}\) is the outward unit normal. Since weak solutions are in general nonunique, additional admissibility conditions, such as entropy inequalities, are required to select the physically relevant solution~\cite{Leveque_ConservationLaws}.

From a numerical standpoint, the solution structure of conservation laws is heterogeneous. In smooth regions, global high-order approximations are extremely effective. For periodic problems with sufficiently smooth solutions, for instance, Fourier spectral methods approximate the solution by a truncated global expansion and satisfy the standard estimate
\[
\|\mathbf{u}-\mathbf{u}_N\|_{L^2(\Omega)} \lesssim N^{-s}\|\mathbf{u}\|_{H^s(\Omega)},
\]
when \(\mathbf{u}(\cdot,t)\in H^s(\Omega)\), with exponential convergence for analytic solutions~\cite{spectral_method}. This advantage, however, relies on regularity. Once shocks or contact discontinuities emerge, the loss of regularity is localized near them, whereas the Fourier basis remains globally supported. Spectral convergence is then lost, and truncated expansions develop Gibbs oscillations whose amplitude stays \(O(1)\) near the discontinuity even as the resolution increases~\cite{Gibbs_filter}. Global spectral methods thus capture smooth coherent structures and long-range interactions efficiently, but resolve sharply localized nonsmooth features poorly.

Classical shock-capturing methods are designed precisely for this regime. Finite volume methods build directly on Eq.~\eqref{eq:integral_conservation}. In one space dimension, let \(I_j=[x_{j-\frac{1}{2}},x_{j+\frac{1}{2}}]\) be a cell of width \(\Delta x\), and let
\[
\overline{\mathbf{u}}_j(t)
=
\frac{1}{\Delta x}\int_{I_j}\mathbf{u}(x,t)\,\mathrm{d}x
\]
denote the cell average. A semi-discrete conservative update takes the form
\begin{align}
\frac{\mathrm{d}}{\mathrm{d}t}\overline{\mathbf{u}}_j(t)
=
-\frac{1}{\Delta x}
\Bigl(
\widehat{\mathbf{f}}_{j+\frac{1}{2}}
-
\widehat{\mathbf{f}}_{j-\frac{1}{2}}
\Bigr),\label{eqn:conservative_form}
\end{align}
where \(\widehat{\mathbf{f}}_{j+\frac{1}{2}}\) is a numerical flux at the cell interface. Godunov-type methods, approximate Riemann solvers, and high-resolution finite-volume schemes differ mainly in how these interfacial fluxes are reconstructed, evaluated, and stabilized~\cite{Leveque_ConservationLaws_fv,Morton2007FV}.

Among high-order shock-capturing methods, ENO and WENO schemes exemplify this local, nonlinear design~\cite{SHU2020ENOWENO}. They build interface values from local stencils and adapt the reconstruction nonlinearly to the local smoothness, yielding high-order accuracy where the solution is smooth while biasing away from nonsmooth stencils near discontinuities to suppress spurious oscillations and robustly capture localized wave structures. Unlike global spectral approximations, their effectiveness comes from localized reconstruction and nonlinear adaptation to the loss of regularity rather than from global smoothness.

This contrast suggests that learned operators for conservation laws should not rely on global spectral information alone, but should combine global mechanisms for smooth coherent structures with local mechanisms for sharply localized nonsmooth features~\cite{Morton2007FV,Shu2016WENODG,SHU2020ENOWENO}.

\subsection{Neural Operators and Operator Learning}\label{sec:neural_operator_background}

Neural operators approximate mappings between function spaces. Let \(\mathcal{X}\) and \(\mathcal{Y}\) be Banach spaces of functions on the spatial domain \(\Omega \subset \mathbb{R}^d\); one seeks to approximate an operator
\begin{equation}
\label{eq:operator_general}
\mathcal{G} : \mathcal{X} \to \mathcal{Y},
\quad
\mathbf{a}(\mathbf{x}) \mapsto \mathbf{u}(\mathbf{x}),
\end{equation}
where \(\mathbf{a}(\mathbf{x})\) is an input function and \(\mathbf{u}(\mathbf{x})\) the corresponding output. For conservation laws, it is often natural to take function spaces adapted to weak solutions, such as subsets of \(L^1(\Omega;\mathbb{R}^m)\) or \(BV(\Omega;\mathbb{R}^m)\).

In practice, neural operators act not on the raw input and output spaces but through higher-dimensional latent feature spaces, giving the common architecture
\[
\mathcal{G}_{\theta}
=
Q \circ \mathcal{L}_{L-1} \circ \cdots \circ \mathcal{L}_{0} \circ P.
\]
Here
\(
P:\mathcal{X}\to \mathcal{Z}_{0}
\)
is a lifting map that embeds the input function into a latent feature space,
\(
\mathcal{L}_{\ell}:\mathcal{Z}_{\ell}\to\mathcal{Z}_{\ell + 1}
\)
are operator layers acting on latent features, and
\(
Q:\mathcal{Z}_{L}\to\mathcal{Y}
\)
is a projection map that returns the latent representation to the output space. The kernel-inspired operator structure resides mainly in the intermediate layers \(\mathcal{L}_{\ell}\), while \(P\) and \(Q\) supply the feature lifting and projection needed for an expressive neural parameterization.

For many PDEs, the solution map can be viewed as an operator induced by an integral kernel. In the linear setting, one may write
\[
\mathbf{u}(\mathbf{x}) = \int_{\Omega} K(\mathbf{x},\mathbf{y})\mathbf{a}(\mathbf{y})\,\mathrm{d}\mathbf{y},
\]
where \(K\) plays the role of a Green's kernel. This motivates neural operator layers that model feature transformations through learned kernel interactions and a nonlinear activation, giving the generic form
\[
(\mathcal{L}_{\ell} \mathbf{v})(\mathbf{x})
=
\sigma_{\ell} \left( W_{\ell}\mathbf{v}(\mathbf{x}) + (\mathcal{K}_{\ell}\mathbf{v})(\mathbf{x}) \right),
\]
where \(W_{\ell}\) is a pointwise linear map on channel features, \(\sigma_{\ell}\) a nonlinear activation, and \(\mathcal{K}_{\ell}\) an integral operator of the form
\[
(\mathcal{K}_{\ell} \mathbf{v})(\mathbf{x})
=
\int_{\Omega}\kappa_{\ell}(\mathbf{x},\mathbf{y};\theta)\,\mathbf{v}(\mathbf{y})\,\mathrm{d}\mathbf{y}.
\]
This kernel formulation abstracts many neural operator architectures, particularly those built on global integral or spectral representations.

Among these, the FNO is one of the most influential. In translation-invariant settings it replaces the kernel integral by a spectral convolution: with \(\mathbf{v}_{\ell}\) the latent state at layer \(\ell\), the update takes the form
\[
\mathbf{v}_{\ell+1}(\mathbf{x})
=
\sigma \left( W_{\ell}\mathbf{v}_{\ell}(\mathbf{x}) + \mathcal{F}^{-1} \bigl( R_{\ell}(\mathbf{k})\,\widehat{\mathbf{v}}_{\ell}(\mathbf{k}) \bigr)(\mathbf{x}) \right),
\]
where \(\widehat{\mathbf{v}}_{\ell}=\mathcal{F}(\mathbf{v}_{\ell})\) denotes the Fourier transform of \(\mathbf{v}_{\ell}\), \(\mathcal{F}^{-1}\) is the inverse transform, and \(R_{\ell}(\mathbf{k})\) is a learnable matrix-valued multiplier in Fourier space. In practice, \(R_{\ell}(\mathbf{k})\) is nonzero only on a prescribed band of low-frequency modes, so the FNO acts as a learned global spectral low-pass filter alongside a local pointwise linear path. This is very efficient for smooth problems and for dynamics dominated by large-scale coherent structures, since the Fourier basis captures long-range interactions naturally. Because the representation is global and truncated in frequency, however, its inductive bias is poorly aligned with the localized nonsmooth structures of shock-dominated conservation laws.

For time-dependent PDEs it helps to distinguish the solution operator from the discrete flow map. Let \(\mathbf{u}(\cdot,t)\) solve the conservation law with initial data \(\mathbf{u}_0\); the continuous solution operator at time \(t\) is
\[
\mathcal{S}_t:\mathcal{X}\to\mathcal{X},
\quad
\mathbf{u}_0 \mapsto \mathbf{u}(\cdot,t).
\]
Equivalently, for a time increment \(\Delta t\), the associated flow map is written formally as
\[
\mathcal{G}_{\Delta t}:\mathcal{X}\to\mathcal{X},
\quad
\mathbf{u}^{n}\mapsto \mathbf{u}^{n+1}.
\]
Operator learning may target either the full solution operator \(\mathcal{S}_t\) or the one-step flow map \(\mathcal{G}_{\Delta t}\).

In this work we adopt the latter viewpoint, with a fixed \(\Delta t\). For a mesh with spacing \(h\), let \(\mathcal{X}_h\) be the finite-dimensional space of grid functions; the exact discrete flow map is
\[
\mathcal{G}_{\Delta t,h}:\mathcal{X}_h\to\mathcal{X}_h,
\quad
\mathbf{u}_h^n\mapsto \mathbf{u}_h^{n+1}.
\]
This grid-function perspective matches time advancement in conventional numerical schemes and is the setting of our numerical experiments. For brevity we drop the subscript \(h\) when no ambiguity arises.

A simple residual parameterization of the learned flow map is
\begin{equation}
\mathcal{G}_{\theta}(\mathbf{u})
=
\mathbf{u}
+
\Delta t\,\Phi_{\theta}(\mathbf{u};\Delta t),
\label{eq::flow_update}
\end{equation}
where \(\Phi_{\theta}\) models the increment of the discrete solution over one time step \(\Delta t\).

Given one-step training data
\[
\mathcal{D}
=
\{(\mathbf{u}_i,\mathbf{v}_i)\}_{i=1}^{N},
\quad
\mathbf{v}_i = \mathcal{G}_{\Delta t}(\mathbf{u}_i),
\]
generated by a high-order, well-resolved numerical scheme \(\mathcal{G}_{\Delta t}\) (whose internal time step due to CFL is typically tens to a hundred times smaller than the sampled increment \(\Delta t\)), we seek a neural approximation \(\mathcal{G}_{\theta}\) of this discrete flow map. The parameters \(\theta\) are determined by minimizing
\[
\mathcal{L}(\theta)
=
\frac{1}{N}
\sum_{i=1}^{N}
d_{\mathcal{X}}
\Bigl(
\mathcal{G}_{\theta}(\mathbf{u}_i),
\mathbf{v}_i
\Bigr),
\]
where \(d_{\mathcal{X}}\) is a distance on the state space, for example an \(L^1\)-type loss. When trajectory data are available, the objective may also include multi-step rollout losses, but the goal is unchanged: to learn an operator \(\mathcal{G}_{\theta}\) that accurately approximates the discrete flow map while respecting the structure of the underlying conservation law.


\section{Local--Global Neural Operator (LGNO) for Conservation Laws}\label{sec:methods}

\subsection{Conservative One-Step Flow Map Parameterization}

We adopt the residual parameterization of the discrete flow map in Eq.~\eqref{eq::flow_update}, in which \(\Phi_\theta\) gives the cell-centered increment over one step. The network output is thus interpreted neither as the physical flux \(\mathbf{F}(\mathbf{u})\) nor as a classical finite-volume numerical flux. We represent \(\Phi_\theta\) with the LGNO architecture introduced below, or, for the baseline, with an FNO.

When the underlying conservation law preserves the spatial average, as under periodic boundary conditions, we enforce the same property by projecting the learned increment onto the mean-zero subspace. Writing \(\langle\cdot\rangle\) for the spatial average, computed over grid cells in the discrete setting and applied componentwise for systems,
\[
\langle \mathbf v\rangle
=
|\Omega|^{-1}\int_{\Omega}\mathbf v(\mathbf x)\,d\mathbf x,
\quad
\Pi_0\Phi_\theta
=
\Phi_\theta-\langle \Phi_\theta\rangle,
\quad
\mathcal G_\theta(\mathbf u)
=
\mathbf u+\Delta t\,\Pi_0\Phi_\theta(\mathbf u;\Delta t).
\]
\begin{proposition}[Exact global conservation]
With the projected update \(\mathcal G_\theta(\mathbf u)=\mathbf u+\Delta t\,\Pi_0\Phi_\theta(\mathbf u)\) above, for any network parameters \(\theta\), time step \(\Delta t\), and input state \(\mathbf u\),
\[
\langle \mathcal G_\theta(\mathbf u)\rangle
=
\langle \mathbf u\rangle,
\]
so the learned flow map preserves the total of each conserved component exactly.
\end{proposition}
This holds by construction, since \(\langle \Pi_0\Phi_\theta\rangle=0\); it is an exact algebraic property, at both the continuous and discrete levels.

\begin{remark}[Flux-form parameterization]
A natural alternative is to learn a flux-form update as in Eq.~\eqref{eqn:conservative_form}, but we found that this did not significantly improve accuracy. Moreover, when only solution snapshots are available, the numerical flux is not uniquely determined, since the discrete divergence has a nontrivial null space; recovering a unique flux then requires an additional constraint, which complicates the framework without a clear benefit. We therefore adopt the simpler residual increment with the global conservation projection above.
\end{remark}

\subsection{Local--Global Operator Architecture}
\label{subsec:hybrid_backbone}

\begin{figure}[t]
    \centering
    \includegraphics[width=0.98\textwidth]{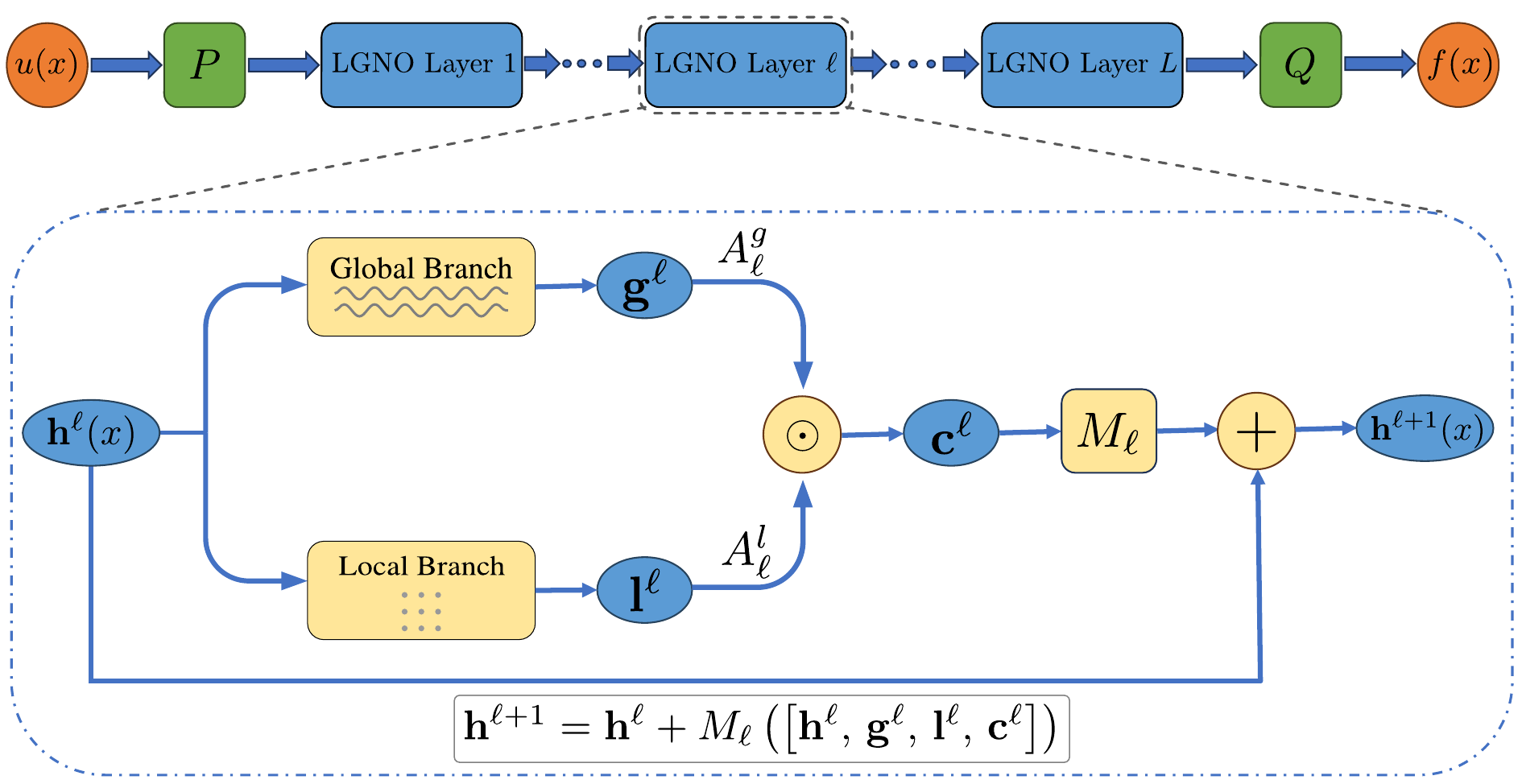}
    \caption{
    \textbf{Schematic of the LGNO architecture.} The input state \(\mathbf{u}^n\) is first lifted to a latent feature field \(\mathbf{h}^1\) by a pointwise map. Each hidden layer has two branches: the global spectral branch applies learnable Fourier multipliers to the low-frequency modes, while the local branch extracts localized features by convolution. Their outputs are coupled multiplicatively and then fused with the incoming hidden state through a pointwise residual mixing map. After \(L\) such layers, a pointwise output network maps the final latent feature \(\mathbf{h}^{L+1}\) to \(\mathbf{f}_{\theta}\), which parameterizes the learned one-step update.
    }
    \label{fig:hybrid_operator}
\end{figure}

We now specify the LGNO architecture that represents the update \(\Phi_\theta\). It combines a global spectral branch, which captures smooth large-scale dynamics, with a local branch, which resolves sharp localized features; its overall structure is summarized in Figure~\ref{fig:hybrid_operator}.

Given a discrete state \(\mathbf{u}^n\), a pointwise lifting map produces an initial latent feature field
\[
\mathbf{h}^{1}= \sigma \left( P_{\theta}(\mathbf{u}^{n}) \right),
\]
where, at each grid point, \(P_{\theta}\) applies the same affine map to the channel vector. Here ``pointwise'' refers to the spatial action, while the channel dimension is lifted from \(d\) to \(d_h\); this is equivalently a convolution with a \(1\times1\) kernel. The nonlinear activation \(\sigma\) then yields the latent feature \(\mathbf{h}^{1}\in\mathbb{R}^{d_h\times N}\). 

Each of the \(L\) hidden layers combines a global spectral branch and a local multiresolution branch, which are coupled multiplicatively and then fused with the incoming state by a residual map. Let \(\mathbf{h}^{\ell}\) denote the hidden feature at layer \(\ell\); for a discretization with \(N\) grid points, \(\mathbf{h}^{\ell}\in\mathbb{R}^{d_h\times N}\), with discrete Fourier transform \(\widehat{\mathbf{h}}^{\ell}\in\mathbb{C}^{d_h\times N}\).

\paragraph{Global branch}
The spectral branch is defined as
\begin{equation}
\label{eq:spectral_branch}
\mathbf{g}^{\ell} = \sigma \left( W_{\ell}^{s}\mathbf{h}^{\ell}  + \mathcal{F}^{-1} \left( R_{\ell}(\mathbf{k})\widehat{\mathbf{h}}^{\ell}(\mathbf{k}) \right) \right).
\end{equation}
Here \(\mathcal{F}\) denotes the discrete Fourier transform. We write \(R_{\ell}(\mathbf{k})\in\mathbb{C}^{d_h\times d_h}\) for the Fourier multiplier at mode \(\mathbf{k}\). For each frequency \(\mathbf{k}\), \(R_{\ell}(\mathbf{k})\) acts on the channel vector \(\widehat{\mathbf{h}}(\mathbf{k})\in\mathbb{C}^{d_h}\) as
\[
\bigl(R_{\ell}(\mathbf{k})\widehat{\mathbf{h}}(\mathbf{k})\bigr)_i
=
\sum_{j=1}^{d_h}
R_{\ell,i,j}(\mathbf{k})\,\widehat{h}_{j}(\mathbf{k}),
\quad i=1,\dots,d_h.
\]
In practice, \(R_{\ell}\) carries learnable entries only on a prescribed band of low-frequency modes and is zero elsewhere, so the spectral branch acts as a trainable low-pass filter. The term \(W_{\ell}^{s}\in\mathbb{R}^{d_h\times d_h}\) is a pointwise linear map on the channel vector at each location. For periodic boundary conditions the spectral convolution is applied directly; for nonperiodic settings, the feature is padded before the transform as a practical boundary treatment.

\paragraph{Local branch}
We propose to include a local branch to incorporate stencil-type locality, a key inductive bias in conventional high-order schemes for conservation laws. This locality is important for shocks, contact discontinuities, and sharp gradients, whose evolution is governed by local wave interactions and may be difficult to capture using a global spectral operator. At layer \(\ell\), the local feature is defined by
\begin{equation}
\label{eq:local_feature}
\mathbf{l}^{\ell}
=
\sigma\left(\mathcal{C}_{\ell}(\mathbf{h}^{\ell})\right),
\quad
\mathbf{l}^{\ell}\in\mathbb{R}^{d_h\times N},
\end{equation}
where \(\mathcal{C}_{\ell}:\mathbb{R}^{d_h\times N}\to\mathbb{R}^{d_h\times N}\) is a local multiresolution operator.

The multiresolution design is inspired by the coarse-to-fine feature-fusion strategy of U-Net~\cite{ronneberger2015u}. Instead of applying local convolutions only on the original grid, we first restrict the hidden feature to a coarser grid, where a fixed-size convolution kernel corresponds to a larger effective neighborhood in the original physical mesh. This enlarges the effective receptive field while keeping the additional computational cost small. The processed coarse feature is then interpolated back to the original resolution and fused with the original feature:
\begin{equation}
\label{eq:local_multiresolution}
\mathcal{C}_{\ell}(\mathbf{h})
=
W_{\ell}^{c}
\left[
\mathbf{h},
\,
\mathcal{I}_{\ell}
\mathcal{K}_{\ell}
\mathcal{P}_{\ell}\mathbf{h}
\right].
\end{equation}
Here \(\mathcal{P}_{\ell}\) denotes average pooling by a factor of two, \(\mathcal{I}_{\ell}\) denotes linear interpolation back to the original resolution, and \(W_{\ell}^{c}\) is a pointwise linear map that fuses the fine- and coarse-resolution features.

In our implementation, \(\mathcal{K}_{\ell}\) is a two-layer convolutional network on the coarse grid:
\begin{equation}
\label{eq:coarse_local_network}
\mathcal{K}_{\ell}(\mathbf{z})
=
\sigma\!\left(
\operatorname{Conv}_{\ell,2}
\left(
\sigma\!\left(
\operatorname{Conv}_{\ell,1}(\mathbf{z})
\right)
\right)
\right),
\quad
\mathbf{z}=\mathcal{P}_{\ell}\mathbf{h}.
\end{equation}
Here \(\operatorname{Conv}_{\ell,1}\) and \(\operatorname{Conv}_{\ell,2}\) denote learned convolutional operations on the coarse grid, each including its convolution kernel and bias term. The convolutions use small kernels with boundary padding chosen according to the prescribed boundary condition, so that the coarse-grid resolution is preserved. Thus, the local branch injects stencil-type information while using the coarse grid to provide a larger effective receptive field at low cost.

\paragraph{LGNO layer}
The global and local features are then coupled through an interaction operator. In the present architecture, we use multiplicative coupling and define
\begin{equation}
\label{eq:multiplicative_coupling}
\mathbf{c}^{\ell} = \left(A_{\ell}^{g}\mathbf{g}^{\ell}\right) \odot \left(A_{\ell}^{l}\mathbf{l}^{\ell}\right),
\quad
\mathbf{c}^{\ell}\in\mathbb{R}^{d_h\times N},
\end{equation}
where \(A_{\ell}^{g},A_{\ell}^{l}\in\mathbb{R}^{d_h\times d_h}\) are pointwise linear maps on the channel vector and \(\odot\) denotes componentwise multiplication. This interaction allows the local response to be modulated by the global spectral feature, while preserving the branchwise separation of the two representations.

The hidden state is updated by
\begin{equation}
\label{eq:hybrid_layer}
\mathbf{h}^{\ell+1} = \mathcal{L}_{\ell}(\mathbf{h}^{\ell}) := \mathbf{h}^{\ell}
+ M_{\ell} 
\left(
\left[ \mathbf{h}^{\ell}, \mathbf{g}^{\ell}, \mathbf{l}^{\ell}, \mathbf{c}^{\ell} \right]
\right),
\quad
\ell=1,\dots,L,
\end{equation}
where \(\mathbf{h}^{\ell+1}\in\mathbb{R}^{d_h\times N}\), \([\mathbf{h}^{\ell},\mathbf{g}^{\ell},\mathbf{l}^{\ell},\mathbf{c}^{\ell}]\in\mathbb{R}^{4d_h\times N}\), and \(M_{\ell}\in\mathbb{R}^{d_h\times 4d_h}\) is a pointwise linear map on the concatenated channel vector. This residual update lets each layer refine the latent representation while retaining information from the incoming state.

Finally, after \(L\) layers, a pointwise output map \(Q_{\theta}\) produces the field required by the update rule:
\begin{equation}
\label{eq:hybrid_head}
\mathbf{f}_{\theta}
=
Q_{\theta}(\mathbf{h}^{L+1}).
\end{equation}
Here \(Q_{\theta}:\mathbb{R}^{d_h}\to\mathbb{R}^{d}\) is applied at each grid point, implemented as a two-layer pointwise network with a nonlinear activation in between. The output \(\mathbf{f}_{\theta}\) parameterizes the learned increment in the discrete flow map.

A concrete implementation of one LGNO layer \(\mathcal{L}_{\ell}\) is summarized in Algorithm~\ref{alg:hybrid_layer}.
\begin{algorithm}[H]
\caption{Implementation of one LGNO layer \(\mathcal{L}_{\ell}\)}
\label{alg:hybrid_layer}
\begin{algorithmic}[1]
\Require Hidden feature field \(\mathbf{h}^{\ell}\in\mathbb{R}^{d_h\times N}\); Fourier multiplier \(R_{\ell}\); pointwise linear maps \(W_{\ell}^{s},W_{\ell}^{c},A_{\ell}^{g},A_{\ell}^{l},M_{\ell}\); average-pooling operator \(\mathcal{P}_{\ell}\); coarse-grid convolutional network \(\mathcal{K}_{\ell}\); interpolation operator \(\mathcal{I}_{\ell}\)
\Ensure Updated hidden feature field \(\mathbf{h}^{\ell+1}\)

\Statex \textbf{Global spectral branch}
\State Compute the global feature
\[
\mathbf{g}^{\ell}
\gets
\sigma\!\left(
W_{\ell}^{s}\mathbf{h}^{\ell} + 
\mathcal{F}^{-1}\!\left(
R_{\ell}(\mathbf{k})\,\widehat{\mathbf{h}}^{\ell}(\mathbf{k})
\right)
\right)
\]

\Statex \textbf{Local multiresolution branch}
\State Compute the local feature
\[
\mathbf{l}^{\ell}
\gets
\sigma\!\left(
W_{\ell}^{c}
\left[
\mathbf{h}^{\ell},
\,
\mathcal{I}_{\ell}
\mathcal{K}_{\ell}
\mathcal{P}_{\ell}\mathbf{h}^{\ell}
\right]
\right)
\]

\Statex \textbf{Global--local interaction and residual update}
\State Compute the coupling
\[
\mathbf{c}^{\ell}
\gets
\left(A_{\ell}^{g}\mathbf{g}^{\ell}\right)
\odot
\left(A_{\ell}^{l}\mathbf{l}^{\ell}\right)
\]
\State Update the hidden state
\[
\mathbf{h}^{\ell+1}
\gets
\mathbf{h}^{\ell}
+
M_{\ell}
\left(
\left[
\mathbf{h}^{\ell},
\mathbf{g}^{\ell},
\mathbf{l}^{\ell},
\mathbf{c}^{\ell}
\right]
\right)
\]
\State \Return \(\mathbf{h}^{\ell+1}\)
\end{algorithmic}
\end{algorithm}

\subsection{Boundary Padding for Outflow}\label{sec:outflow_padding}

For two-dimensional outflow problems, waves and coherent structures may reach the computational boundary during autoregressive rollout. In this regime, simple boundary copying or fixed extrapolation can provide inaccurate ghost values and introduce artifacts into local convolutional stencils. We therefore use a learned ghost-cell padding operator for the local convolutional operations under outflow boundary conditions. This padding is applied to latent feature maps in the local branch. 

Let \(\mathbf{x}\in\mathbb{R}^{B\times d_h\times H\times W}\) be a latent feature map before a local convolution, where \(B\) is the batch size, \(d_h\) is the number of channels, and \(H\times W\) is the spatial resolution. For an odd convolution kernel of size \(k\), we set the ghost width to
\[
g_w=\frac{k-1}{2}.
\]
The padding operator constructs ghost features from boundary-adjacent interior information before applying the convolution.

Given a context width \(c\), we first predict the left and right ghost columns from boundary-adjacent interior strips:
\[
\mathbf{g}_{L}
=
\mathcal{G}_{L}(\mathbf{x}_{:, :, :, 1:c}),
\qquad
\mathbf{g}_{R}
=
\mathcal{G}_{R}(\mathbf{x}_{:, :, :, W-c+1:W}),
\]
where \(\mathcal{G}_{L}\) and \(\mathcal{G}_{R}\) are shallow local convolutional networks whose outputs satisfy
\[
\mathbf{g}_{L},\mathbf{g}_{R}
\in
\mathbb{R}^{B\times d_h\times H\times g_w}.
\]
The horizontally padded feature map is then
\[
\mathbf{x}_{LR}
=
\left[
\mathbf{g}_{L},\,
\mathbf{x},\,
\mathbf{g}_{R}
\right]
\in
\mathbb{R}^{B\times d_h\times H\times (W+2g_w)}.
\]

We next predict the top and bottom ghost rows from boundary-adjacent strips of \(\mathbf{x}_{LR}\):
\[
\mathbf{g}_{T}
=
\mathcal{G}_{T}((\mathbf{x}_{LR})_{:, :, 1:c, :}),
\qquad
\mathbf{g}_{B}
=
\mathcal{G}_{B}((\mathbf{x}_{LR})_{:, :, H-c+1:H, :}),
\]
where \(\mathcal{G}_{T}\) and \(\mathcal{G}_{B}\) are defined analogously and output
\[
\mathbf{g}_{T},\mathbf{g}_{B}
\in
\mathbb{R}^{B\times d_h\times g_w\times (W+2g_w)}.
\]
Thus the corner ghost features are generated as part of the top and bottom ghost rows. The final padded feature map is
\[
\operatorname{Pad}_{\rm out}(\mathbf{x})
=
\left[
\mathbf{g}_{T};\,
\mathbf{x}_{LR};\,
\mathbf{g}_{B}
\right]
\in
\mathbb{R}^{B\times d_h\times (H+2g_w)\times (W+2g_w)}.
\]
Here \([\cdot,\cdot]\) denotes horizontal concatenation along the width dimension, and \([\cdot;\cdot]\) denotes vertical concatenation along the height dimension.

The operator \(\operatorname{Pad}_{\rm out}\) is applied before each local convolution in the local branch. The subsequent convolution is performed with stride one and without additional zero padding, so that the output has the same spatial resolution as the input feature map. This learned padding provides a lightweight treatment of outflow boundaries for the neural local stencil, allowing the local branch to maintain a complete receptive field near the boundary without injecting zero-padding artifacts.

\subsection{Training Loss}
\label{subsec:training_objective}

The architecture above induces a discrete flow map \(\mathcal{G}_{\theta}\) for the one-step update, which we train in a supervised manner on snapshot pairs generated by a high-order numerical scheme. Our training objective combines a physical-space prediction loss with a high-frequency spectral penalty that suppresses oscillatory error near discontinuities.

Let \(\mathbf{u}^{n+1}\) and \(\widehat{\mathbf{u}}^{\,n+1}_{\theta}\) denote the reference and predicted next states on the discrete grid \(\Omega_h\), and let \(d\) be the number of physical components. The physical-space loss is
\begin{equation}
\label{eq:phys_loss}
\mathcal{L}_{\mathrm{phys}}
=
\begin{cases}
\displaystyle
\frac{1}{|\Omega_h|}
\left\|
\widehat{\mathbf{u}}^{\,n+1}_{\theta}
-
\mathbf{u}^{n+1}
\right\|_{L^1},
& d=1, \\[1.2em]
\displaystyle
\frac{1}{d}
\sum_{m=1}^{d}
\frac{
\left\|
\widehat{\mathbf{u}}^{\,n+1}_{\theta,m}
-
\mathbf{u}^{n+1}_{m}
\right\|_{L^1}
}{
\left\|
\mathbf{u}^{n+1}_{m}
\right\|_{L^1}
+\varepsilon
},
& d>1,
\end{cases}
\end{equation}
where \(\varepsilon>0\) is a small constant for numerical stability. Thus, scalar problems use the standard mean absolute error, while multi-component systems use a componentwise relative \(L^1\) loss to balance variables with different physical scales.

To suppress high-frequency oscillatory error, let \(\mathbf{e}=\widehat{\mathbf{u}}^{\,n+1}_{\theta}-\mathbf{u}^{n+1}\) be the one-step prediction error and \(\widehat{\mathbf{e}}(k)=\mathcal{F}(\mathbf{e})(k)\) its discrete Fourier coefficient at mode \(k\). With \(K\) the number of retained Fourier modes and \(k_{0}=\lfloor \kappa K \rfloor\) for a prescribed fraction \(\kappa\in(0,1)\), we define
\begin{equation}
\label{eq:high_freq_loss}
\mathcal{L}_{\mathrm{hf}}
=
\frac{1}{K-k_{0}}
\sum_{k=k_{0}}^{K-1}
\left\|
\widehat{\mathbf{e}}(k)
\right\|_{2}^{2}.
\end{equation}
Here \(\|\widehat{\mathbf{e}}(k)\|_2^2=\sum_{m=1}^{d}|\widehat{e}_m(k)|^2\) for multi-component systems, and reduces to \(|\widehat{e}(k)|^2\) when \(d=1\). The total training loss is
\begin{equation}
\label{eq:total_training_loss}
\mathcal{L}
=
\mathcal{L}_{\mathrm{phys}}
+
\lambda_{\mathrm{hf}}\mathcal{L}_{\mathrm{hf}},
\end{equation}
where \(\lambda_{\mathrm{hf}}>0\) controls the strength of the spectral penalty.

\begin{remark}[Symmetry-based extension of the state range]
Since the learned flow map is applied at inference time, exact scaling symmetries of the governing equations can be used as preprocessing and postprocessing transformations. Such symmetries are present in several homogeneous hyperbolic conservation laws, including inviscid Burgers, the shallow-water equations, and the compressible Euler equations, although the precise scaling weights are equation-dependent. When such a symmetry is available, a test state whose magnitude lies outside the training range can be mapped into a comparable range, advanced by the learned operator in the rescaled variables, and then mapped back, preserving the corresponding dynamics. We use this idea in the two-dimensional Euler shock tests in~\ref{app:euler2d_shock_tests}, where the detailed transformation is given.
\end{remark}


\section{Numerical Results}\label{sec:results}

In this section, we evaluate the proposed LGNO on several benchmark problems: linear advection, Burgers, shallow water, and compressible Euler equations. All problems use periodic boundary conditions except the final two-dimensional Euler case, which considers both periodic and outflow boundary conditions. Additional two-dimensional Euler shock tests, based on standard two-dimensional Riemann problems not included in the training data, are reported in~\ref{app:euler2d_shock_tests} as a further robustness check in strongly discontinuous regimes.

We compare LGNO against a parameter-matched FNO baseline, using comparable parameter counts on every benchmark except the two-dimensional Burgers case, on which the FNO baseline is substantially larger. At the single-step level, LGNO already reduces the mean relative \(L^1\) error by factors of \(2.18\) to \(5.09\) (Table~\ref{tab:onestep_summary}). These per-step gains compound under autoregressive rollout, so that the long-time advantage of LGNO is far larger than the one-step errors alone would suggest. LGNO remains stable and accurate over long horizons, whereas the FNO baseline degrades steadily and, in some of the two-dimensional Euler cases, fails entirely by producing unphysical solutions. We examine this long-horizon behavior and the solution quality for each benchmark in turn. 

A wall-clock cost comparison against WENO-Z, reported in~\ref{apx:computational_cost}, shows that the learned models substantially reduce inference cost relative to the WENO-Z scheme. Model, training, and benchmark configurations are detailed in~\ref{apx:implementation_details} and~\ref{apx:numerical_setup}.

\begin{table}[tbp]
\centering
\caption{One-step test relative \(L^1\) errors for the proposed LGNO model and the FNO baseline, reported as mean \(\pm\) standard deviation over the test set. The last column gives the mean-error ratio \(\mathrm{FNO}/\mathrm{LGNO}\), where values larger than \(1\) indicate better performance of LGNO. LGNO consistently outperforms the FNO baseline across all reported benchmarks, reducing the mean error by factors between \(2.18\) and \(5.09\).}
\label{tab:onestep_summary}
\begin{tabular}{lccc}
\toprule
Problem & LGNO & FNO & \(\text{FNO}/\text{LGNO}\) \\
\midrule
Linear advection
& \(4.57\times 10^{-4}\pm 1.16\times 10^{-3}\)
& \(1.99\times 10^{-3}\pm 2.91\times 10^{-3}\)
& \(4.36\) \\
Burgers 1D
& \(2.38\times 10^{-4}\pm 3.94\times 10^{-4}\)
& \(5.19\times 10^{-4}\pm 5.96\times 10^{-4}\)
& \(2.18\) \\
SWE 1D
& \(5.27\times 10^{-4}\pm 5.41\times 10^{-4}\)
& \(1.16\times 10^{-3}\pm 1.23\times 10^{-3}\)
& \(2.20\) \\
Euler 1D
& \(1.23\times 10^{-3}\pm 1.30\times 10^{-3}\)
& \(3.85\times 10^{-3}\pm 3.74\times 10^{-3}\)
& \(3.13\) \\
Burgers 2D
& \(3.53\times 10^{-4}\pm 2.69\times 10^{-4}\) 
& \(9.81\times 10^{-4}\pm 8.63\times 10^{-4}\) 
& \(2.78\) \\
Euler 2D periodic
& \(3.87\times 10^{-3}\pm 2.31\times 10^{-3}\)
& \(1.97\times 10^{-2}\pm 8.86\times 10^{-3}\)
& \(5.09\) \\
Euler 2D outflow
& \(1.67\times 10^{-3}\pm 1.07\times 10^{-3}\)
& \(7.33\times 10^{-3}\pm 3.08\times 10^{-3}\)
& \(4.39\) \\
\bottomrule
\end{tabular}
\end{table}


\subsection{Linear Advection Equation}
\label{subsec:results_convection}

We begin with the linear advection equation
\begin{equation}
u_t + u_x = 0, \quad x\in[0,1].
\label{eq:convection_pde}
\end{equation}
Although the dynamics are simple, accurate long-time prediction remains nontrivial for discontinuous data: the solution is transported without intrinsic smoothing, so any representation error in the profile propagates directly over time. For this reason it is a standard benchmark for assessing both accuracy and numerical dissipation in scheme design, and here it provides a clean test of whether LGNO improves on a purely spectral baseline when the target is poorly represented by a truncated Fourier series.

\begin{figure}[tbp]
    \centering
    \begin{subfigure}[b]{0.49\textwidth}
        \centering
        \includegraphics[width=\textwidth]{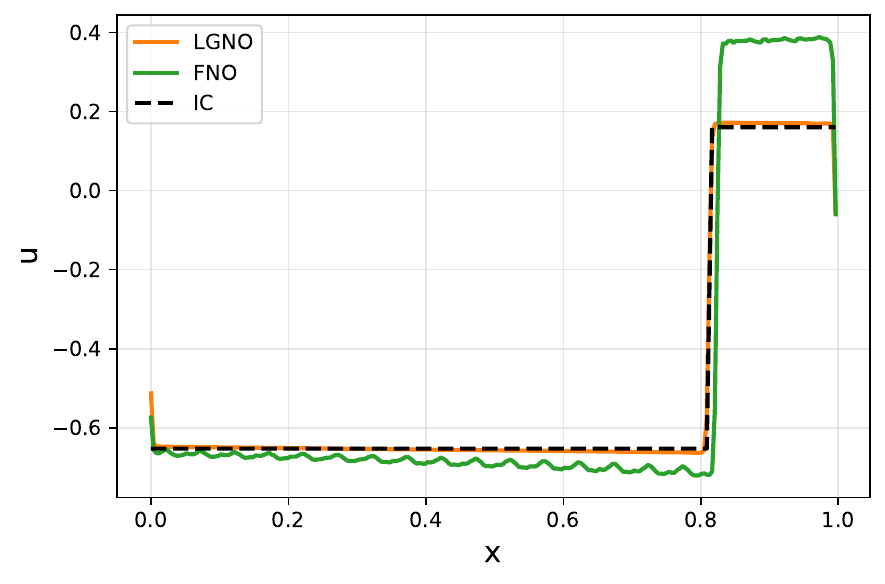}
        \caption{Solution at \(T=3\)}
    \end{subfigure}
    \hfill
    \begin{subfigure}[b]{0.49\textwidth}
        \centering
        \includegraphics[width=\textwidth]{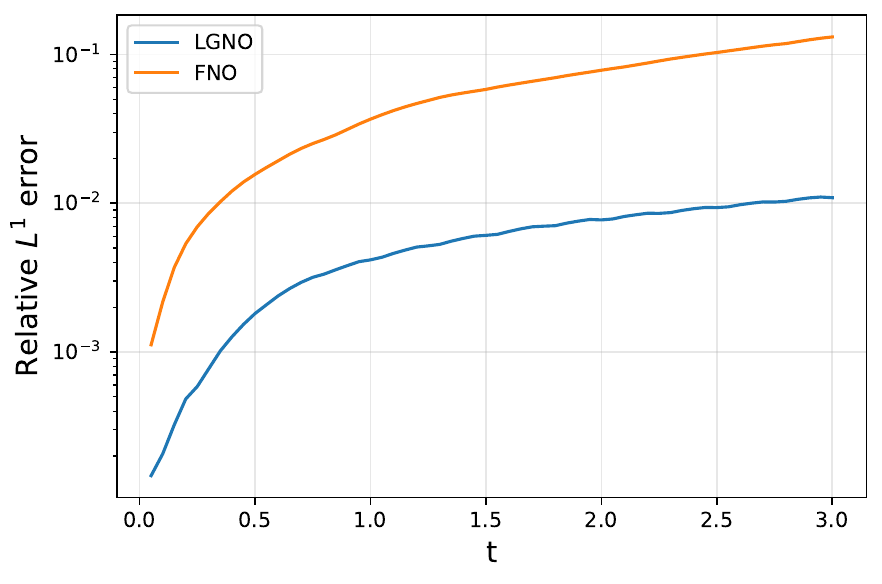}
        \caption{Relative \(L^1\) error versus time}
    \end{subfigure}
    \caption{\textbf{Long-horizon rollout performance for linear advection.} Panel (a) compares the predictions at \(T=3\) with the reference solution, which coincides with the initial condition due to periodic unit-speed transport. Panel (b) shows the rollout relative \(L^1\) error over time for the same representative discontinuous initial condition.}
    \label{fig:convection_rollout}
\end{figure}

We evaluate the LGNO and FNO models under autoregressive rollout of a representative discontinuous initial condition up to \(T=3\). With the learned time step \(\Delta t=0.05\), this amounts to a long rollout of \(60\) steps. Since \(T=3\) corresponds to three full periods of the unit-speed equation, the reference solution coincides with the initial condition. As shown in Figure~\ref{fig:convection_rollout}, although both models use a conservative update, the FNO baseline accumulates noticeable amplitude drift and fails to maintain the flat states, whereas LGNO stays closely aligned with the reference. This drift is not inconsistent with mass conservation, since the conservative update preserves the total mass but not the profile shape. The FNO jump is also noticeably displaced from its exact location, whereas LGNO matches it well. The relative \(L^1\) error curve confirms this, remaining consistently below the FNO baseline throughout the rollout.

\begin{figure}[tbp]
    \centering
    \begin{subfigure}[b]{0.31\textwidth}
        \centering
        \includegraphics[width=\textwidth]{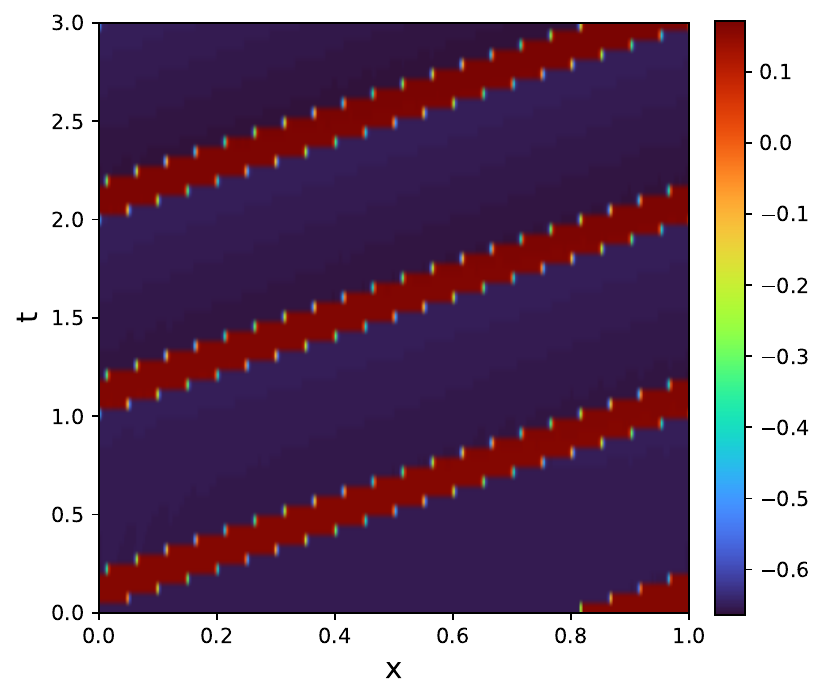}
        \caption{LGNO solution in the \(x\)-\(t\) plane}
    \end{subfigure}
    \hfill
    \begin{subfigure}[b]{0.31\textwidth}
        \centering
        \includegraphics[width=\textwidth]{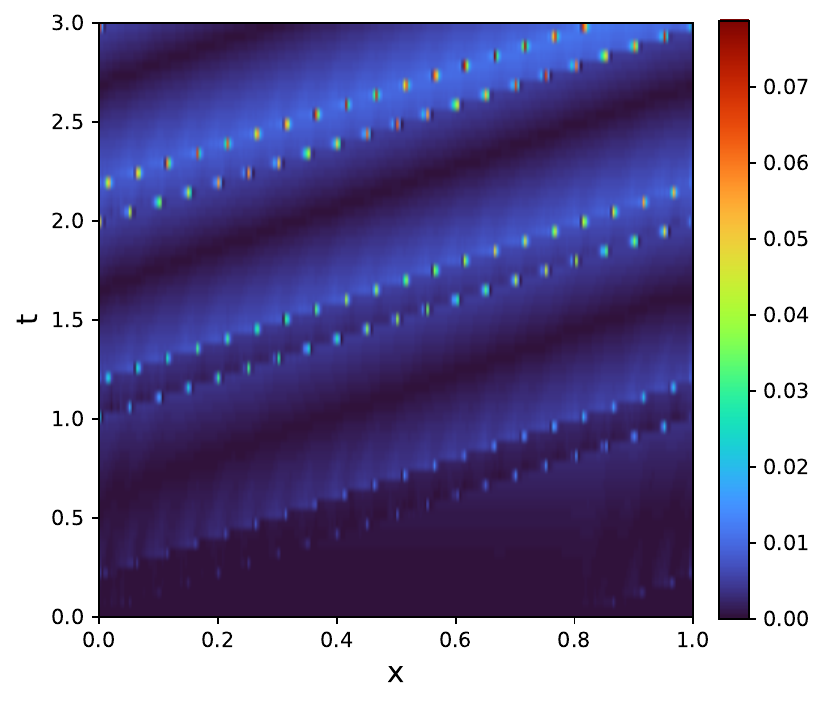}
        \caption{LGNO absolute error in the \(x\)-\(t\) plane}
    \end{subfigure}
    \hfill
    \begin{subfigure}[b]{0.31\textwidth}
        \centering
        \includegraphics[width=\textwidth]{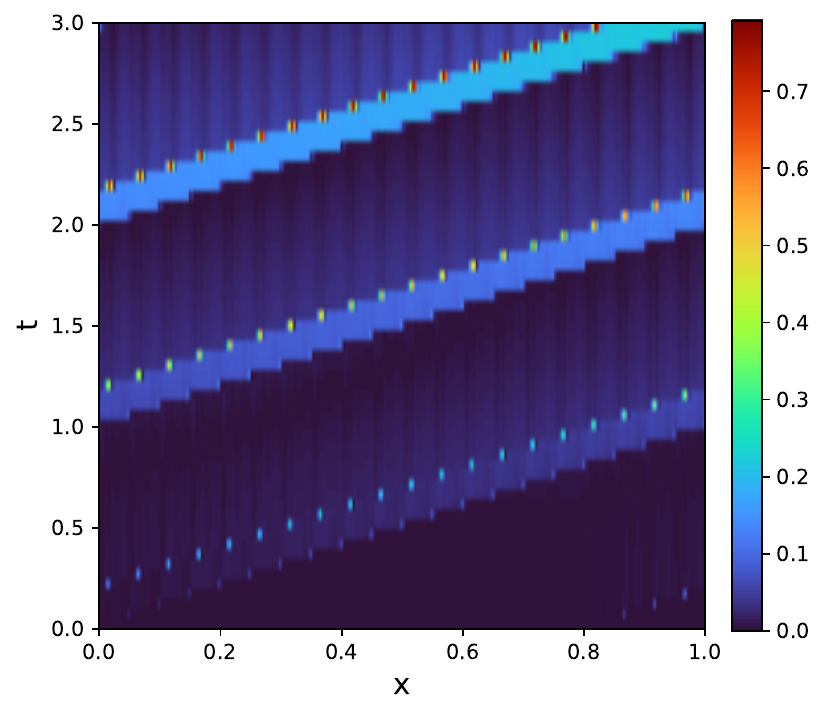}
        \caption{FNO absolute error in the \(x\)-\(t\) plane}
    \end{subfigure}
    \caption{\textbf{Spatiotemporal solution and error for linear advection.} Panel (a) shows the LGNO solution in the \(x\)-\(t\) plane. Panels (b) and (c) compare the spatiotemporal absolute errors of LGNO and the FNO baseline. LGNO preserves the translation pattern better and exhibits significantly smaller errors over the full rollout horizon.}
    \label{fig:convection_xt}
\end{figure}

The spatiotemporal plots in Figure~\ref{fig:convection_xt} reinforce this over the full rollout: the LGNO \(x\)-\(t\) solution preserves the translation pattern with little visible distortion and its error map stays uniformly small, whereas the FNO baseline shows much larger error accumulation.

\begin{figure}[tbp]
    \centering
    \begin{subfigure}[b]{0.49\textwidth}
        \centering
        \includegraphics[width=\textwidth]{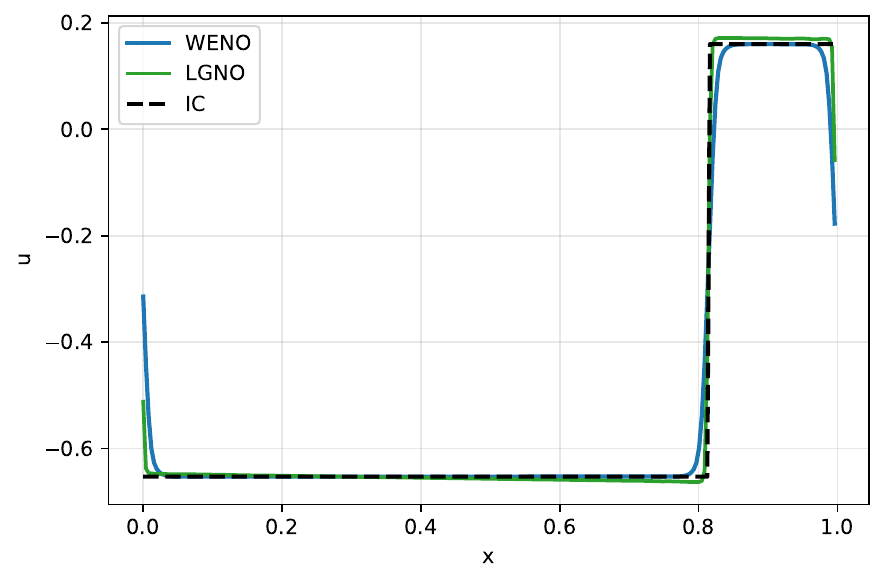}
        \caption{\(T=3\)}
    \end{subfigure}
    \hfill
    \begin{subfigure}[b]{0.49\textwidth}
        \centering
        \includegraphics[width=\textwidth]{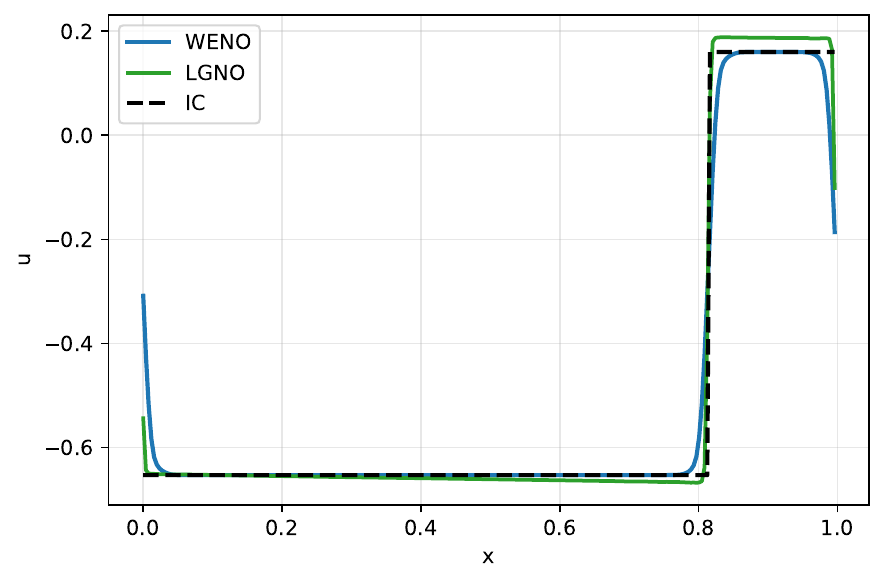}
        \caption{\(T=6\)}
    \end{subfigure}
    \caption{\textbf{Comparison of LGNO and WENO-Z in the linear advection problem.} All methods are evaluated on the same \(256\)-cell mesh and compared at two rollout times. LGNO preserves sharper transitions than the same-mesh WENO-Z solution, suggesting less diffusive behavior in this example, though over long times it shows a larger amplitude error than WENO-Z.}
    \label{fig:pureconvection_same_mesh_dissipation}
\end{figure}

Figure~\ref{fig:pureconvection_same_mesh_dissipation} compares LGNO and WENO-Z on the same \(256\)-cell mesh at two rollout times, \(T=3\) and \(T=6\), to contrast their dissipation behavior. For this representative discontinuous profile, LGNO is less diffusive than the same-mesh WENO-Z solution. At \(T=6\) (a \(120\)-step rollout), however, LGNO exhibits a small amplitude drift that could be further controlled with a conservative bound-preserving correction if needed.

\subsection{One-Dimensional Burgers Equation}
\label{subsec:results_burgers1d}

We next consider the one-dimensional inviscid Burgers equation
\begin{equation}
u_t + \left(\frac{u^2}{2}\right)_x = 0, \quad x\in[0,1].
\label{eq:burgers1d_pde}
\end{equation}
Nonlinear transport and progressive steepening make this a standard test of whether the learned model maintains accuracy as the solution develops increasingly localized structures.

\begin{figure}[tbp]
    \centering
    \begin{subfigure}[b]{0.49\textwidth}
        \centering
        \includegraphics[width=\textwidth]{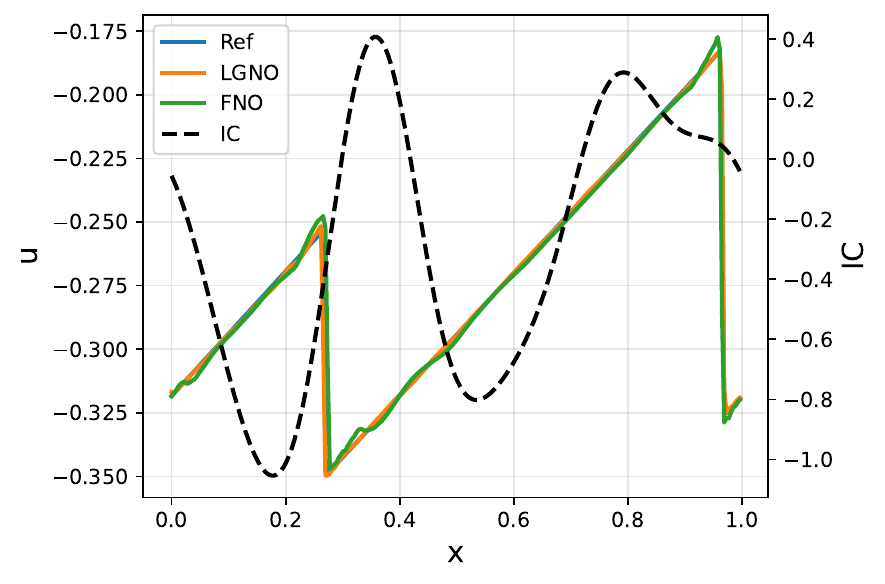}
        \caption{Solution at \(T=4\)}
    \end{subfigure}
    \hfill
    \begin{subfigure}[b]{0.49\textwidth}
        \centering
        \includegraphics[width=\textwidth]{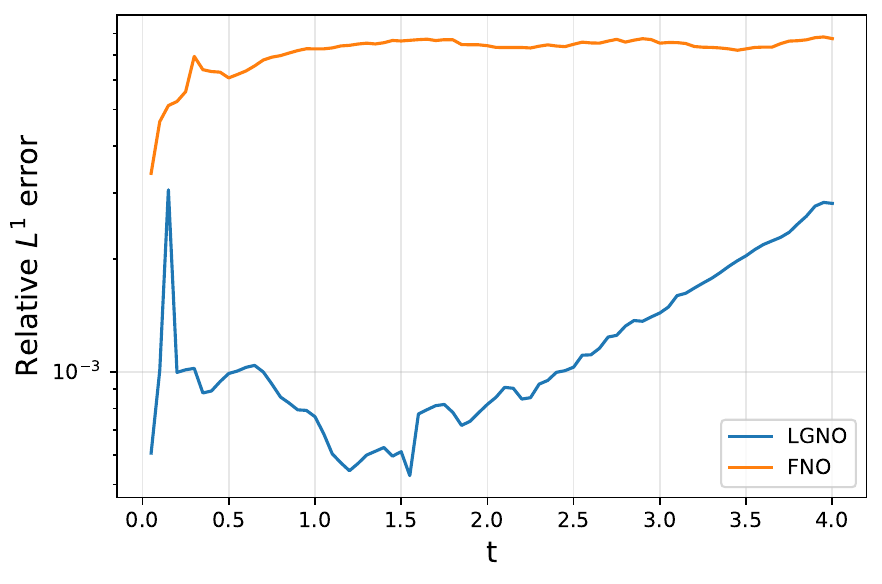}
        \caption{Relative \(L^1\) error versus time}
    \end{subfigure}
    \caption{\textbf{Long-horizon rollout accuracy for the one-dimensional Burgers equation.} Panel (a) compares the final-time solution at \(T=4\), with the randomly selected initial condition overlaid as a dashed line on a separate right axis. Panel (b) shows the rollout relative \(L^1\) error over time. LGNO preserves the shock and ramp structures more accurately and maintains a substantially lower rollout error than the FNO baseline, which exhibits larger local deviations and spurious oscillations near steep fronts.}
    \label{fig:burgers1d_rollout_curves}
\end{figure}

We evaluate both models under autoregressive rollout up to \(T=4\) for a representative, randomly selected initial condition. In the left panel of Figure~\ref{fig:burgers1d_rollout_curves}, this initial condition is overlaid as a dashed line on a separate right axis; since its amplitude is much larger than that of the steepened solution at \(T=4\), plotting it on its own scale keeps the final-time comparison among the reference, LGNO, and FNO more apparent. LGNO stays close to the reference and preserves the main shock and ramp structures, whereas the FNO baseline develops larger local deviations and spurious oscillations near steep fronts. The rollout relative \(L^1\) error of LGNO also remains consistently and substantially lower throughout.

\begin{figure}[tbp]
    \centering
    \begin{subfigure}[b]{0.32\textwidth}
        \centering
        \includegraphics[width=\textwidth]{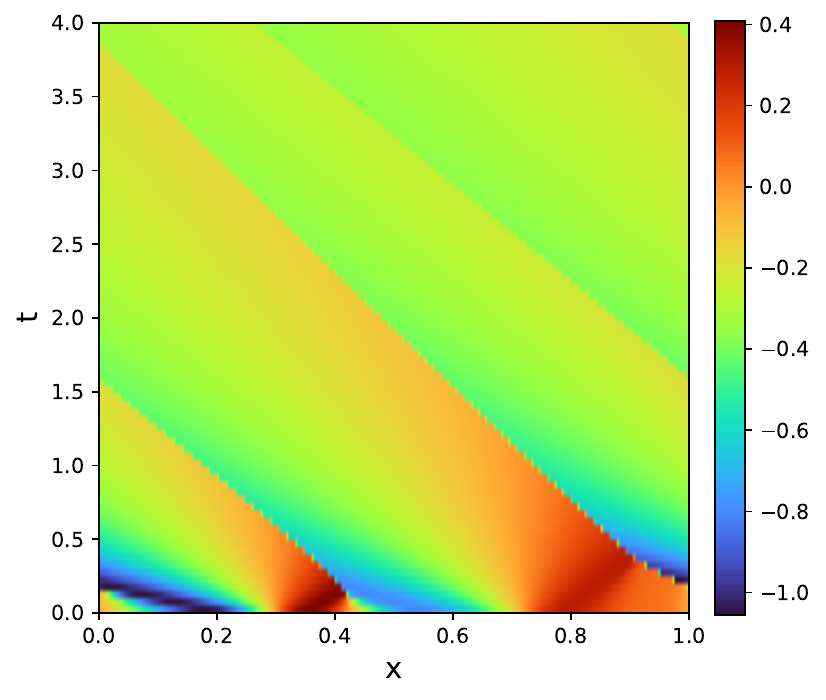}
        \caption{LGNO solution in the \(x\)-\(t\) plane}
    \end{subfigure}
    \begin{subfigure}[b]{0.32\textwidth}
        \centering
        \includegraphics[width=\textwidth]{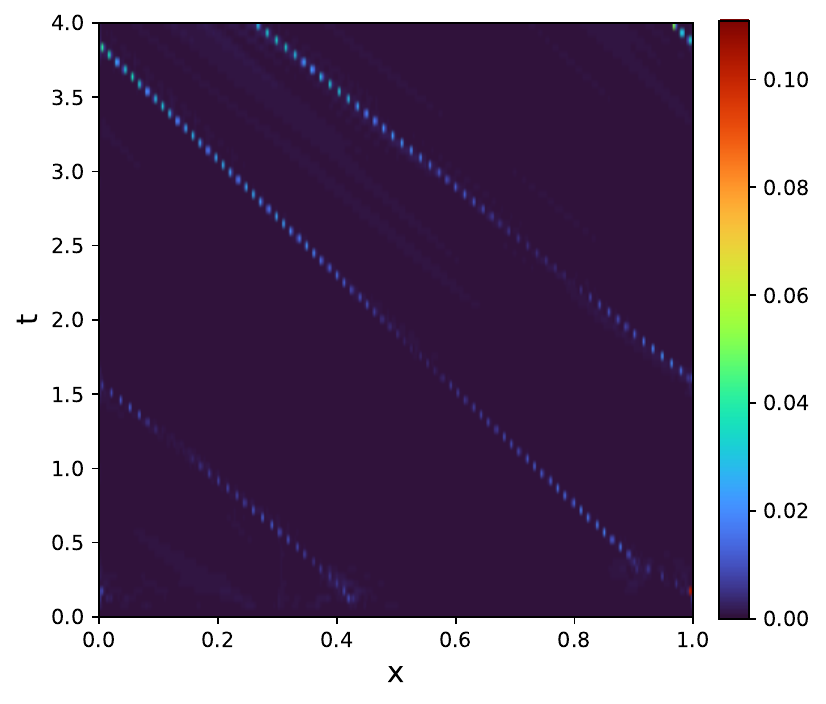}
        \caption{LGNO absolute error in the \(x\)-\(t\) plane}
    \end{subfigure}
    \begin{subfigure}[b]{0.32\textwidth}
        \centering
        \includegraphics[width=\textwidth]{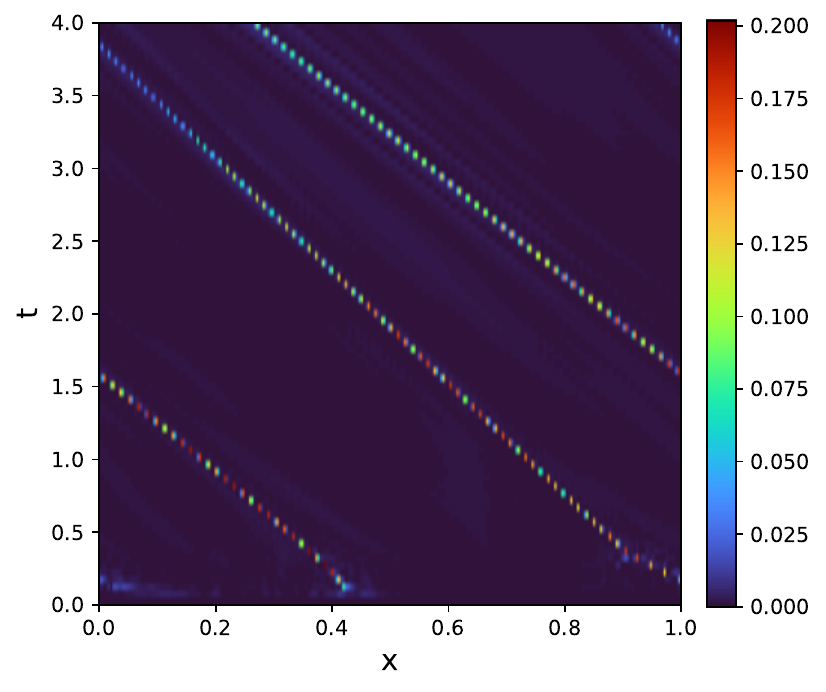}
        \caption{FNO absolute error in the \(x\)-\(t\) plane}
    \end{subfigure}
    \caption{\textbf{Spatiotemporal rollout behavior for the one-dimensional Burgers equation.} Panel (a) shows the LGNO rollout in the \(x\)-\(t\) plane. Panels (b) and (c) compare the spatiotemporal absolute errors of LGNO and the FNO baseline, respectively. Errors are mainly concentrated near the moving shocks in both models, but LGNO yields smaller error bands than the FNO baseline.}
    \label{fig:burgers1d_xt_error}
\end{figure}

Figure~\ref{fig:burgers1d_xt_error} examines the rollout in the \(x\)-\(t\) plane. LGNO preserves the moving shock and ramp structures over the full interval, and although the dominant errors of both models localize near shocks, LGNO produces smaller error bands than the FNO baseline. Overall, LGNO is more effective for long-time prediction of nonlinear wave dynamics.

\subsection{One-Dimensional Shallow Water Equations}
\label{subsec:results_swe1d}

The shallow water system is given by
\begin{equation}
\mathbf{U}_t+\mathbf{F}(\mathbf{U})_x=0,
\quad
\mathbf{U}=(h,hu)^{\top},\quad
\mathbf{F}(\mathbf{U})=\left(hu,\,hu^2+\frac{1}{2}gh^2\right)^{\top},
\quad x\in[0,1],\quad g=9.80665,
\label{eq:swe_pde}
\end{equation}
where \(h\) is the fluid height and \(u\) the velocity. Compared with the scalar Burgers equation, this system has multiple coupled fields and richer wave interactions, testing whether LGNO remains effective for a multi-component hyperbolic system.

\begin{figure}[tbp]
    \centering
    \begin{subfigure}[b]{0.49\textwidth}
        \centering
        \includegraphics[width=\textwidth]{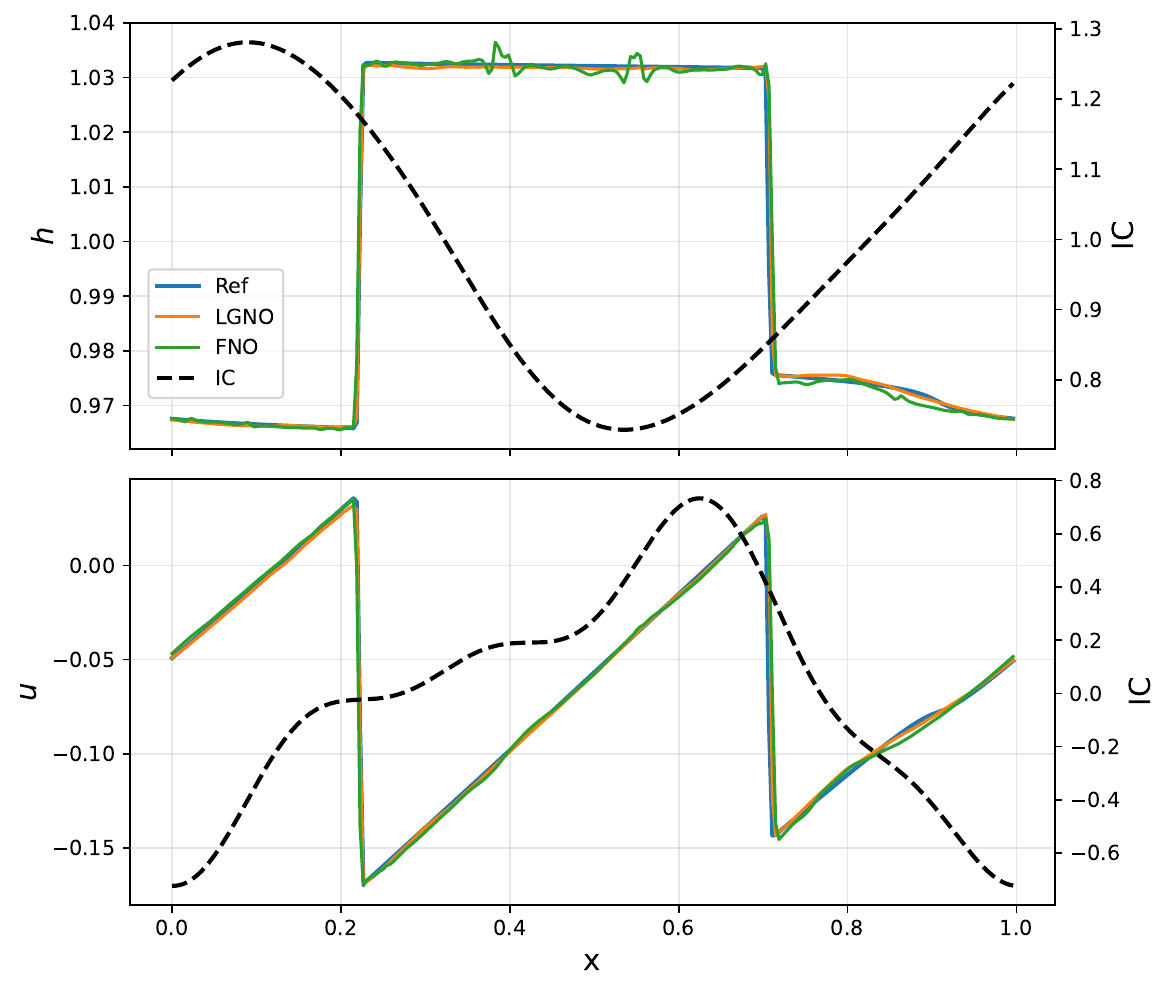}
        \caption{Representative test example 1 at \(T=3\)}
    \end{subfigure}
    \hfill
    \begin{subfigure}[b]{0.49\textwidth}
        \centering
        \includegraphics[width=\textwidth]{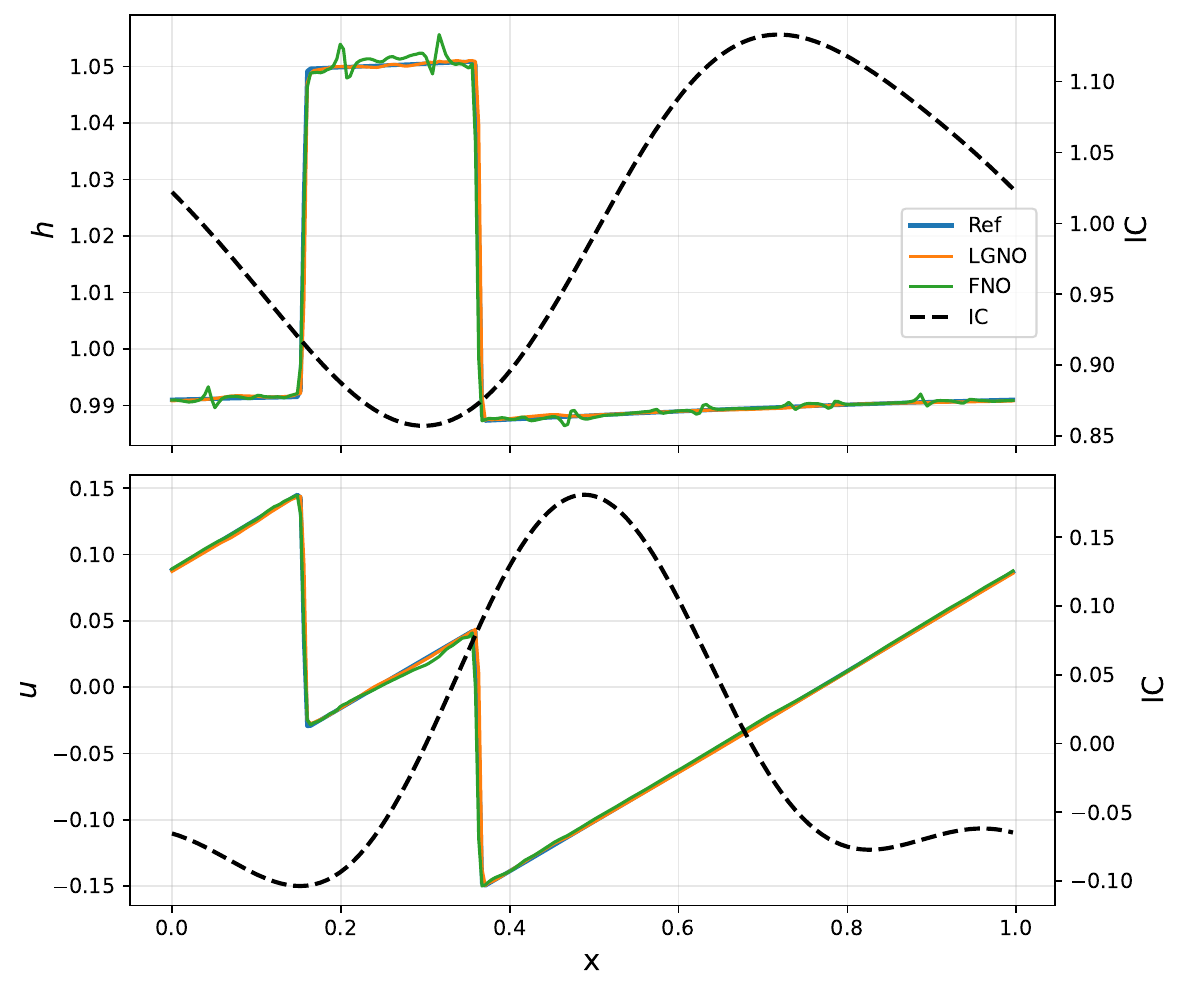}
        \caption{Representative test example 2 at \(T=3\)}
    \end{subfigure}
    \caption{\textbf{Final-time solution profiles for the one-dimensional shallow water equations.} Panels (a) and (b) show two representative test examples at \(T=3\). Each panel reports the fluid height \(h\) on top and the velocity \(u\) on the bottom, with the initial condition shown as a dashed line on a separate right axis. LGNO remains close to the reference solution in both variables, while the FNO baseline shows larger deviations, including unphysical oscillations in \(h\) and visible deviations in \(u\).}
    \label{fig:swe1d_examples}
\end{figure}

Figure~\ref{fig:swe1d_examples} shows predicted solutions at \(T=3\) for two representative test examples; the prediction is performed on the conserved variables, from which the physical variables are recovered at the final time. In each panel, the randomly selected initial condition is overlaid as a dashed line on a separate right axis for both the height and velocity subplots; because it spans a much larger range than the final-time solution, this separate scaling keeps the comparison among the reference, LGNO, and FNO clear.

In both examples, LGNO accurately captures the wave structures in the height and velocity fields, remaining close to the reference in smooth regions and near sharp transitions. The FNO baseline captures the coarse propagation pattern but produces unphysical oscillations in \(h\), including in smooth regions away from any discontinuity, together with visible deviations in \(u\). Its error therefore spreads into the smooth parts of the solution rather than staying localized at the shocks and contact discontinuities. LGNO thus yields a more stable and accurate prediction of the coupled shallow water dynamics.

\subsection{One-Dimensional Euler Equations}
\label{subsec:results_euler1d}

We then consider the one-dimensional compressible Euler equations for an ideal gas,
\begin{equation}
\mathbf{U}_t + \mathbf{F}(\mathbf{U})_x =0,\quad x\in[0,1],
\quad
\mathbf{U}=(\rho,m,E)^{\top},
\label{eq:euler1d_pde}
\end{equation}
where
\[
\mathbf{F}(\mathbf{U})=
\left(m,\frac{m^2}{\rho}+p,(E+p)u\right)^{\top},
\quad m=\rho u,\quad
p=(\gamma-1)\left(E-\frac{1}{2}\rho u^2\right).
\]
This is the most demanding one-dimensional example: its coupled conserved variables, nonlinear wave interactions, and sharp features make it a critical test of whether LGNO improves over the FNO baseline.

\begin{figure}[tbp]
    \centering
    \begin{subfigure}[b]{0.49\textwidth}
        \centering
        \includegraphics[width=\textwidth]{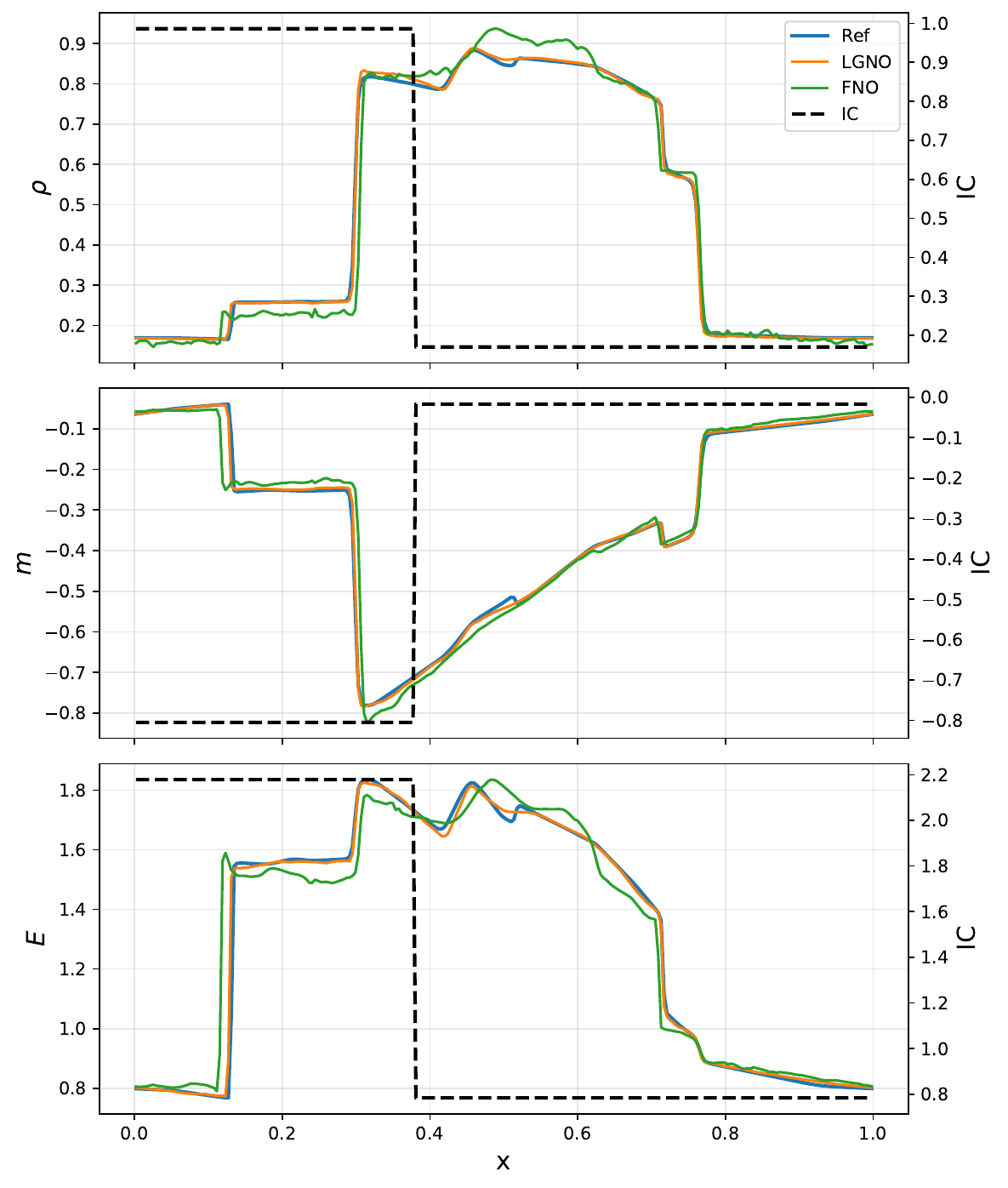}
        \caption{Final-time solution in conserved variables \((\rho,m,E)\)}
    \end{subfigure}
    \hfill
    \begin{subfigure}[b]{0.49\textwidth}
        \centering
        \includegraphics[width=\textwidth,trim={0 0 0 213.6pt},clip]{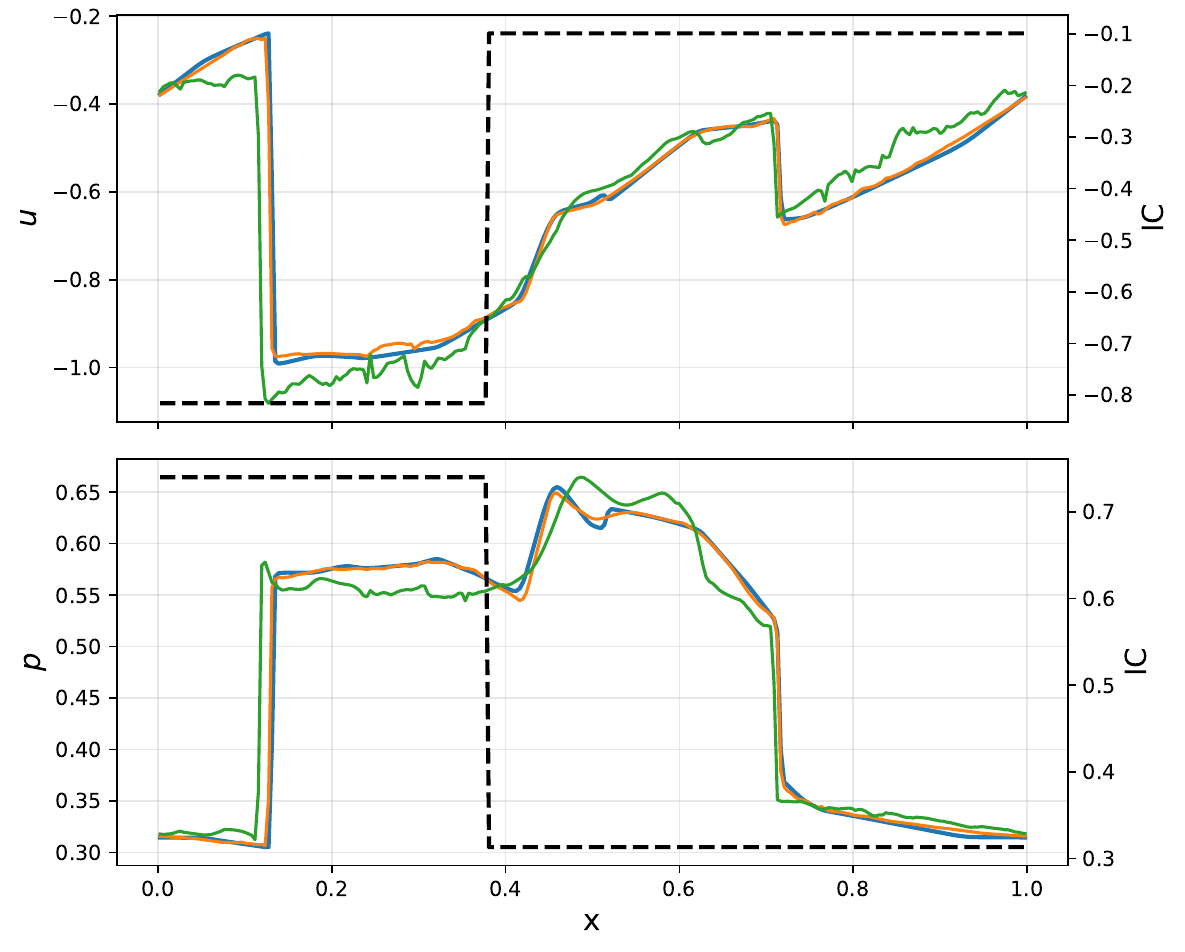}
        \caption{Final-time velocity and pressure \((u,p)\)}
    \end{subfigure}
    \caption{\textbf{Solution at \(T=1\) for the one-dimensional Euler equations.} Panel (a) compares the reference solution with the LGNO and FNO predictions in the conserved variables \((\rho,m,E)\), and panel (b) shows the recovered velocity and pressure \((u,p)\). In each subplot, the initial condition is shown as a dashed line on a separate right axis. LGNO remains consistently closer to the reference solution and places the shocks and contact discontinuities at their correct positions, whereas the FNO baseline exhibits larger oscillatory errors and fails to propagate the shocks at the correct speed.}
    \label{fig:euler1d_profiles}
\end{figure}

Figure~\ref{fig:euler1d_profiles} compares the conserved variables \((\rho,m,E)\) and the recovered velocity and pressure \((u,p)\) at the final time \(T=1\); in each subplot, the initial condition is overlaid as a dashed line on a separate right axis, plotted on its own scale. LGNO stays close to the reference across the conserved variables and yields more accurate velocity and pressure fields, capturing the shocks and contact discontinuities at their correct positions. The FNO baseline, by contrast, develops larger oscillatory errors and does not propagate the shocks at the correct speed. These oscillations are amplified in the velocity and pressure recovered from the conserved fields, while LGNO remains well behaved in both representations, giving sharper transitions with far fewer spurious oscillations.

\begin{figure}[tbp]
    \centering
    \begin{subfigure}[b]{0.31\textwidth}
        \centering
        \includegraphics[width=\textwidth]{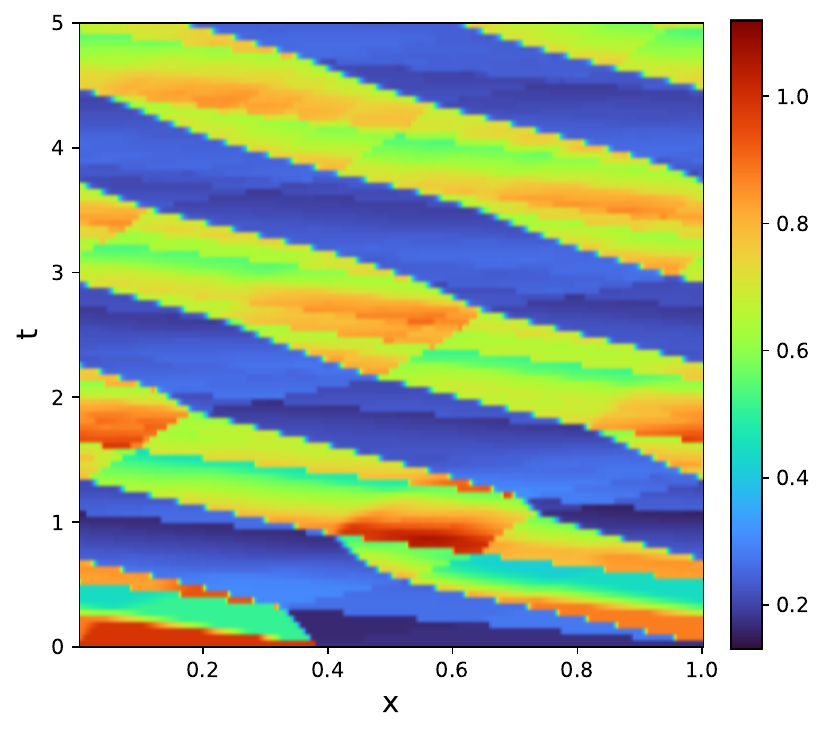}
         \vspace{-1.5em}
        \caption{Ref.~\(\rho(x,t)\)}
    \end{subfigure}
    \begin{subfigure}[b]{0.315\textwidth}
        \centering
        \includegraphics[width=\textwidth]{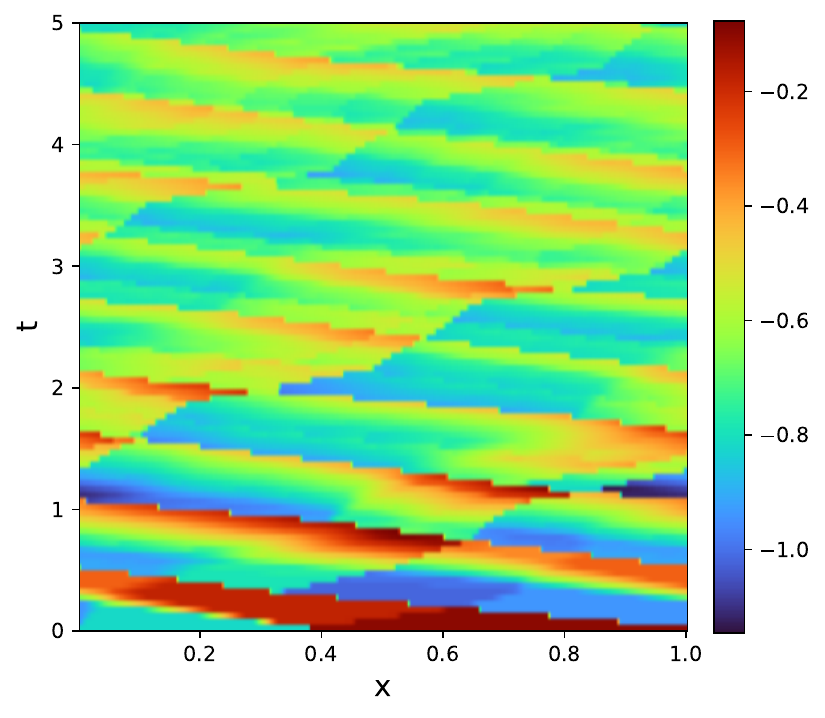}
        \vspace{-1.5em}
        \caption{Ref.~\(u(x,t)\)}
    \end{subfigure}
    \begin{subfigure}[b]{0.31\textwidth}
        \centering
        \includegraphics[width=\textwidth]{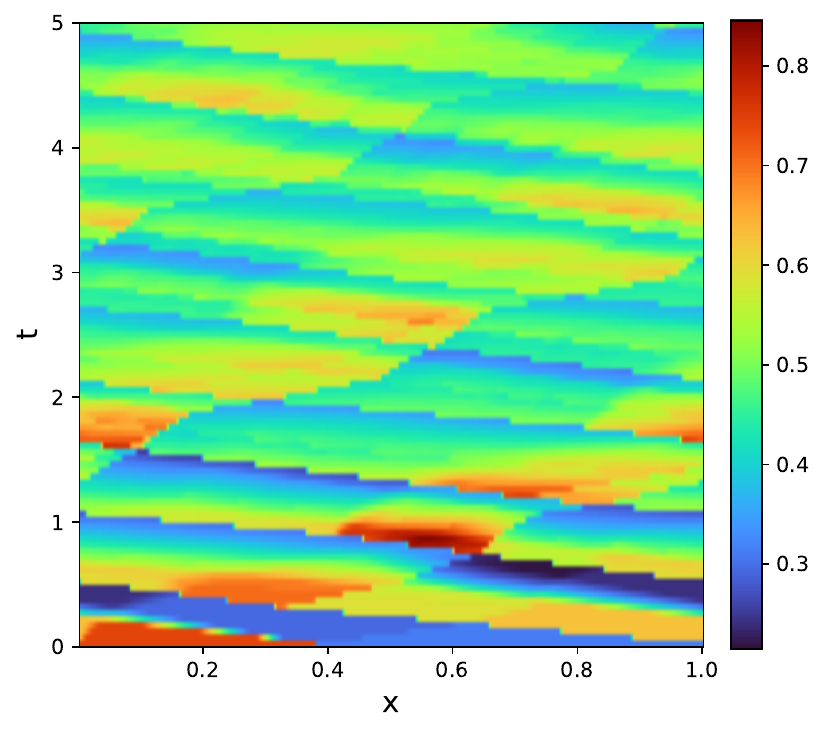}
         \vspace{-1.5em}
        \caption{Ref.~\(p(x,t)\)}
    \end{subfigure}
    \\[.5em]
    \begin{subfigure}[b]{0.31\textwidth}
        \centering
        \includegraphics[width=\textwidth]{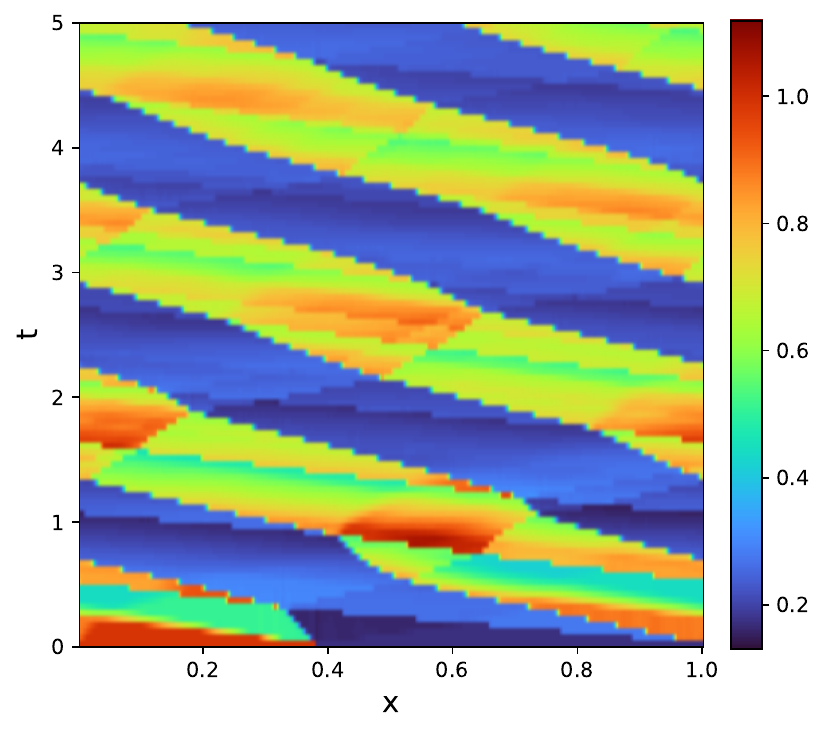}
        \vspace{-1.5em}
        \caption{LGNO \(\rho(x,t)\)}
    \end{subfigure}
    \begin{subfigure}[b]{0.315\textwidth}
        \centering
        \includegraphics[width=\textwidth]{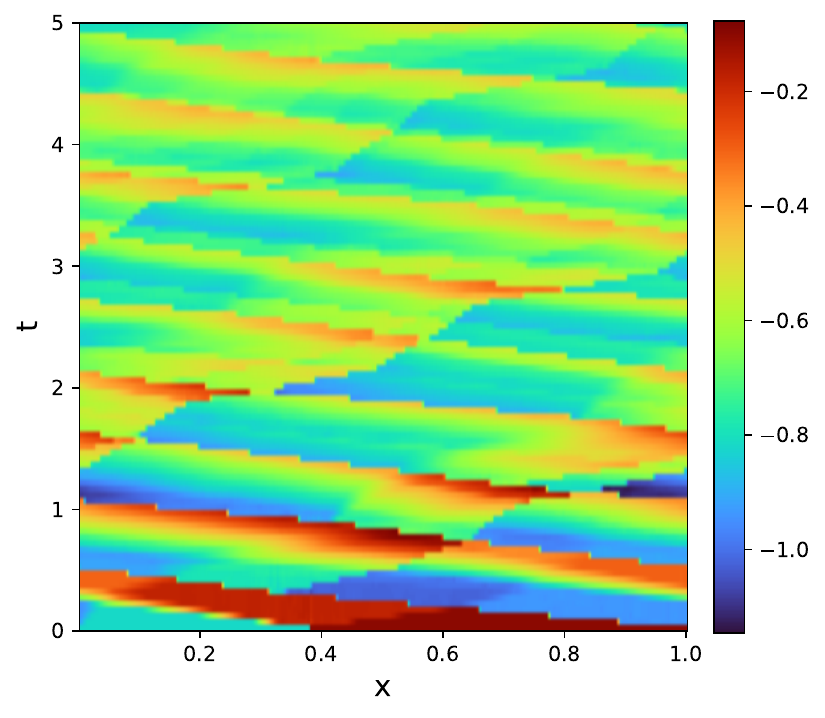}
        \vspace{-1.5em}
        \caption{LGNO \(u(x,t)\)}
    \end{subfigure}
    \begin{subfigure}[b]{0.31\textwidth}
        \centering
        \includegraphics[width=\textwidth]{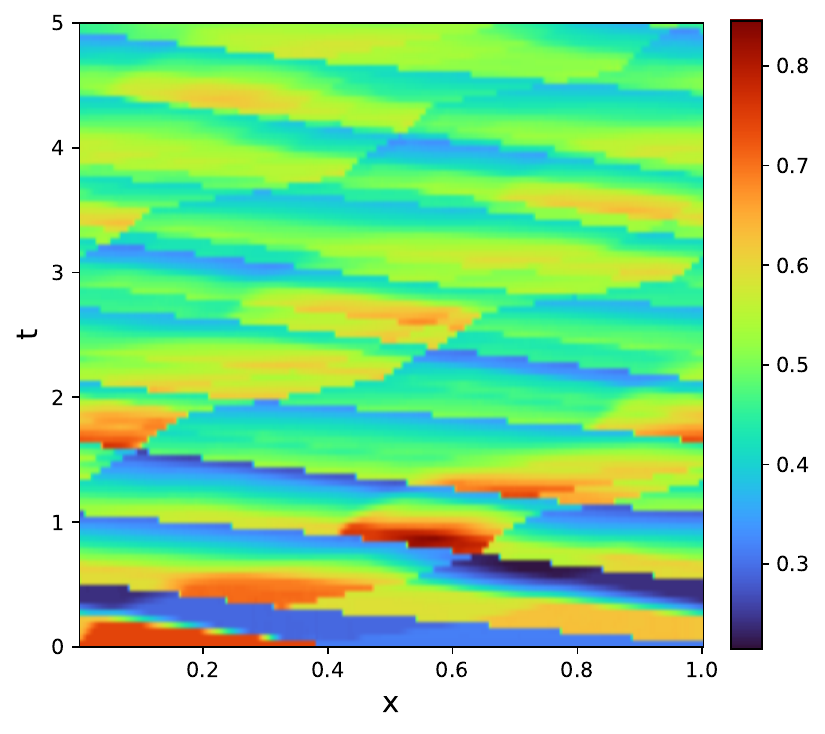}
        \vspace{-1.5em}
        \caption{LGNO \(p(x,t)\)}
    \end{subfigure}
    \\[.5em]
    \begin{subfigure}[b]{0.31\textwidth}
        \centering
        \includegraphics[width=\textwidth]{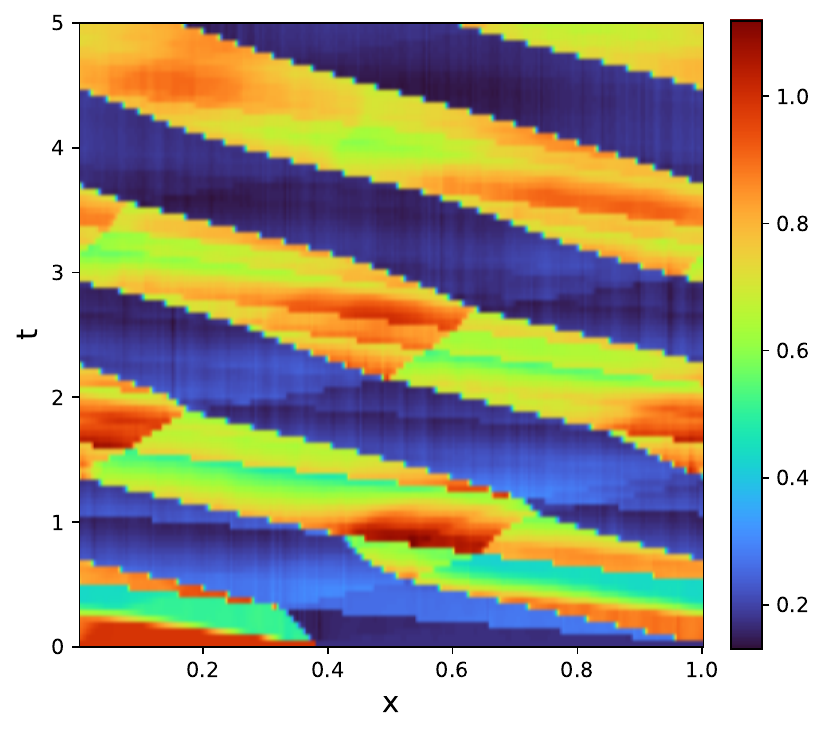}
        \vspace{-1.5em}
        \caption{FNO \(\rho(x,t)\)}
    \end{subfigure}
    \begin{subfigure}[b]{0.315\textwidth}
        \centering
        \includegraphics[width=\textwidth]{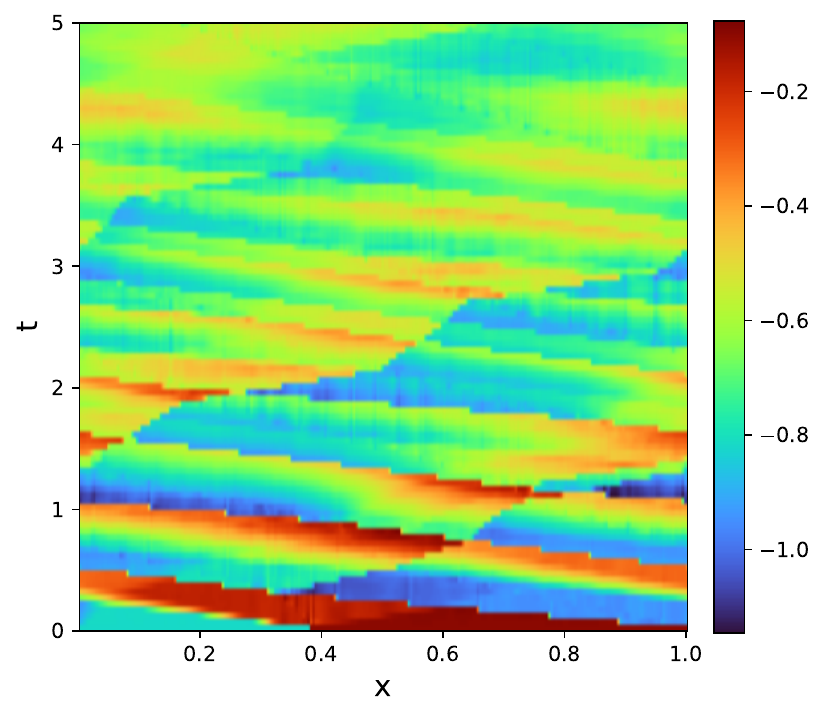}
        \vspace{-1.5em}
        \caption{FNO \(u(x,t)\)}
    \end{subfigure}
    \begin{subfigure}[b]{0.31\textwidth}
        \centering
        \includegraphics[width=\textwidth]{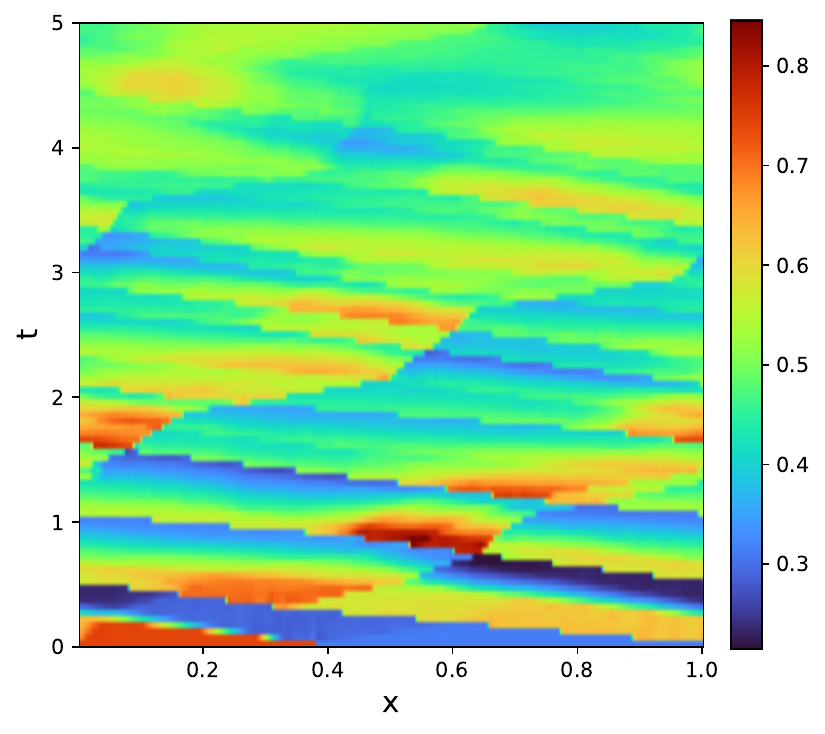}
        \vspace{-1.5em}
        \caption{FNO \(p(x,t)\)}
    \end{subfigure}
    \caption{\textbf{Spatiotemporal evolution for the one-dimensional Euler equations.} Panels (a)--(c) show the reference solution, panels (d)--(f) show LGNO, and panels (g)--(i) show the FNO baseline. Within each row, the columns correspond to the density \(\rho(x,t)\), velocity \(u(x,t)\), and pressure \(p(x,t)\), and all panels share a common contour range so that the three solutions are shown on the same scale. LGNO closely tracks the reference wave evolution across all three primitive variables, with the shock and contact fronts advancing at the correct speeds. The FNO baseline shows much larger errors than LGNO and becomes overly diffusive over the long rollout.}
    \label{fig:euler1d_xt}
\end{figure}

The spatiotemporal plots in Figure~\ref{fig:euler1d_xt} extend this to the full rollout up to \(T=5\), where all panels share a common contour range so that the three solutions are compared on the same scale. LGNO closely follows the reference wave structure and keeps the shocks and contacts on their correct trajectories. Under this common scale, the FNO baseline shows much larger errors than LGNO, with its discontinuity fronts advancing at incorrect speeds. The difference is more pronounced in the velocity and pressure fields, where the FNO baseline becomes overly diffusive, smearing sharp transitions in \(u\) and \(p\) into broad bands and washing out fine structures by the end of the rollout.

\subsection{Two-Dimensional Burgers Equation}
\label{subsec:results_burgers2d}

We then consider the two-dimensional Burgers equation
\begin{equation}
u_t + \left(\frac{u^2}{2}\right)_x + \left(\frac{u^2}{2}\right)_y = 0,
\quad (x,y)\in[0,1]^2.
\label{eq:burgers2d_pde}
\end{equation}
This extends the study to a genuinely two-dimensional setting, assessing whether LGNO retains its advantage when the solution develops nontrivial two-dimensional structure over long rollouts.

\begin{figure}[tbp]
    \centering
    \begin{subfigure}[b]{0.49\textwidth}
        \centering
        \includegraphics[width=\textwidth]{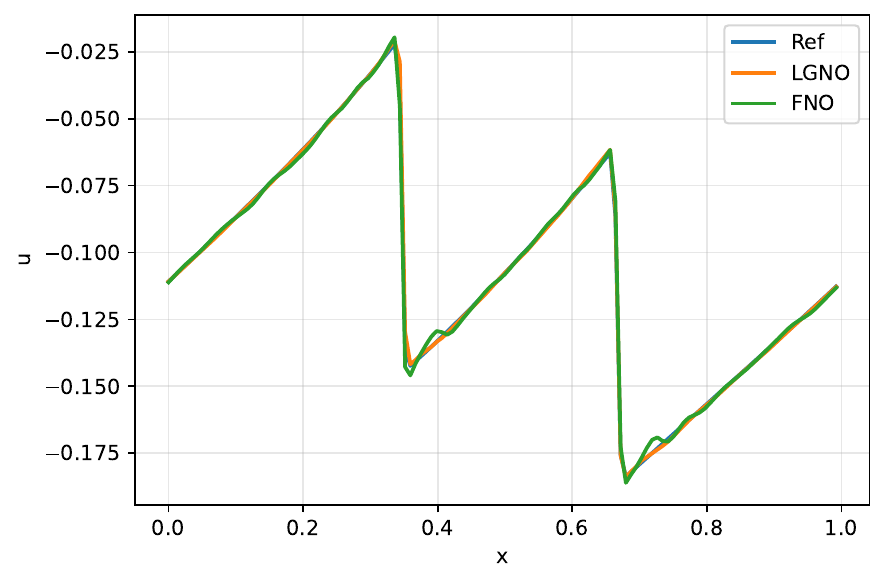}
        \caption{Solution slice at \(y=0.3\) and \(T=3\)}
    \end{subfigure}
    \hfill
    \begin{subfigure}[b]{0.49\textwidth}
        \centering
        \includegraphics[width=\textwidth]{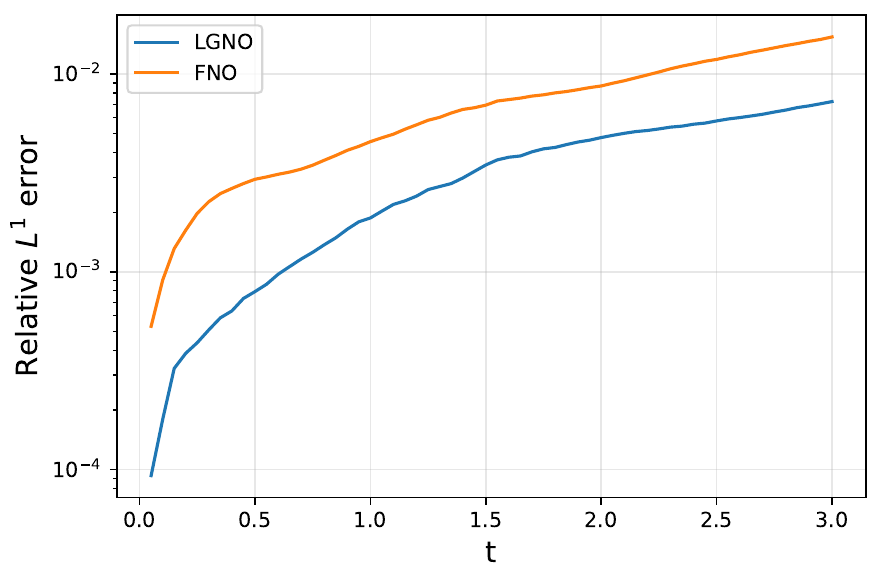}
        \caption{Relative \(L^1\) error versus time}
    \end{subfigure}
    \caption{\textbf{Long-horizon rollout diagnostics for the two-dimensional Burgers equation.} Panel (a) compares a representative one-dimensional solution slice at \(y=0.3\) and \(T=3\), and panel (b) shows the rollout relative \(L^1\) error over time. LGNO remains closer to the reference near steep transitions and achieves consistently smaller accumulated error than the FNO baseline.}
    \label{fig:burgers2d_rollout}
\end{figure}

Figure~\ref{fig:burgers2d_rollout} compares a representative one-dimensional slice of the solution at \(T=3\) with the rollout relative \(L^1\) error. LGNO remains close to the reference across the domain, including near steep transitions, whereas the FNO baseline shows larger deviations both at the sharp features and in neighboring smooth regions. The error curve confirms this, with LGNO staying consistently below the FNO baseline throughout the rollout.

\begin{figure}[tbp]
    \centering
    \begin{subfigure}[b]{0.315\textwidth}
        \centering
        \includegraphics[width=\textwidth]{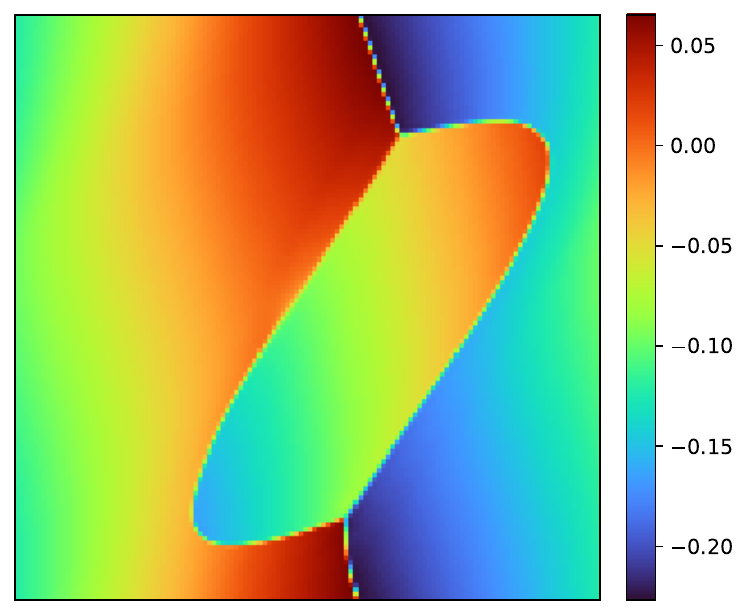}
        \caption{LGNO solution at \(T=3.0\)}
    \end{subfigure}
    \quad
    \begin{subfigure}[b]{0.31\textwidth}
        \centering
        \includegraphics[width=\textwidth]{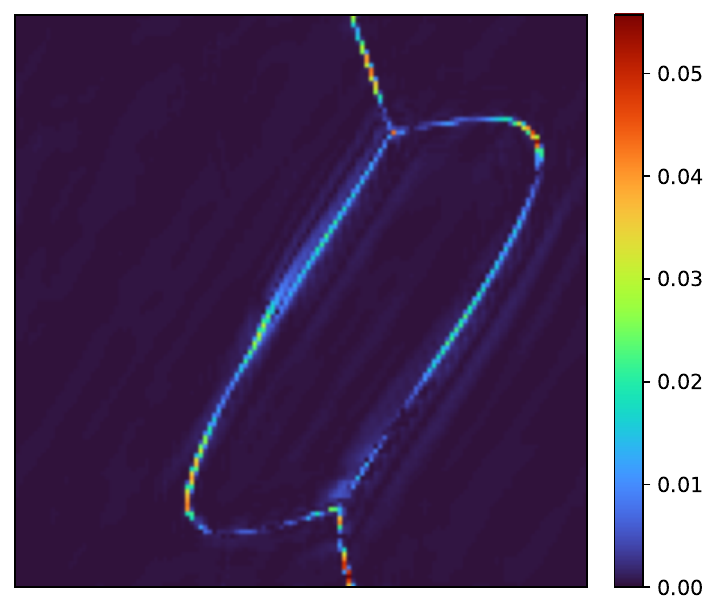}
        \caption{LGNO absolute error}
    \end{subfigure}
    \quad
    \begin{subfigure}[b]{0.31\textwidth}
        \centering
        \includegraphics[width=\textwidth]{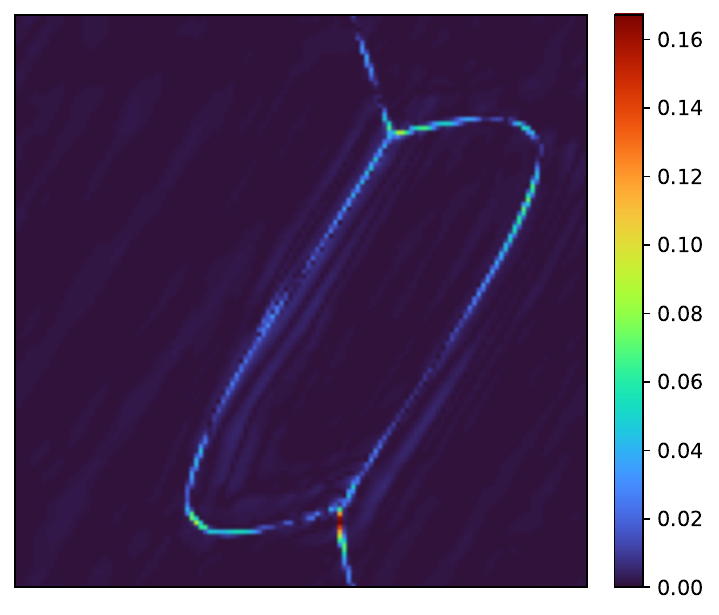}
        \caption{FNO absolute error}
    \end{subfigure}
    \caption{\textbf{Solution and error snapshots for the two-dimensional Burgers equation.} Panels (a)--(c) show the results at \(T=3.0\): panel (a) presents the LGNO solution, while panels (b) and (c) compare the absolute error fields of LGNO and the FNO baseline. The dominant errors concentrate near the evolved sharp fronts, but LGNO shows a visibly smaller error magnitude than the FNO baseline.}
    \label{fig:burgers2d_snapshots}
\end{figure}

Figure~\ref{fig:burgers2d_snapshots} gives a two-dimensional view at the final time \(T=3\). LGNO captures the evolved structures and curved fronts well; the dominant errors of both models concentrate near these sharp interfaces, but LGNO maintains a visibly smaller error magnitude than the FNO baseline. LGNO thus retains its advantage in the genuinely two-dimensional regime.

\subsection{Two-Dimensional Euler Equations with Periodic Boundary Conditions}
\label{subsec:euler2d_periodic}

We next consider the two-dimensional compressible Euler equations
\begin{equation}
\label{eq:euler2d_system}
\mathbf{U}_t+\mathbf{F}(\mathbf{U})_x+\mathbf{G}(\mathbf{U})_y=0,
\quad (x,y)\in[0,1]^2,
\end{equation}
where
\[
\mathbf{U}=(\rho,\rho u,\rho v,E)^{\top},
\quad 
\mathbf{F}(\mathbf{U})
=
(\rho u,\rho u^2+p,\rho uv,(E+p)u)^{\top},
\quad
\mathbf{G}(\mathbf{U})
=
(\rho v,\rho uv,\rho v^2+p,(E+p)v)^{\top}.
\]
The pressure is given by the ideal-gas equation of state
\begin{equation}
\label{eq:euler2d_pressure}
p=(\gamma-1)\left(E-\frac12\rho(u^2+v^2)\right),
\quad \gamma=1.4 .
\end{equation}
We first study this system under periodic boundary conditions. Table~\ref{tab:onestep_summary} already shows that LGNO achieves a substantially smaller one-step mean relative \(L^1\) error than the FNO baseline. To examine the rollout behavior, we show a representative autoregressive prediction to \(T=0.5\) with \(\Delta t=0.01\), corresponding to \(50\) rollout steps.

\begin{figure}[tbp]
    \centering
    \begin{subfigure}[b]{0.24\textwidth}
        \centering
        \includegraphics[width=\textwidth,trim={0 0 0 19pt},clip]{euler2d_ref_periodic_conserved_rho_WENO256.pdf}
        \vspace{-1.5em}
        \caption{WENO-Z \(256^2\) \(\rho\)}
    \end{subfigure}
    \begin{subfigure}[b]{0.24\textwidth}
        \centering
        \includegraphics[width=\textwidth,trim={0 0 0 19pt},clip]{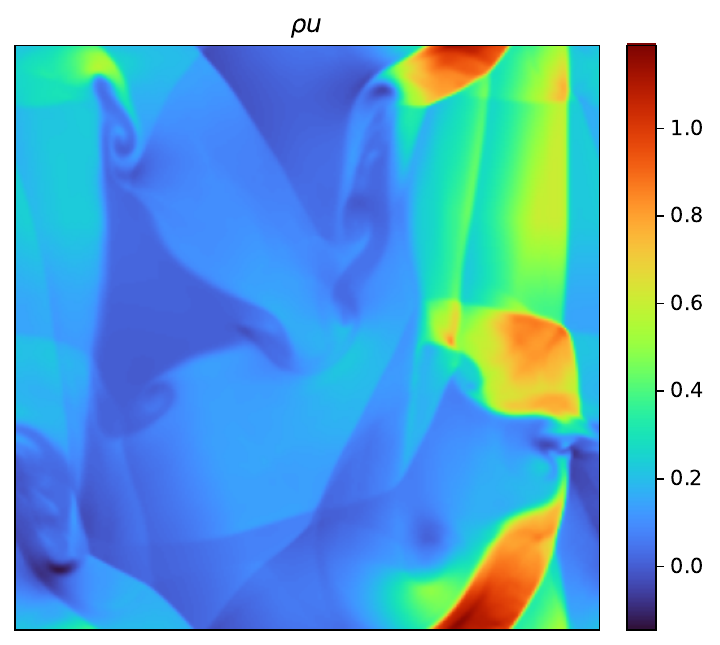}
        \vspace{-1.5em}
        \caption{WENO-Z \(256^2\) \(\rho u\)}
    \end{subfigure}
    \begin{subfigure}[b]{0.245\textwidth}
        \centering
        \includegraphics[width=\textwidth,trim={0 0 0 19pt},clip]{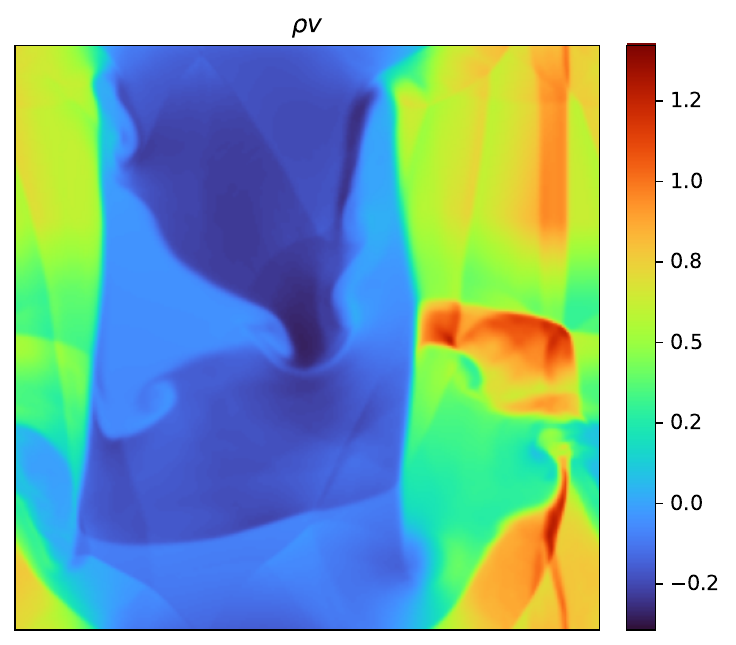}
        \vspace{-1.5em}
        \caption{WENO-Z \(256^2\) \(\rho v\)}
    \end{subfigure}
    \begin{subfigure}[b]{0.24\textwidth}
        \centering
        \includegraphics[width=\textwidth,trim={0 0 0 19pt},clip]{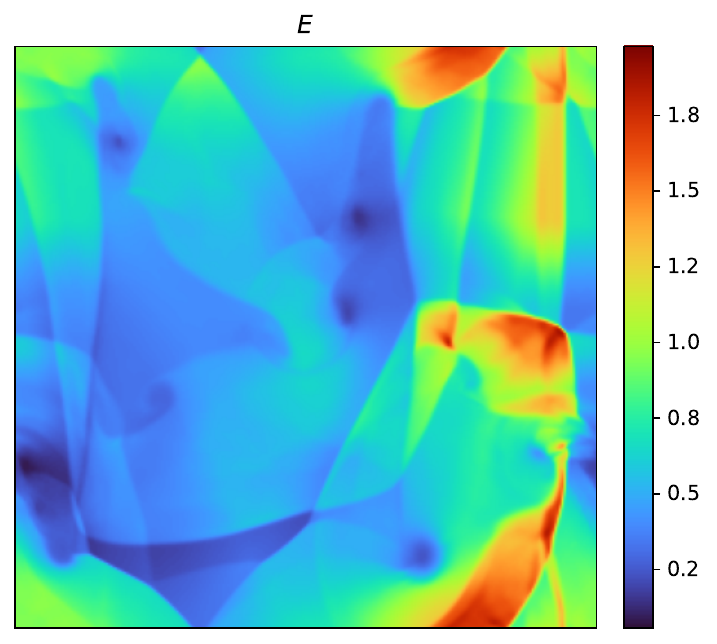}
        \vspace{-1.5em}
        \caption{WENO-Z \(256^2\) \(E\)}
    \end{subfigure}
    \\[.5em]
    \begin{subfigure}[b]{0.24\textwidth}
        \centering
        \includegraphics[width=\textwidth,trim={0 0 0 19pt},clip]{euler2d_ref_periodic_conserved_rho.pdf}
        \vspace{-1.5em}
        \caption{WENO-Z \(512^2\) \(\rho\)}
    \end{subfigure}
    \begin{subfigure}[b]{0.24\textwidth}
        \centering
        \includegraphics[width=\textwidth,trim={0 0 0 19pt},clip]{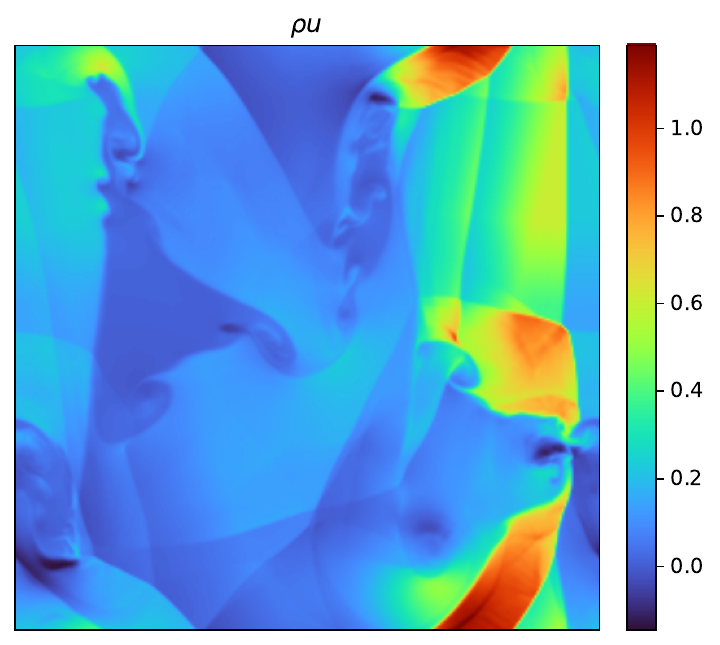}
        \vspace{-1.5em}
        \caption{WENO-Z \(512^2\) \(\rho u\)}
    \end{subfigure}
    \begin{subfigure}[b]{0.24\textwidth}
        \centering
        \includegraphics[width=\textwidth,trim={0 0 0 19pt},clip]{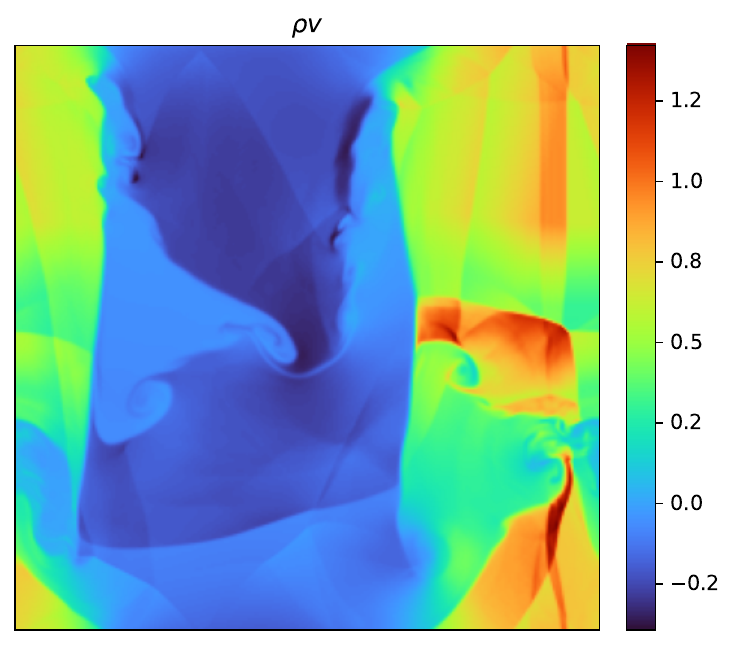}
        \vspace{-1.5em}
        \caption{WENO-Z \(512^2\) \(\rho v\)}
    \end{subfigure}
    \begin{subfigure}[b]{0.245\textwidth}
        \centering
        \includegraphics[width=\textwidth,trim={0 0 0 19pt},clip]{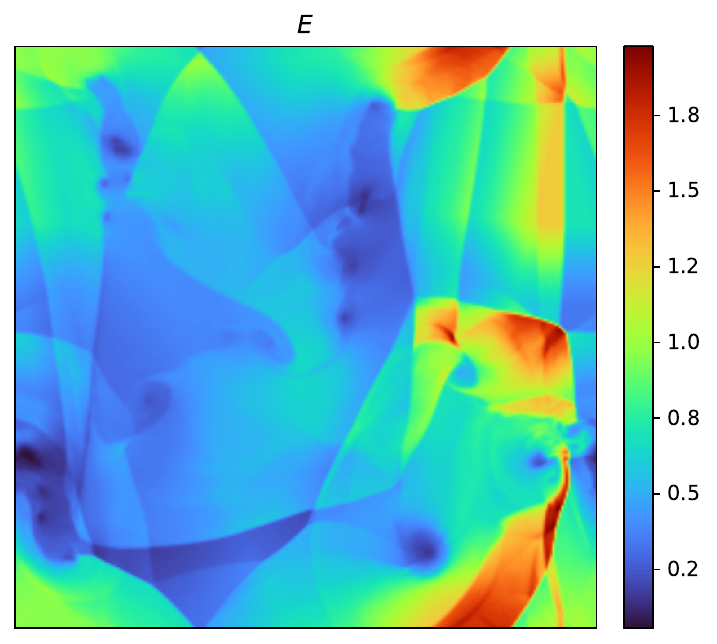}
        \vspace{-1.5em}
        \caption{WENO-Z \(512^2\) \(E\)}
    \end{subfigure}
    \\[.5em]
    \begin{subfigure}[b]{0.24\textwidth}
        \centering
        \includegraphics[width=\textwidth,trim={0 0 0 19pt},clip]{euler2d_hybrid_periodic_conserved_rho.pdf}
        \vspace{-1.5em}
        \caption{LGNO \(256^2\) \(\rho\)}
    \end{subfigure}
    \begin{subfigure}[b]{0.24\textwidth}
        \centering
        \includegraphics[width=\textwidth,trim={0 0 0 19pt},clip]{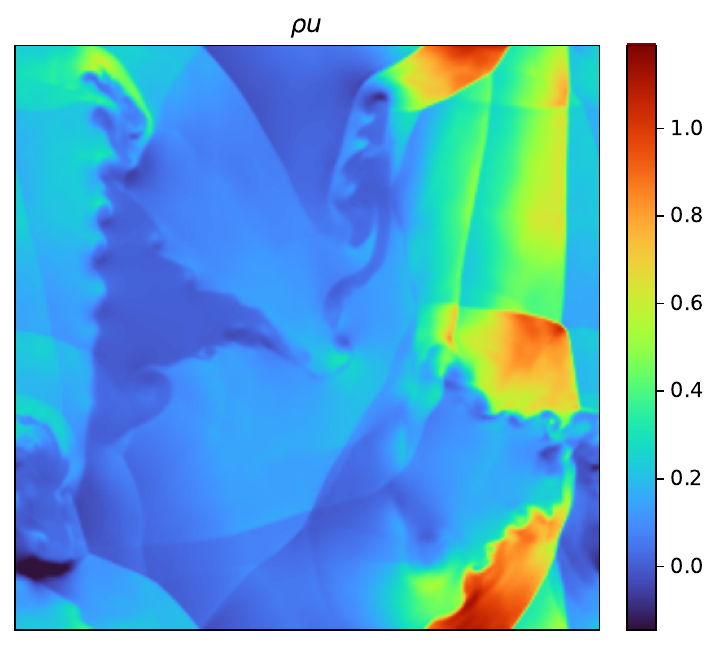}
        \vspace{-1.5em}
        \caption{LGNO \(256^2\) \(\rho u\)}
    \end{subfigure}
    \begin{subfigure}[b]{0.245\textwidth}
        \centering
        \includegraphics[width=\textwidth,trim={0 0 0 19pt},clip]{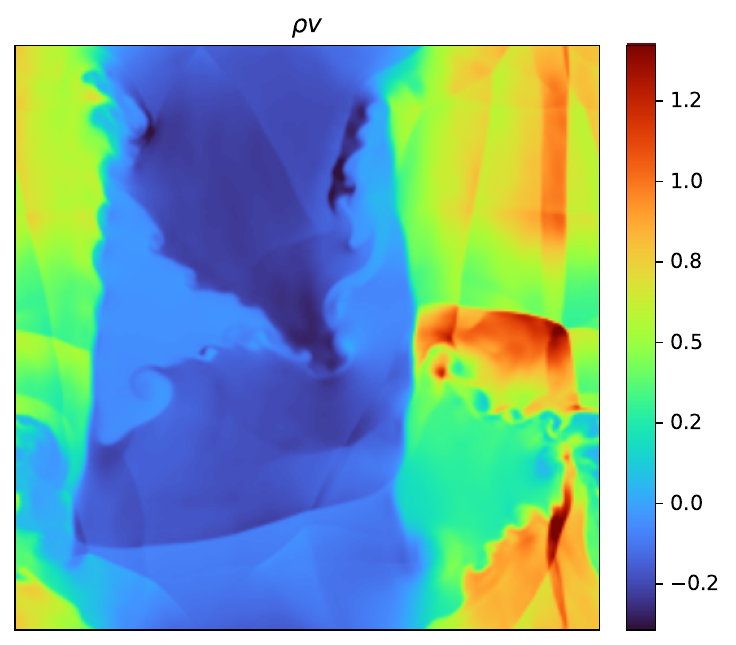}
        \vspace{-1.5em}
        \caption{LGNO \(256^2\) \(\rho v\)}
    \end{subfigure}
    \begin{subfigure}[b]{0.24\textwidth}
        \centering
        \includegraphics[width=\textwidth,trim={0 0 0 19pt},clip]{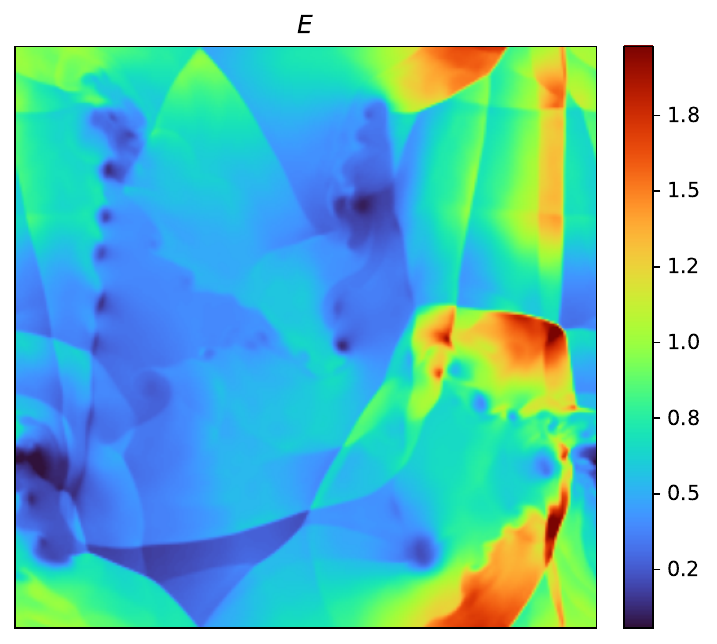}
        \vspace{-1.5em}
        \caption{LGNO \(256^2\) \(E\)}
    \end{subfigure}
    \\[.5em]
    \begin{subfigure}[b]{0.24\textwidth}
        \centering
        \includegraphics[width=\textwidth,trim={0 0 0 19pt},clip]{euler2d_fno_periodic_conserved_rho.pdf}
        \vspace{-1.5em}
        \caption{FNO \(256^2\) \(\rho\)}
    \end{subfigure}
    \begin{subfigure}[b]{0.24\textwidth}
        \centering
        \includegraphics[width=\textwidth,trim={0 0 0 19pt},clip]{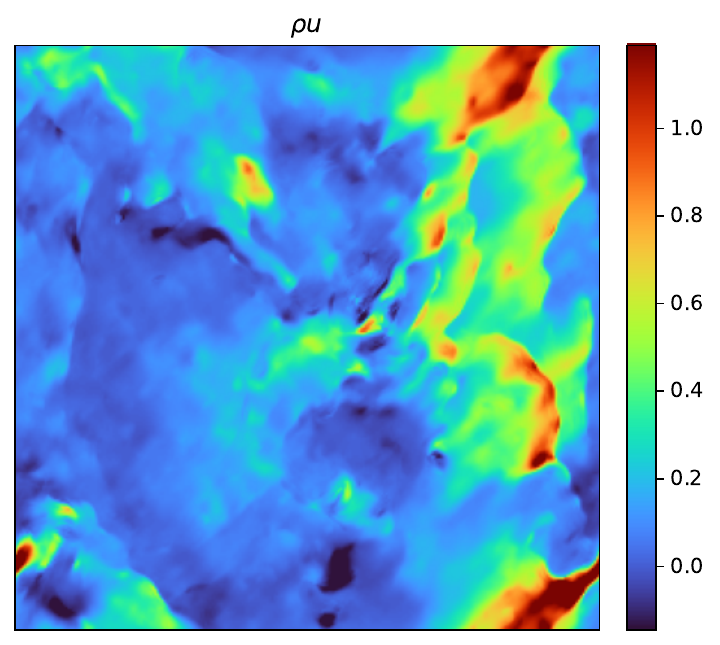}
        \vspace{-1.5em}
        \caption{FNO \(256^2\) \(\rho u\)}
    \end{subfigure}
    \begin{subfigure}[b]{0.245\textwidth}
        \centering
        \includegraphics[width=\textwidth,trim={0 0 0 19pt},clip]{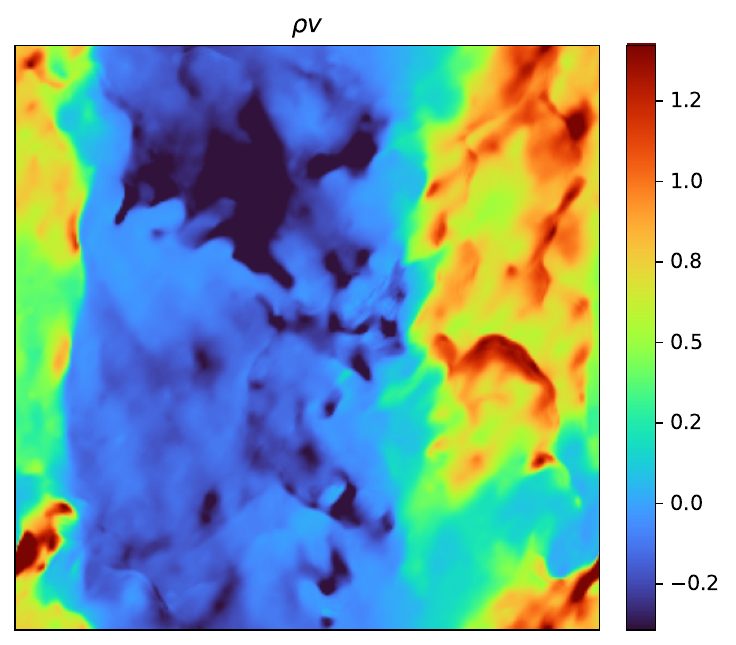}
        \vspace{-1.5em}
        \caption{FNO \(256^2\) \(\rho v\)}
    \end{subfigure}
    \begin{subfigure}[b]{0.24\textwidth}
        \centering
        \includegraphics[width=\textwidth,trim={0 0 0 19pt},clip]{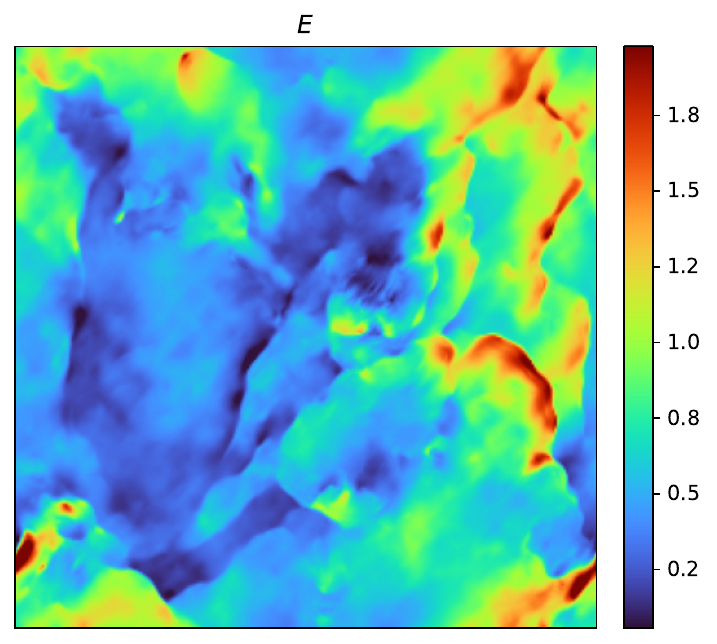}
        \vspace{-1.5em}
        \caption{FNO \(256^2\) \(E\)}
    \end{subfigure}
    \caption{\textbf{Representative rollout results for the two-dimensional Euler equations with periodic boundary conditions at \(T=0.5\).} The first row shows the WENO-Z solution on the same \(256^2\) mesh as the learned models, and the second row shows the high-resolution WENO-Z solution on a \(512^2\) mesh. The third and fourth rows show the LGNO and FNO predictions on the \(256^2\) mesh. From left to right, the columns correspond to density, \(x\)-momentum, \(y\)-momentum, and total energy. For each variable, the four panels in the corresponding column use the same colorbar range, allowing direct comparison across methods. LGNO better preserves the localized structures than the FNO baseline, which loses coherent wave patterns and develops stronger oscillatory artifacts. Notably, LGNO on the \(256^2\) mesh also retains fine vortical and roll-up structures that are smeared in the WENO-Z solutions even on the finer \(512^2\) mesh, reflecting its lower long-time numerical dissipation.}
    \label{fig:euler2d_periodic_conserved}
\end{figure}

Figure~\ref{fig:euler2d_periodic_conserved} compares the conserved variables from WENO-Z on the \(256^2\) mesh, the high-resolution \mbox{WENO-Z} reference on the \(512^2\) mesh, and the learned predictions on the \(256^2\) mesh. For each variable, all panels use the colorbar limits determined by the high-resolution reference. LGNO stays qualitatively closer to the reference wave pattern and preserves localized structures more clearly than the FNO baseline; notably, the density and momentum fields retain fine vortical and roll-up structures that are weakened in the WENO-Z solutions even on the \(512^2\) mesh. This indicates a relatively low-dissipation character and points to LGNO's potential as a learned solver. It should not be read as LGNO being more accurate than the high-resolution reference, but rather as its preserving sharper localized structures under long rollout. The FNO prediction, by contrast, loses much of the coherent structure and develops stronger oscillatory artifacts. 

We also verified the primitive variables \((\rho,u,v,p)\) recovered from the conserved fields (not shown for brevity): LGNO remains stable after this conversion, whereas the FNO baseline becomes severely unstable in the velocity fields, with isolated unphysical values of \(u\) and \(v\) reaching magnitudes of order \(10^3\).

\subsection{Two-Dimensional Euler Equations with Outflow Boundary Conditions}
\label{subsec:euler2d_outflow}

We finally consider the same two-dimensional Euler system \eqref{eq:euler2d_system}--\eqref{eq:euler2d_pressure} under outflow boundary conditions. This setting is more challenging due to the boundary condition. In this experiment, the network is trained to predict the primitive variables \((\rho,u,v,p)\) rather than the conserved variables \((\rho,\rho u,\rho v,E)\). Near the outflow boundaries, the local branch uses the learned ghost-padding strategy of Section~\ref{sec:outflow_padding} in place of periodic extension, and the positivity of \(\rho\) and \(p\) is enforced through a generalized softplus parametrization with \(\beta=100\) on the corresponding output channels. All models are trained with learning rate \(5\times10^{-4}\), and additional implementation details are provided in~\ref{app:euler2d_outflow_details}.

\begin{figure}[tbp]
    \centering
    \begin{subfigure}[b]{0.24\textwidth}
        \centering
        \includegraphics[width=\textwidth]{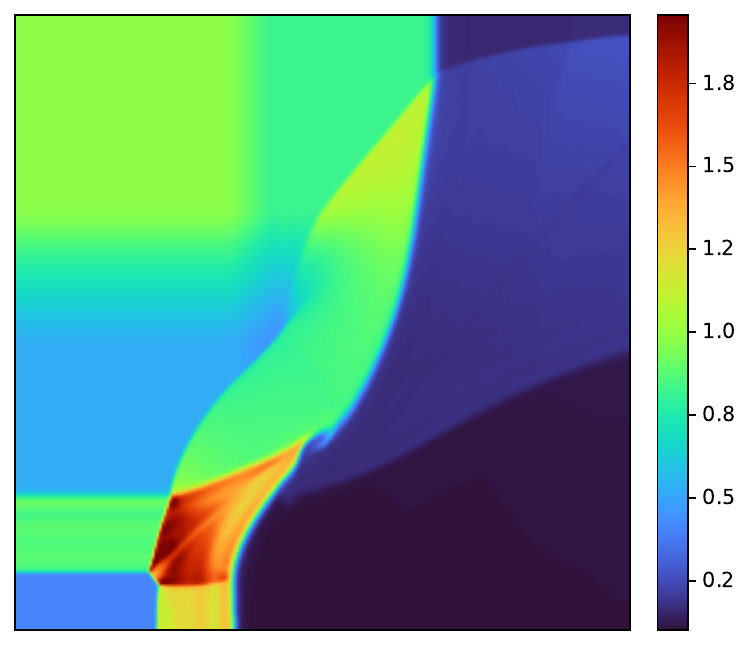}
        \vspace{-1.5em}
        \caption{WENO-Z \(256^2\) \(\rho\)}
    \end{subfigure}
    \begin{subfigure}[b]{0.24\textwidth}
        \centering
        \includegraphics[width=\textwidth]{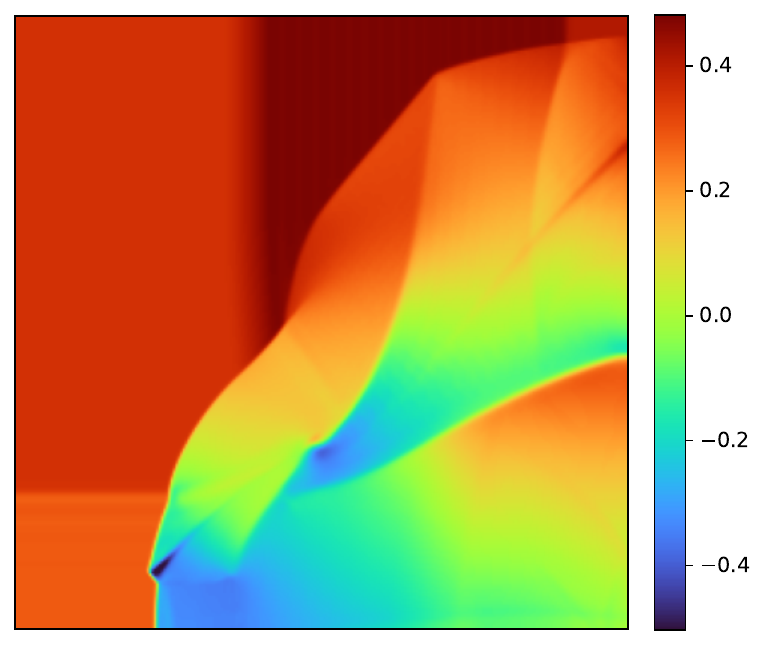}
        \vspace{-1.5em}
        \caption{WENO-Z \(256^2\) \(u\)}
    \end{subfigure}
    \begin{subfigure}[b]{0.24\textwidth}
        \centering
        \includegraphics[width=\textwidth]{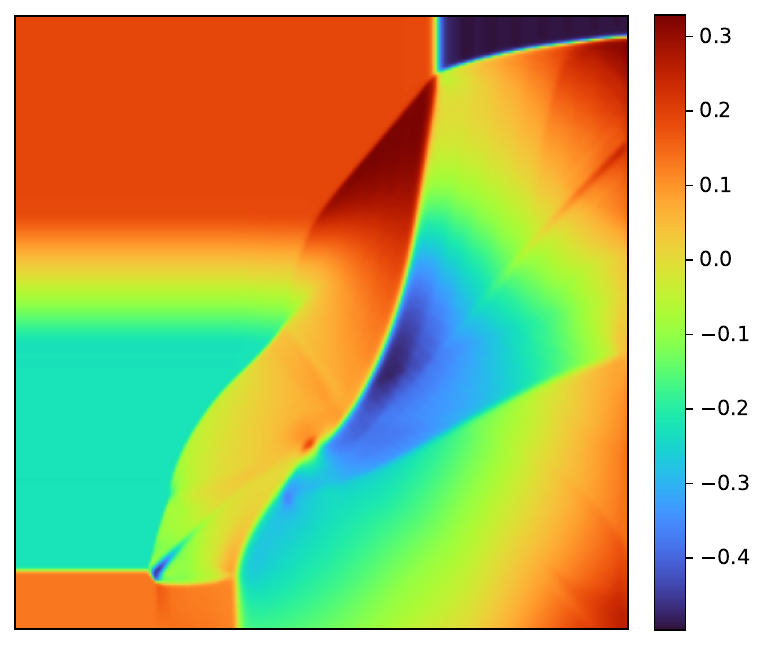}
        \vspace{-1.5em}
        \caption{WENO-Z \(256^2\) \(v\)}
    \end{subfigure}
    \begin{subfigure}[b]{0.24\textwidth}
        \centering
        \includegraphics[width=\textwidth]{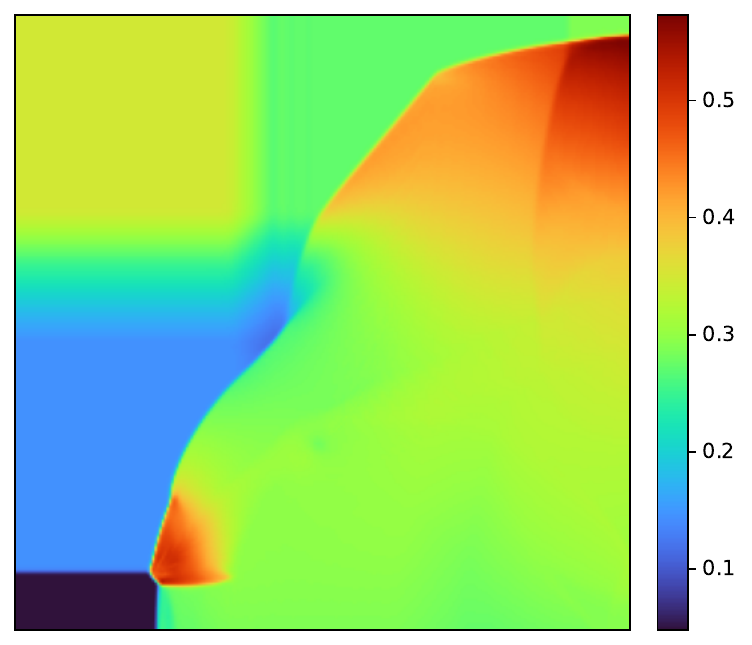}
        \vspace{-1.5em}
        \caption{WENO-Z \(256^2\) \(p\)}
    \end{subfigure}
    \\[.5em]
    \begin{subfigure}[b]{0.24\textwidth}
        \centering
        \includegraphics[width=\textwidth]{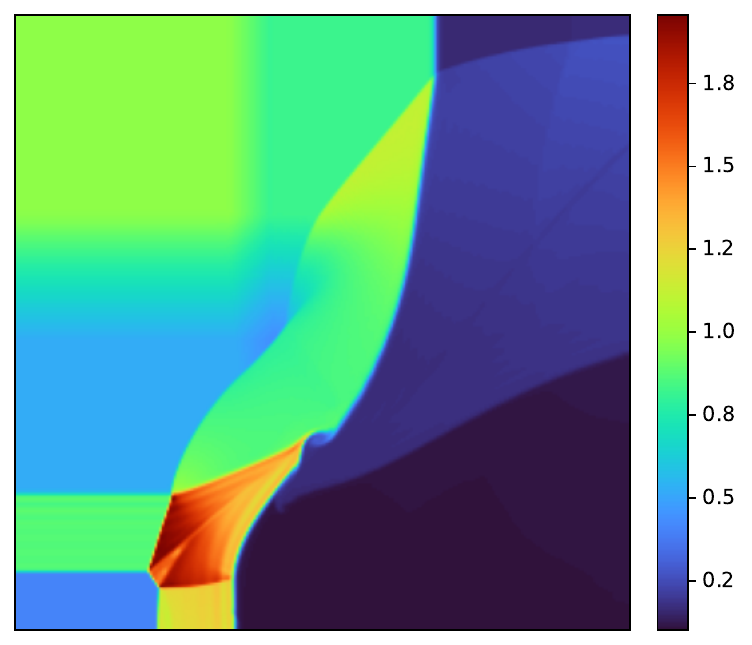}
        \vspace{-1.5em}
        \caption{WENO-Z \(512^2\) \(\rho\)}
    \end{subfigure}
    \begin{subfigure}[b]{0.24\textwidth}
        \centering
        \includegraphics[width=\textwidth]{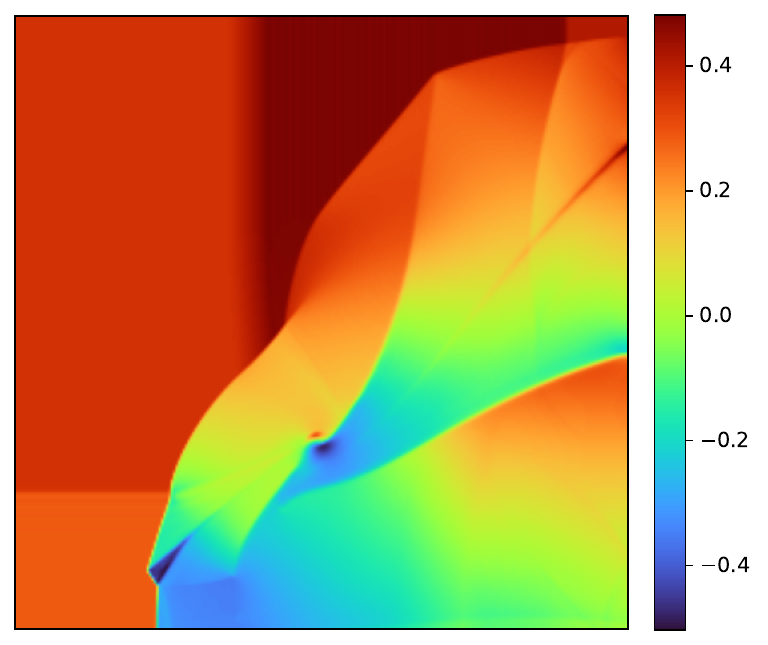}
        \vspace{-1.5em}
        \caption{WENO-Z \(512^2\) \(u\)}
    \end{subfigure}
    \begin{subfigure}[b]{0.24\textwidth}
        \centering
        \includegraphics[width=\textwidth]{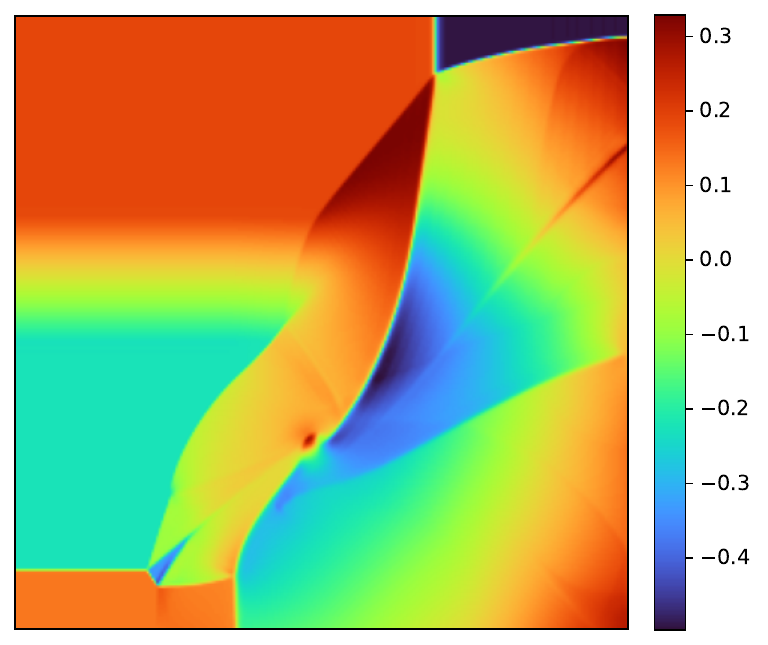}
        \vspace{-1.5em}
        \caption{WENO-Z \(512^2\) \(v\)}
    \end{subfigure}
    \begin{subfigure}[b]{0.24\textwidth}
        \centering
        \includegraphics[width=\textwidth]{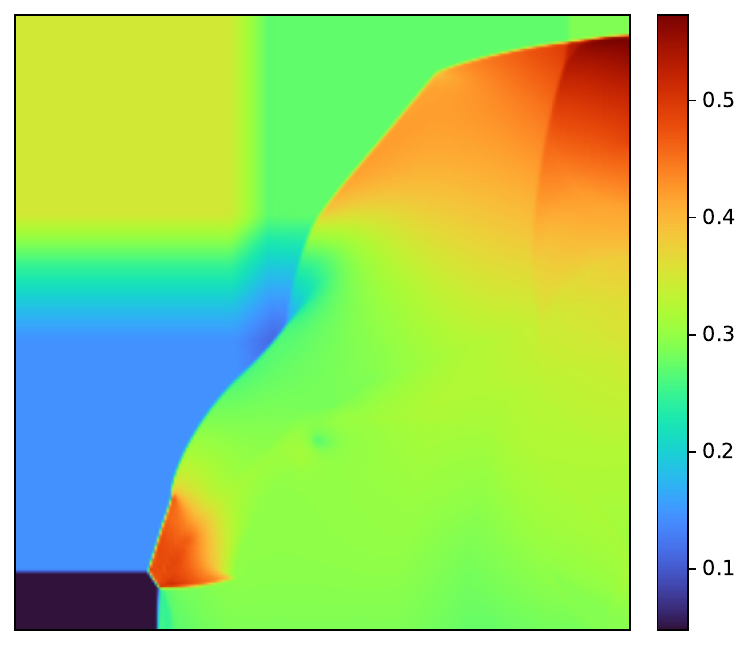}
        \vspace{-1.5em}
        \caption{WENO-Z \(512^2\) \(p\)}
    \end{subfigure}
    \\[.5em]
    \begin{subfigure}[b]{0.24\textwidth}
        \centering
        \includegraphics[width=\textwidth]{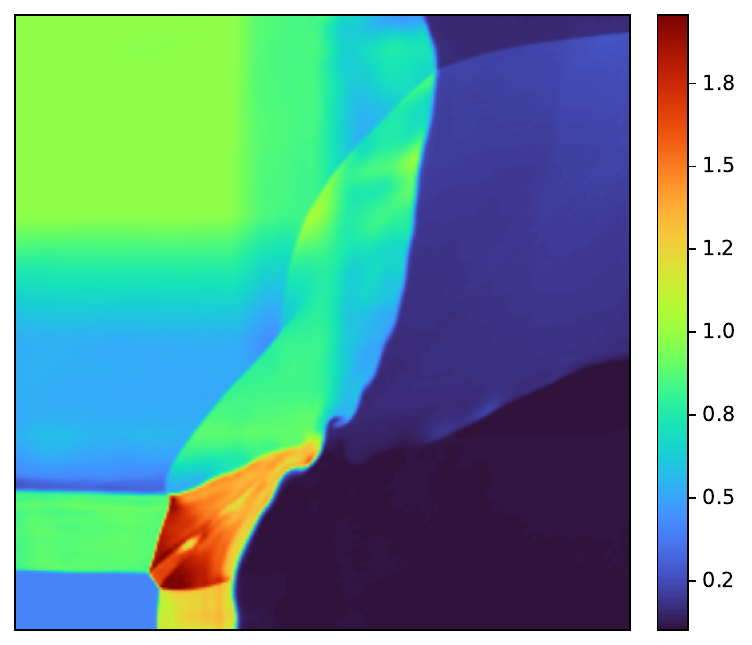}
        \vspace{-1.5em}
        \caption{LGNO \(256^2\) \(\rho\)}
    \end{subfigure}
    \begin{subfigure}[b]{0.24\textwidth}
        \centering
        \includegraphics[width=\textwidth]{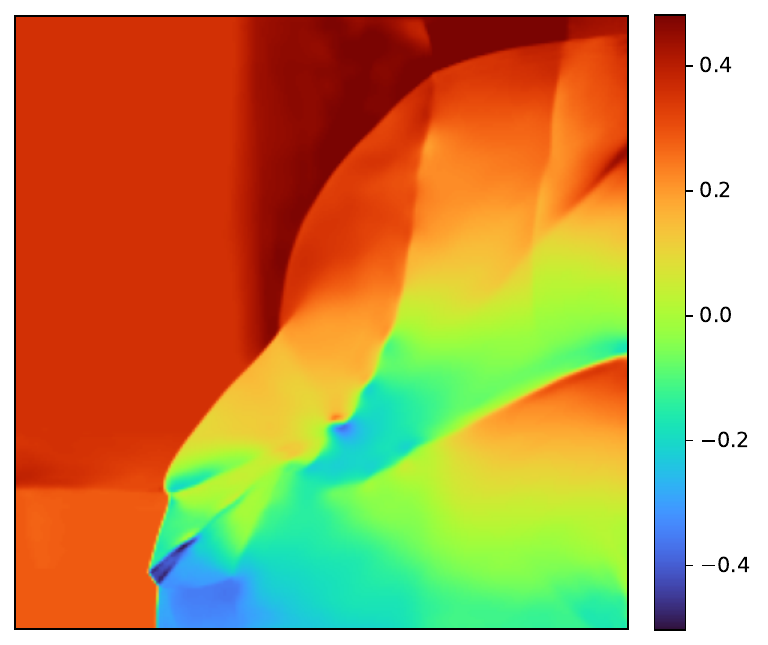}
        \vspace{-1.5em}
        \caption{LGNO \(256^2\) \(u\)}
    \end{subfigure}
    \begin{subfigure}[b]{0.24\textwidth}
        \centering
        \includegraphics[width=\textwidth]{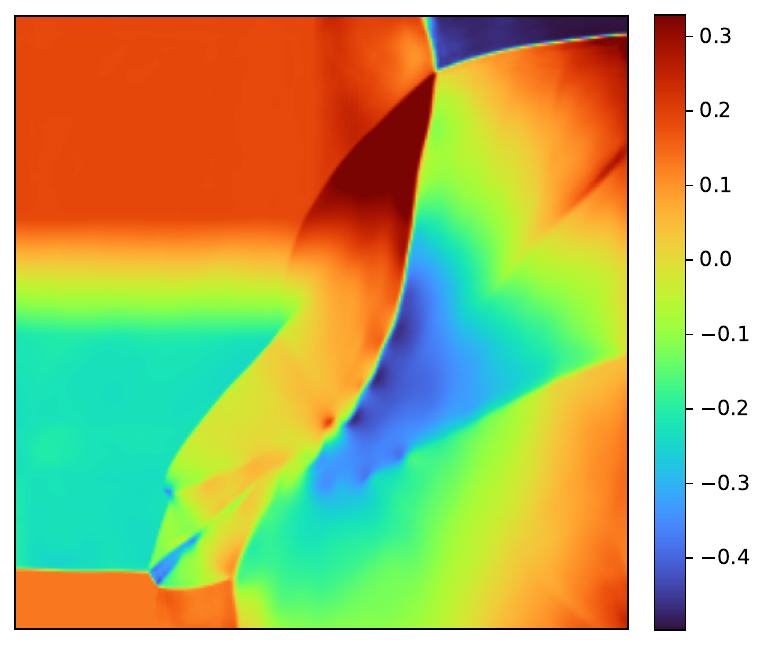}
        \vspace{-1.5em}
        \caption{LGNO \(256^2\) \(v\)}
    \end{subfigure}
    \begin{subfigure}[b]{0.24\textwidth}
        \centering
        \includegraphics[width=\textwidth]{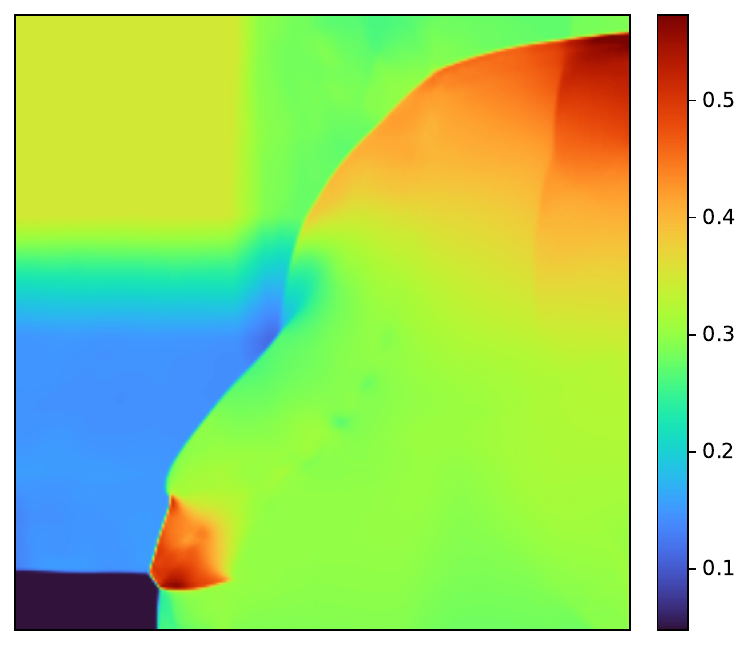}
        \vspace{-1.5em}
        \caption{LGNO \(256^2\) \(p\)}
    \end{subfigure}
    \\[.5em]
    \begin{subfigure}[b]{0.24\textwidth}
        \centering
        \includegraphics[width=\textwidth]{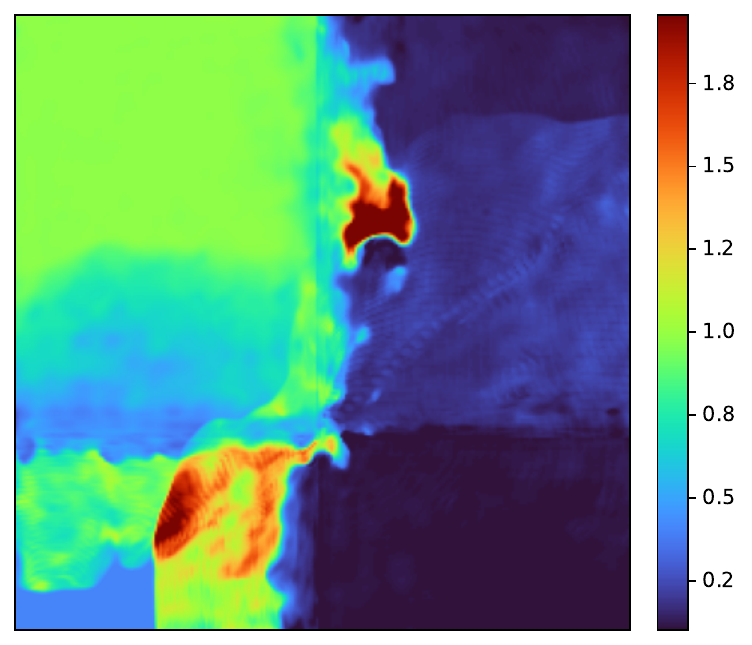}
        \vspace{-1.5em}
        \caption{FNO \(256^2\) \(\rho\)}
    \end{subfigure}
    \begin{subfigure}[b]{0.24\textwidth}
        \centering
        \includegraphics[width=\textwidth]{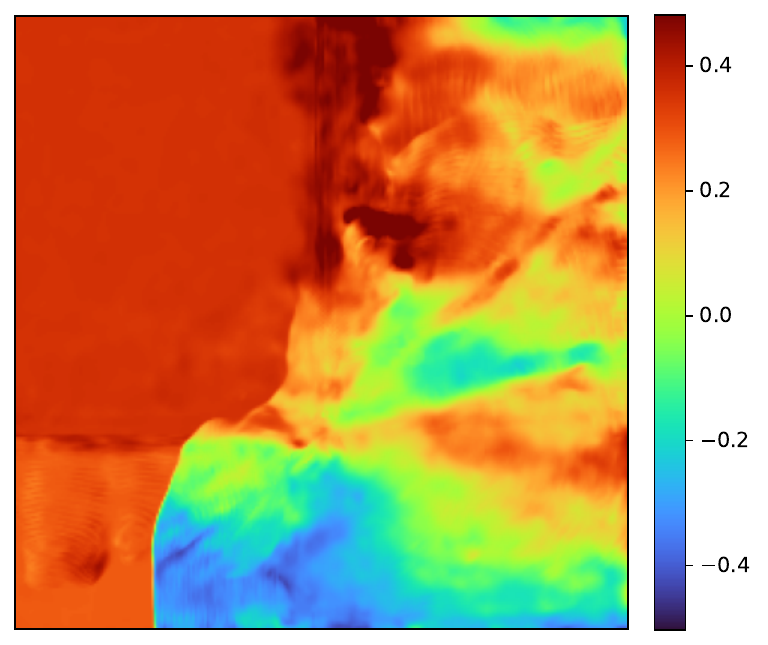}
        \vspace{-1.5em}
        \caption{FNO \(256^2\) \(u\)}
    \end{subfigure}
    \begin{subfigure}[b]{0.24\textwidth}
        \centering
        \includegraphics[width=\textwidth]{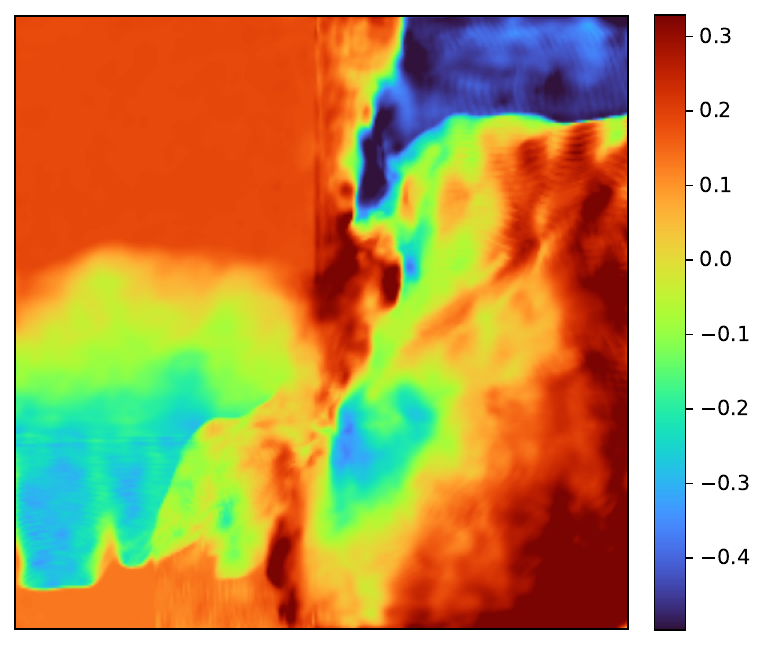}
        \vspace{-1.5em}
        \caption{FNO \(256^2\) \(v\)}
    \end{subfigure}
    \begin{subfigure}[b]{0.24\textwidth}
        \centering
        \includegraphics[width=\textwidth]{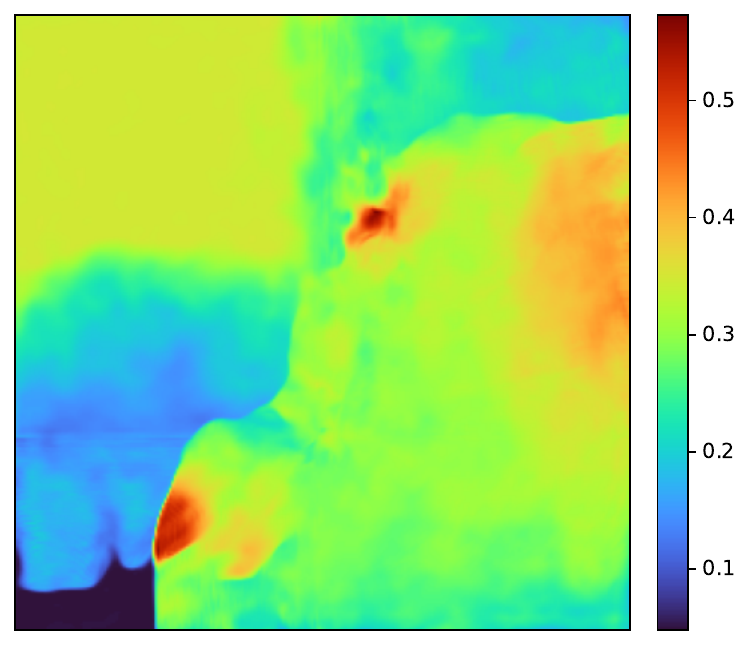}
        \vspace{-1.5em}
        \caption{FNO \(256^2\) \(p\)}
    \end{subfigure}
    \caption{\textbf{Representative rollout results for the two-dimensional Euler equations with outflow boundary conditions at \(T=0.4\).} Panels (a)--(d) show the WENO-Z solution on the \(256^2\) mesh, panels (e)--(h) show the WENO-Z reference solution on the \(512^2\) mesh, panels (i)--(l) show the LGNO prediction on the \(256^2\) mesh, and panels (m)--(p) show the FNO prediction on the \(256^2\) mesh, all after \(40\) autoregressive rollout steps. The columns correspond to the primitive variables \(\rho\), \(u\), \(v\), and \(p\). For each primitive variable, the four panels in the corresponding column use the same colorbar range, allowing direct comparison across methods. LGNO better preserves the main wave structures and fine patterns, while the FNO baseline develops stronger distortions and oscillatory artifacts during rollout.}
    \label{fig:euler2d_outflow_rollout}
\end{figure}

Figure~\ref{fig:euler2d_outflow_rollout} shows a representative autoregressive rollout at \(T=0.4\), corresponding to \(40\) rollout steps. For each primitive variable, the four method panels share the same colorbar range and are therefore directly comparable. LGNO preserves the main wave interactions and discontinuity structures more accurately than the FNO baseline: the density field retains fine roll-up patterns near the interacting wave fronts, whereas the FNO prediction loses much of the coherent structure and develops noisy artifacts. LGNO thus provides a more stable and less dissipative rollout for this challenging outflow problem.

\section{Conclusion}
\label{sec:conc}

We have proposed the Local--Global Neural Operator (LGNO), a neural operator for hyperbolic conservation laws that learns a one-step discrete flow map by coupling a global Fourier branch, which captures smooth large-scale dynamics and long-range interactions, with a local multiresolution branch, which resolves sharp localized structures such as shocks and contact discontinuities. The two branches are combined through a multiplicative interaction, and the model is trained with a one-step loss that augments a physical-space prediction term with a high-frequency spectral penalty to suppress spurious oscillations near steep fronts. For average-preserving boundary conditions, a mean-zero projection makes the learned update exactly conservative, and a simple learned padding extends the local branch to outflow boundaries.

Across several one- and two-dimensional benchmark problems, LGNO consistently outperforms parameter-matched FNO baselines, reducing one-step errors by factors of two to five and remaining markedly more accurate over long autoregressive rollouts. It advances shocks and contact discontinuities at the correct speeds and preserves fine vortical and roll-up structures that the FNO baseline smears or destabilizes. Most notably, although it is trained only on short-time data from a high-order WENO-Z scheme, the long-time rollout of LGNO on a coarse grid exhibits lower numerical dissipation than the same scheme run on a finer grid, while remaining orders of magnitude cheaper to evaluate. This behavior arises because a high-order shock-capturing scheme injects numerical dissipation at every step of a long integration, whereas the learned flow map carries only the low dissipation present in its short-time training data. These results suggest that, with an appropriate architecture and training objective, learned operators have the potential to control long-time numerical dissipation better than the conventional shock-capturing schemes that generate their training data.

Several directions remain open. On the theoretical side, it would be valuable to analyze the approximation and stability properties of the local--global architecture and to characterize when the learned flow map stays low-dissipative over long horizons. On the practical side, we plan to extend LGNO to more complex systems and higher dimensions, to enforce local conservation and bound-preserving structure, to study robustness under extrapolation beyond the training distribution, and to develop efficient training strategies for large-scale problems.


\section{Reproducibility}
\label{sec:reproducibility}
To support reproducibility, we release our implementation of LGNO, the training and evaluation scripts, the WENO-Z data-generation code, and the trained model weights for all benchmark experiments reported in this paper. The code and pretrained checkpoints are available at \url{https://github.com/shanxue-w/ConservationLaws}. The repository includes the configuration files needed to reproduce the reported autoregressive rollouts and figures.


\section*{Acknowledgment}
This research used resources of the National Energy Research Scientific Computing Center, a DOE Office of Science User Facility supported by the Office of Science of the U.S. Department of Energy under Contract No. DE-AC02-05CH11231 using NERSC award ASCR-ERCAP0037224.


\appendix
\newcommand{\shockquad}{%
  \begingroup
  \setlength{\arraycolsep}{3pt}%
  \renewcommand{\arraystretch}{0.85}%
  \raisebox{-0.7ex}{$
  \begin{array}{c|c}
  \scriptstyle \mathcal R_2 & \scriptstyle \mathcal R_1\\[-0.2em]\hline
  \scriptstyle \mathcal R_3 & \scriptstyle \mathcal R_4
  \end{array}$}%
  \endgroup
}

\newcommand{\shocktitle}[1]{%
  \par\medskip
  \noindent\textsc{#1.}\quad \shockquad \quad: the initial data are
}

\newenvironment{shockstates}{%
  \[
  \begin{array}{r@{\;}c@{\;}l@{\qquad}
                r@{\;}c@{\;}l@{\qquad\qquad}
                r@{\;}c@{\;}l@{\qquad}
                r@{\;}c@{\;}l}
}{%
  \end{array}
  \]
}

\newcommand{\shockrhofig}[6]{%
\begin{figure}[!htbp]
\centering
\setlength{\tabcolsep}{2pt}
\renewcommand{\arraystretch}{0.92}
\begin{tabular}{c c c c}
 & Ref. \(\rho\) & LGNO \(\rho\) & FNO \(\rho\) \\
S1 &
\includegraphics[width=0.29\textwidth]{euler2d_rho_ref_#2_#1_01_#6.pdf} &
\includegraphics[width=0.29\textwidth]{euler2d_rho_hybrid_#2_#1_01_#6.pdf} &
\includegraphics[width=0.29\textwidth]{euler2d_rho_fno_#2_#1_01_#6.pdf}
\\[-0.25em]
S2 &
\includegraphics[width=0.29\textwidth]{euler2d_rho_ref_#2_#1_02_#6.pdf} &
\includegraphics[width=0.29\textwidth]{euler2d_rho_hybrid_#2_#1_02_#6.pdf} &
\includegraphics[width=0.29\textwidth]{euler2d_rho_fno_#2_#1_02_#6.pdf}
\\[-0.25em]
S3 &
\includegraphics[width=0.29\textwidth]{euler2d_rho_ref_#2_#1_03_#6.pdf} &
\includegraphics[width=0.29\textwidth]{euler2d_rho_hybrid_#2_#1_03_#6.pdf} &
\includegraphics[width=0.29\textwidth]{euler2d_rho_fno_#2_#1_03_#6.pdf}
\\[-0.25em]
S4 &
\includegraphics[width=0.29\textwidth]{euler2d_rho_ref_#2_#1_04_#6.pdf} &
\includegraphics[width=0.29\textwidth]{euler2d_rho_hybrid_#2_#1_04_#6.pdf} &
\includegraphics[width=0.29\textwidth]{euler2d_rho_fno_#2_#1_04_#6.pdf}
\end{tabular}
\caption{#5}
\label{fig:euler2d_#1_shock_#2_#6}
\end{figure}
}

\section{Generalization Test: Two-Dimensional Euler Shock Tests}
\label{app:euler2d_shock_tests}

To further study the robustness of the learned flow map for the two-dimensional Euler equations, we include four shock tests based on standard quadrant Riemann benchmarks~\cite{euler_2d_example}. \emph{These tests do not appear in the training set and thus probe generalization to unseen initial data; since the states remain within the training distribution, this is an interpolation rather than an extrapolation test.} They assess whether the learned rollouts reproduce the main shock and contact structures, especially their locations, under strong wave interactions.

\paragraph{Test setup}
The initial data are prescribed by piecewise constant quadrant states in primitive variables. We divide the domain into four regions,
\[
\mathcal R_1=\{x\ge 1/2,\ y\ge 1/2\},\qquad
\mathcal R_2=\{x<1/2,\ y\ge 1/2\},
\]
\[
\mathcal R_3=\{x<1/2,\ y<1/2\},\qquad
\mathcal R_4=\{x\ge 1/2,\ y<1/2\}.
\]
On each region \(\mathcal R_i\), the primitive variables are set to a constant state,
\[
(p,\rho,u,v)(x,y)=(p_i,\rho_i,u_i,v_i),
\qquad (x,y)\in \mathcal R_i,\quad i=1,2,3,4.
\]

We label the four cases S1--S4 and run each under both outflow and periodic boundary conditions. With outflow conditions, the tests are the usual quadrant shock-interaction setup, evolved to \(T=0.25\) so the main waves stay away from the boundary. With periodic conditions, the same quadrant data are periodized; we therefore use a longer rollout to \(T=0.5\) to further stress the learned periodic flow map under strong discontinuities and interacting waves. The periodic case is a stress test of the learned flow map, not a canonical self-similar Riemann problem.

\paragraph{Shock-state definitions}
The four shock tests are specified below. In each display, the small \(2\times2\) block indicates the spatial arrangement of the four regions: the upper row corresponds to \(y\ge 1/2\), and the lower row corresponds to \(y<1/2\). The primitive states are listed in the order \((p,\rho,u,v)\).

\shocktitle{Shock test S1}
\begin{shockstates}
p_2&=&0.4     & \rho_2&=&0.5197 & p_1&=&1   & \rho_1&=&1 \\
u_2&=&-0.7259 & v_2&=&0        & u_1&=&0   & v_1&=&0 \\[0.35em]
p_3&=&1       & \rho_3&=&1      & p_4&=&0.4 & \rho_4&=&0.5197 \\
u_3&=&-0.7259 & v_3&=&-0.7259  & u_4&=&0   & v_4&=&-0.7259
\end{shockstates}

\shocktitle{Shock test S2}
\begin{shockstates}
p_2&=&0.4     & \rho_2&=&0.5197 & p_1&=&1   & \rho_1&=&1 \\
u_2&=&-0.6259 & v_2&=&0.1       & u_1&=&0.1 & v_1&=&0.1 \\[0.35em]
p_3&=&0.4     & \rho_3&=&0.8    & p_4&=&0.4 & \rho_4&=&0.5197 \\
u_3&=&0.1     & v_3&=&0.1       & u_4&=&0.1 & v_4&=&-0.6259
\end{shockstates}

\shocktitle{Shock test S3}
\begin{shockstates}
p_2&=&1       & \rho_2&=&1      & p_1&=&0.4 & \rho_1&=&0.5197 \\
u_2&=&-0.6259 & v_2&=&0.1      & u_1&=&0.1 & v_1&=&0.1 \\[0.35em]
p_3&=&1       & \rho_3&=&0.8    & p_4&=&1   & \rho_4&=&1 \\
u_3&=&0.1     & v_3&=&0.1      & u_4&=&0.1 & v_4&=&-0.6259
\end{shockstates}

\shocktitle{Shock test S4}
\begin{shockstates}
p_2&=&1      & \rho_2&=&1      & p_1&=&0.4 & \rho_1&=&0.5313 \\
u_2&=&0.7276 & v_2&=&0        & u_1&=&0   & v_1&=&0 \\[0.35em]
p_3&=&1      & \rho_3&=&0.8    & p_4&=&1   & \rho_4&=&1 \\
u_3&=&0      & v_3&=&0        & u_4&=&0   & v_4&=&0.7276
\end{shockstates}

\paragraph{Outflow boundary condition}
We first evaluate the four tests with outflow boundary conditions, rolling out the learned models autoregressively to \(T=0.25\) and comparing against WENO-Z reference solutions. We report density contours, which emphasize shock and contact locations, together with density fields, which give a complementary view of the global distribution.

Figure~\ref{fig:euler2d_outflow_shock_contour_25} compares the density contours, a direct check of shock and contact locations, since a learned flow map that does not explicitly enforce local conservation may in principle produce incorrect shock speeds. Across S1--S4, LGNO has a comparable contour plot as the WENO-Z reference, including the dominant shock and contact locations, with discrepancies confined to strong wave-interaction regions where minor contour distortion and mild smoothing appear. The FNO baseline, in contrast, develops substantial spurious oscillations and distorted wave patterns, especially near the central interaction zones, so LGNO avoids the shock-location errors seen in FNO.

Figure~\ref{fig:euler2d_outflow_shock_colormap_25} shows the corresponding density fields. LGNO stays close to the reference in the large-scale distribution and retains the main high- and low-density regions across all four tests, whereas the FNO baseline shows clear unphysical artifacts: localized overshoots, oscillatory structures, and large-scale distortion. LGNO is thus markedly more robust than FNO for shock-dominated outflow rollouts of the two-dimensional Euler equations.

\paragraph{Symmetry-consistent scaling}
The compressible Euler equations admit a two-parameter scaling invariance. We use this invariance to align the prescribed shock states with the magnitudes seen during training, while preserving the fact that the rescaled variables still solve the same Euler system. Specifically, for any constants \(s>0\) and \(k>0\), suppose that \((\rho,u,v,p)(\mathbf{x},t)\) solves the system~\eqref{eq:euler2d_system}--\eqref{eq:euler2d_pressure}. Define the rescaled variables
\[
\tilde\rho(\mathbf{x},t)=s\rho(\mathbf{x},t/k),\qquad
\tilde u(\mathbf{x},t)=\frac{1}{k}u(\mathbf{x},t/k),\qquad
\tilde v(\mathbf{x},t)=\frac{1}{k}v(\mathbf{x},t/k),\qquad
\tilde p(\mathbf{x},t)=\frac{s}{k^2}p(\mathbf{x},t/k).
\]
Then \((\tilde\rho,\tilde u,\tilde v,\tilde p)\) also satisfies the same two-dimensional compressible Euler equations on the same spatial domain. Indeed, with \(\tau=t/k\), the continuity equation is multiplied by the common factor \(s/k\). The two momentum equations are multiplied by the common factor \(s/k^2\), since
\[
\tilde\rho\tilde u=\frac{s}{k}\rho u,\qquad
\tilde\rho\tilde v=\frac{s}{k}\rho v,
\]
and
\[
\tilde\rho\tilde u^2+\tilde p
=
\frac{s}{k^2}(\rho u^2+p),\qquad
\tilde\rho\tilde u\tilde v
=
\frac{s}{k^2}\rho uv,\qquad
\tilde\rho\tilde v^2+\tilde p
=
\frac{s}{k^2}(\rho v^2+p).
\]
Similarly, the total energy rescales as
\[
\tilde E
=
\frac{\tilde p}{\gamma-1}
+
\frac12\tilde\rho(\tilde u^2+\tilde v^2)
=
\frac{s}{k^2}E,
\]
so the energy equation is multiplied by the common factor \(s/k^3\). Therefore all rescaled conservation laws reduce to the original Euler equations evaluated at \((\mathbf{x},\tau)\), and hence the rescaled variables remain an exact solution of the Euler system.

Under this scaling, the factor \(s\) rescales density and pressure together, while the factor \(k\) divides the velocities and the sound speed \(c=\sqrt{\gamma p/\rho}\) by \(k\). Consequently, the Mach number \(\sqrt{u^2+v^2}/c\) and the wave structure are preserved, while the flow evolves \(k\) times more slowly in time. This allows the prescribed states to be rescaled to the training range without changing the underlying Euler dynamics.

We use this invariance in the outflow tests at magnitudes compatible with the learned primitive-variable models. Density and pressure are first multiplied by a common positive factor \(s_{\rho p}\),
\[
(\rho,u,v,p)\mapsto (s_{\rho p}\rho,u,v,s_{\rho p}p),
\]
and, for the learned models, we additionally introduce a positive velocity scale \(k\) and evaluate on
\[
(s_{\rho p}\rho,u,v,s_{\rho p}p)
\mapsto
\left(
s_{\rho p}\rho,\frac{u}{k},\frac{v}{k},
\frac{s_{\rho p}p}{k^2}
\right).
\]
At the continuous level, this scaling is valid for any \(k>0\): the rescaled flow evolves more slowly by the time factor \(k\). In the discrete rollout experiments, we choose \(k\) to be an integer so that the corresponding time-rescaled trajectory can be evaluated by taking \(k\) times as many internal time steps with the same step size. The final prediction is then mapped back to the original primitive variables before visualization. This brings the prescribed states into the range seen during data generation while, by the invariance above, leaving the underlying Euler dynamics unchanged, so it does not impact the qualitative behavior of the prediction.

No such scaling is used in the periodic tests, where the models are evaluated directly on the original variables; the periodic results therefore give a cleaner assessment of generalization to these unseen configurations.

\paragraph{Periodic boundary condition}
We next repeat the four cases with periodic boundary conditions. Periodic extension of the quadrant states introduces additional discontinuities across the domain boundaries, making these rollouts especially challenging; with no boundary to avoid, we evolve them to the longer final time \(T=0.5\) to make it even more challenging.

Figure~\ref{fig:euler2d_periodic_shock_contour_50} compares the density contours. The longer horizon produces more complex wave interactions than the outflow tests, yet LGNO still captures the dominant contour geometry of the reference, including the principal shock and contact structures; local distortion and fine-scale artifacts appear near strong interaction regions, but the main patterns remain recognizable. The FNO baseline produces much denser spurious contours and loses much of the reference shock geometry.

Figure~\ref{fig:euler2d_periodic_shock_colormap_50} shows the corresponding density fields. LGNO retains the large-scale density organization across S1--S4, with only localized deviations after the longer rollout, while the FNO baseline shows stronger oscillations, distorted density regions, and a more severe loss of coherent structure. Consistent with the low-dissipation behavior reported in Section~\ref{subsec:euler2d_periodic}, along several contact discontinuities and shear layers LGNO develops finer roll-up and vortical structures than the WENO-Z reference on the same mesh; as there, this should be read not as LGNO being more accurate than the reference but as its lower long-time numerical dissipation preserving sharper localized features. These tests reinforce that LGNO is more stable than FNO and, even on these unseen shock-interaction configurations, less dissipative than the WENO-Z scheme that generated its training data.

\clearpage
\newpage

\shockrhofig
{outflow}
{contour}
{outflow}
{contours}
{\textbf{Density contours at \(T=0.25\) for the four outflow shock tests.} Rows S1--S4 are the four standard quadrant shock-interaction configurations, and columns show the WENO-Z reference, LGNO, and the FNO baseline. The contours give a direct check of shock and contact locations, which is relevant for assessing possible shock-speed errors in learned rollouts: LGNO closely matches the reference shock locations, whereas the FNO baseline develops spurious contour oscillations and distorted patterns.}
{25}

\shockrhofig
{outflow}
{colormap}
{outflow}
{colormaps}
{\textbf{Density fields at \(T=0.25\) for the four outflow shock tests.} Rows are S1--S4 and columns show the WENO-Z reference, LGNO, and the FNO baseline. LGNO remains close to the reference, while the FNO baseline exhibits visible unphysical artifacts, including localized overshoots, oscillatory structures, and large-scale distortions of the density field.}
{25}

\shockrhofig
{periodic}
{contour}
{periodic}
{contours}
{\textbf{Density contours at \(T=0.5\) for the four periodic shock tests.} Rows S1--S4 are the four periodized quadrant configurations, and columns show the WENO-Z reference, LGNO, and the FNO baseline. After the longer rollout, LGNO retains the dominant shock and contact geometry, whereas the FNO baseline produces much denser spurious contours. Notably, along contacts and shear layers LGNO is visibly less dissipative than the WENO-Z reference, resolving finer roll-up and vortical structures.}
{50}

\shockrhofig
{periodic}
{colormap}
{periodic}
{colormaps}
{\textbf{Density fields at \(T=0.5\) for the four periodic shock tests.} Rows are S1--S4 and columns show the WENO-Z reference, LGNO, and the FNO baseline. The panels use individual color scales, so the comparison focuses on the geometry of the density field and the presence of unphysical oscillations rather than on direct amplitude matching. Along several contacts and shear layers, LGNO resolves finer roll-up and vortical structures than the WENO-Z reference, reflecting its lower long-time numerical dissipation.}
{50}

\FloatBarrier

\section{Implementation Details and Computational Cost}
\label{apx:implementation_and_cost}

\subsection{Model and Training Configuration}
\label{apx:implementation_details}

The neural operators are parameterized in a residual form,
\[
u^{n+1}=u^n+\Delta t\,\Phi_\theta(u^n),
\]
so that the network represents the discrete increment over one time step. In the periodic case, the predicted increment is projected onto the zero-mean subspace before the update for global conservations.

Unless otherwise specified, all models use width \(64\) and four operator blocks. Fixed coordinate channels are appended to the physical input fields. For periodic problems, the coordinates are represented by periodic features, e.g., \(\sin(2\pi x)\) and \(\cos(2\pi x)\) in one space dimension, with the analogous features used for each coordinate direction in two space dimensions. For non-periodic boundary conditions, the raw spatial coordinates are used directly. The retained Fourier modes and the corresponding numbers of trainable parameters are summarized in Table~\ref{tab:implementation_model_sizes}. The benchmark-specific numerical setups, including the reference solvers, initial data, and snapshot spacing, are given in~\ref{apx:numerical_setup}.

\begin{table}[htbp]
\centering
\caption{Model sizes, retained Fourier modes, and local convolution kernel sizes used in the numerical comparisons.}
\label{tab:implementation_model_sizes}
\small
\begin{tabular}{lccccc}
\toprule
Problem & LGNO modes & FNO modes & LGNO kernel & LGNO params & FNO params \\
\midrule
Linear advection & \(16\) & \(24\) & \(5\) & \(841{,}729\) & \(811{,}970\) \\
Burgers & \(16\) & \(24\) & \(5\) & \(841{,}729\) & \(811{,}777\) \\
SWE & \(16\) & \(24\) & \(5\) & \(841{,}858\) & \(811{,}970\) \\
Euler & \(16\) & \(24\) & \(5\) & \(841{,}987\) & \(812{,}163\) \\
Burgers2D & \(16\times 16\) & \(24\times 24\) & \(5\) & \(17{,}750{,}017\) & \(37{,}774{,}209\) \\
Euler2D periodic & \(24\times 24\) & \(24\times 24\) & \(5\) & \(39{,}508{,}356\) & \(37{,}774{,}788\) \\
Euler2D outflow (primitive) & \(24\times 24\) & \(24\times 24\) & \(7\) & \(41{,}744{,}772\) & \(37{,}774{,}660\) \\
\bottomrule
\end{tabular}
\end{table}

Unless otherwise noted, the LGNO and FNO models are trained with the AdamW optimizer~\cite{adamw} using an initial learning rate of \(10^{-3}\), weight decay \(10^{-4}\), and a cosine annealing learning-rate schedule~\cite{cosine_annealing}; benchmark-specific deviations are noted in~\ref{apx:numerical_setup}. All experiments are run on NVIDIA A100 GPUs on the NERSC Perlmutter system.

\subsection{Computational Cost Comparison}
\label{apx:computational_cost}

We compare the wall-clock cost of WENO-Z, LGNO, and the FNO baseline on the same evaluation meshes, to give a rough sense of inference cost rather than a careful benchmark. The neural models use the same configurations as in the main text and Table~\ref{tab:implementation_model_sizes}. Timings are measured on a single core of an AMD EPYC 7713 CPU for WENO-Z and on an NVIDIA A100 GPU for the neural models. The WENO-Z solver is a vectorized NumPy implementation; it is not performance-optimized (no multithreading or compiled kernels) and uses a component-wise reconstruction without characteristic decomposition. A characteristic-wise WENO-Z, as commonly used for nonlinear systems, would be a few times more expensive, whereas a compiled or parallel implementation would be substantially faster; the reported WENO-Z times therefore do not represent an optimized solver.

For each problem, we advance one representative test sample to the same final time: WENO-Z uses CFL-limited explicit time stepping, while the learned models are applied autoregressively with the fixed learned time step. The evaluation settings are summarized in Table~\ref{tab:computational_cost_setup}. WENO-Z is timed on the same spatial mesh as the learned models.

Table~\ref{tab:computational_cost} reports the corresponding wall-clock times and speedups. Both learned models advance the solution faster than this WENO-Z implementation on the same mesh, with the gap growing from modest factors for the inexpensive one-dimensional scalar problems to one to three orders of magnitude for the two-dimensional Euler tests. Several factors contribute to these numbers and pull in different directions: the learned models run on a GPU while WENO-Z runs on a single CPU core; for the two Euler2D cases the neural models are evaluated in single precision with TF32 tensor-core acceleration, which raises A100 throughput by roughly an order of magnitude relative to double precision, whereas WENO-Z remains in double precision; and the WENO-Z baseline is unoptimized and component-wise, as noted above. The reported speedups should therefore be read as practical, order-of-magnitude guidance on deployment cost, not as a controlled comparison between two optimized codes, since neither code was tuned for performance. A careful, hardware- and precision-normalized comparison is left to future work.

The FNO baseline is generally faster than LGNO in raw inference time, since LGNO adds local operations to the global spectral branch. As shown in the preceding experiments, this additional cost buys improved stability and accuracy in shock-dominated rollouts. Both learned models remain considerably cheaper to evaluate than the present WENO-Z implementation.

\begin{table}[!htbp]
\centering
\caption{Evaluation settings used in the computational cost comparison. The learned models are evaluated autoregressively with the fixed time step \(\Delta t\) shown below. WENO-Z is advanced to the same final time using CFL-limited explicit time stepping with the listed CFL number.}
\label{tab:computational_cost_setup}
\small
\begin{tabular}{llccccc}
\toprule
Problem & Mesh & Boundary & Learned \(\Delta t\) & Learned steps & Final time & WENO CFL \\
\midrule
Linear advection 
& \(256\) & periodic & \(0.05\) & \(60\) & \(3.0\) & \(0.4\) \\

Burgers 
& \(256\) & periodic & \(0.05\) & \(80\) & \(4.0\) & \(0.4\) \\

Shallow water 
& \(256\) & periodic & \(0.05\) & \(60\) & \(3.0\) & \(0.4\) \\

Euler 
& \(256\) & periodic & \(0.05\) & \(20\) & \(1.0\) & \(0.4\) \\

Burgers2D 
& \(128^2\) & periodic & \(0.05\) & \(60\) & \(3.0\) & \(0.4\) \\

Euler2D
& \(256^2\) & outflow & \(0.01\) & \(40\) & \(0.4\) & \(0.45\) \\

Euler2D
& \(256^2\) & periodic & \(0.01\) & \(50\) & \(0.5\) & \(0.45\) \\
\bottomrule
\end{tabular}
\end{table}

\begin{table}[!htbp]
\centering
\caption{Wall-clock times for advancing one test sample to the specified final time using WENO-Z, LGNO, and the FNO baseline on the same evaluation meshes. WENO-Z runs as a single-core, unoptimized, component-wise NumPy implementation (no characteristic decomposition) in double precision, while the neural models run on an A100 GPU in double precision, except for the two Euler2D tests, which use single precision with TF32. Speedups are relative to this WENO-Z implementation on the same mesh and time interval, computed from the unrounded timings. Because neither code is performance-optimized and the settings differ in hardware and precision, these numbers should be read as general guidance on deployment cost rather than a careful benchmark.}
\label{tab:computational_cost}
\small
\begin{tabular}{lrrrrr}
\toprule
Problem & WENO-Z (s) & LGNO (s) & FNO (s) 
& LGNO speedup & FNO speedup \\
\midrule
Linear advection 
& \(0.927\) & \(0.317\) & \(0.179\) 
& \(2.9\times\) & \(5.2\times\) \\

Burgers 
& \(0.661\) & \(0.471\) & \(0.267\) 
& \(1.4\times\) & \(2.5\times\) \\

Shallow water 
& \(9.545\) & \(0.325\) & \(0.180\) 
& \(29.4\times\) & \(53.0\times\) \\

Euler 
& \(3.862\) & \(0.115\) & \(0.066\) 
& \(33.7\times\) & \(58.8\times\) \\

Burgers2D 
& \(4.021\) & \(0.471\) & \(0.278\) 
& \(8.5\times\) & \(14.4\times\) \\

Euler2D outflow 
& \(511.263\) & \(0.698\) & \(0.244\) 
& \(732.8\times\) & \(2098.0\times\) \\

Euler2D periodic 
& \(809.503\) & \(0.484\) & \(0.310\) 
& \(1672.5\times\) & \(2611.3\times\) \\
\bottomrule
\end{tabular}
\end{table}

\FloatBarrier

\subsection{Numerical Setups}
\label{apx:numerical_setup}

The training datasets are generated from reference trajectories computed on fine grids (e.g., \(512^2\)) and then transferred to the learning grids (e.g., \(256^2\)) by conservative averaging. All benchmarks use periodic boundary conditions except the two-dimensional Euler outflow case.

\subsubsection{Linear Advection Equation}

For the linear advection equation~\eqref{eq:convection_pde}, the training data consist of \(10{,}000\) input-target pairs on a \(256\)-cell grid (\(\Delta t=0.05\)), generated from \(1{,}000\) trajectories using \(10\) consecutive transitions each; batch size \(64\).

The initial condition is sampled from a mixture of discontinuous and smooth periodic profiles, with \(60\%\) discontinuous and \(40\%\) smooth. The discontinuous component is piecewise constant,
\[
u_0(x) = \sum_{j=1}^{m} c_j \mathbf{1}_{I_j}(x),
\]
where \(m\in\{2,3\}\). The periodic partition \(\{I_j\}_{j=1}^{m}\) is determined by one or two random breakpoints, with a minimum cyclic separation of \(0.15\), and the segment values \(c_j\) are sampled independently from the uniform distribution on \([-0.8,0.8]\). The smooth component is a Fourier series,
\[
\tilde u_0(x) = \sum_{k=1}^{K} \bigl(a_k \cos(2\pi kx) + b_k \sin(2\pi kx)\bigr),
\quad K \in \{2,3\},
\]
with raw coefficients decaying as \(k^{-2}\); each profile is then normalized by its maximum absolute value and multiplied by an amplitude sampled uniformly from \([0.5,1.0]\).

\subsubsection{One-Dimensional Burgers Equation}

For the Burgers equation~\eqref{eq:burgers1d_pde}, the training data consist of \(10{,}000\) input-target pairs on a \(256\)-cell grid (\(\Delta t=0.05\)), generated from \(500\) trajectories evolved to \(T=1\) using \(20\) consecutive transitions each; batch size \(64\).

The initial condition is generated from a Fourier series,
\[
\tilde u_0(x) =
\sum_{k=1}^{5}
\bigl(a_k \cos(2\pi kx) + b_k \sin(2\pi kx)\bigr),
\]
where \(a_k\) and \(b_k\) are independent Gaussian random variables with mean zero and variances proportional to \(k^{-4}\). The profile is then normalized by its maximum absolute value and multiplied by an amplitude \(A\) sampled uniformly from \([0.5,1.0]\). A constant mean shift \(B\), sampled uniformly from \([-\tfrac12 A,\tfrac12 A]\), is then added. This yields smooth initial data that develop shocks during the evolution.

\subsubsection{One-Dimensional Shallow Water Equations}
For the one-dimensional shallow water equations~\eqref{eq:swe_pde}, the training data consist of \(10{,}000\) input-target pairs on a \(256\)-cell grid (\(\Delta t=0.05\)), generated from \(500\) trajectories evolved to \(T=1\) using \(20\) consecutive transitions each; batch size \(64\).

The initial condition is sampled from a mixture of three components. For clarity, we describe the sampling procedure in primitive variables \((h,u)\). The first component (\(30\%\)) is a smooth bidirectional wave, obtained by perturbing the Riemann invariants of the shallow water system,
\[
w_{\pm}(x)=u(x)\pm 2\sqrt{g h(x)}.
\]
Let \((h_0,u_0)=(1,0)\), \(c_0=\sqrt{g h_0}\), and base invariants \(w_{\pm}^{(0)}=u_0\pm 2c_0\). We then sample two independent smooth Fourier profiles \(s_+(x)\) and \(s_-(x)\), normalized to have unit maximum absolute value, and define
\[
w_{\pm}(x)=w_{\pm}^{(0)}+A s_{\pm}(x),
\qquad
A/c_0\in[0.08,0.22].
\]
The number of retained Fourier modes is sampled from \(\{2,3,4\}\), with higher modes damped in the random series. The primitive variables are recovered from the perturbed invariants by
\[
h(x)=\frac{1}{g}\left(\frac{w_+(x)-w_-(x)}{4}\right)^2,
\qquad
u(x)=\frac{w_+(x)+w_-(x)}{2}.
\]
A mild positivity safeguard is applied to avoid near-dry states.

The second component (\(40\%\)) is a coupled smooth state of the form
\[
h(x) = h_0 + \varepsilon_h s_h(x),\qquad
u(x) = u_0 + A_u s_u(x),
\]
where \(\varepsilon_h\) is sampled uniformly from \([0.10,0.22]\) and \(A_u/c_0\) from \([0.15,0.45]\). The functions \(s_h\) and \(s_u\) are Fourier series with decay slopes sampled from \([1,3]\), and the largest retained Fourier mode is sampled between \(2\) and \(4\).

The third component (\(30\%\)) is a periodic smoothed Riemann state,
\[
h(x) = h_{\mathrm{out}} + (h_{\mathrm{in}}-h_{\mathrm{out}})\,\chi(x),\qquad
u(x) = u_{\mathrm{out}} + (u_{\mathrm{in}}-u_{\mathrm{out}})\,\chi(x),
\]
where
\[
\chi(x)
=
\frac12
\left(
1+\tanh\left(\frac{w-|d(x,x_c)|}{\delta}\right)
\right).
\]
Here \(d(x,x_c)\) is the periodic distance from \(x\) to the center \(x_c\). The parameters are sampled uniformly: \(h_{\mathrm{in}}\in[1.2,2.0]\), \(h_{\mathrm{out}}\in[0.2,0.9]\), \(u_{\mathrm{in}},u_{\mathrm{out}}\in[-0.8,0.8]\), \(w\in[0.08,0.22]\), and \(\delta\in[0.01,0.04]\).

\subsubsection{One-Dimensional Euler Equations}

For the one-dimensional Euler equations~\eqref{eq:euler1d_pde}, the training data consist of \(10{,}000\) input-target pairs on a \(256\)-cell grid (\(\Delta t=0.05\)), generated from \(500\) trajectories evolved to \(T=1\) using \(20\) consecutive transitions each; batch size \(64\).

The initial condition is sampled from a mixture of random two-state Riemann data and smooth perturbations, with \(80\%\) Riemann and \(20\%\) smooth. The Riemann component has constant left and right states meeting at a discontinuity \(x_0\) sampled from the middle \(60\%\) of the domain, with \(\rho,p\in[0.05,1.20]\) and \(u\in[-1,1]\) on each side; the two states must differ by at least \(30\%\) in density and pressure and by at least \(0.3\) in velocity. The smooth component has
\[
\rho(x) = \rho_0\bigl(1+\alpha_{\rho} S_{\rho}(x)\bigr),
\quad
p(x) = p_0\bigl(1+\alpha_{p} S_{p}(x)\bigr),
\quad
u(x) = u_0 + \alpha_{u} S_{u}(x),
\]
where \(S_{\rho},S_p,S_u\) are Fourier series with \(\|S_\bullet\|_{L^\infty}\le 1\) and \(K\in\{1,2,3\}\) modes, base state \(\rho_0,p_0\in[0.30,1.10]\), \(u_0\in[-1,1]\), and amplitudes \(\alpha_{\rho},\alpha_p\in[0.12,0.65]\), \(\alpha_u\in[0.10,0.70]\).

\subsubsection{Two-Dimensional Burgers Equation}

For the two-dimensional Burgers equation~\eqref{eq:burgers2d_pde}, the training data consist of \(10{,}000\) input-target pairs on a \(128\times128\) grid (\(\Delta t=0.05\)), generated from \(500\) trajectories evolved to \(T=1\) using \(20\) consecutive transitions each; batch size \(16\).

The initial condition is generated from a  Fourier series,
\[
\tilde u_0(x,y)
= \sum_{k_x=0}^{K_x} \sum_{k_y=0}^{K_y}
\widehat{u}_{k_x,k_y}
\cos\bigl(2\pi(k_x x + k_y y) + \phi_{k_x,k_y}\bigr),
\]
where \(K_x\) and \(K_y\) are sampled independently from \(\{0,1,2,3,4\}\). For each nonzero mode \((k_x,k_y)\), the coefficient magnitude is proportional to \((k_x^2+k_y^2)^{-0.75}\) and the phase \(\phi_{k_x,k_y}\) is sampled uniformly from \([0,2\pi]\). The resulting field is normalized by its maximum absolute value and multiplied by an amplitude \(A\) sampled uniformly from \([0.4,1.0]\). It is then shifted by a constant \(B\), sampled uniformly from \([-\tfrac12 A,\tfrac12 A]\).

\subsubsection{Two-Dimensional Euler Equations with Outflow Boundary Conditions}
\label{app:euler2d_outflow_details}

For the two-dimensional Euler equations~\eqref{eq:euler2d_system}--\eqref{eq:euler2d_pressure} with outflow boundary conditions, the training data consist of \(10{,}000\) input-target pairs on a \(256\times256\) grid (\(\Delta t=0.01\)), generated from \(250\) trajectories evolved to \(T=0.4\) using \(40\) consecutive transitions each; the reference solutions are computed on a \(512\times512\) grid and transferred to the \(256\times256\) grid. The models use batch size \(16\), initial learning rate \(5\times10^{-4}\), and cosine decay over \(500\) epochs.

The initial condition is a piecewise constant quadrant state,
\[
U_0(x,y)=
\begin{cases}
U_1, & x \ge x_0,\ y \ge y_0,\\
U_2, & x < x_0,\ y \ge y_0,\\
U_3, & x \ge x_0,\ y < y_0,\\
U_4, & x < x_0,\ y < y_0,
\end{cases}
\]
where the split point \((x_0,y_0)\) is sampled in the interior of the domain. Each quadrant state is sampled independently in primitive variables \((\rho,u,v,p)\). The density \(\rho_j\) is sampled from the log-uniform distribution on \([0.08,1.0]\), and the pressure \(p_j\) is sampled from the log-uniform distribution on \([0.04,0.5]\). The velocity components \(u_j,v_j\) are sampled uniformly from \([-0.5,0.5]\), and the state is required to satisfy the Mach bound \(\sqrt{u_j^2+v_j^2}/\sqrt{\gamma p_j/\rho_j} < 3.5\); if violated, the velocity vector is rescaled until the bound holds.

For the outflow experiment, the models are trained in primitive variables \((\rho,u,v,p)\), in which case we do not enforce global conservation. Instead, we add a softplus activation to guarantee positivity. 
Given \(\mathbf u^n=(\rho^n,u^n,v^n,p^n)\), the model first produces an unconstrained prediction
\[
\widetilde{\mathbf u}^{n+1}
=
\mathbf u^n+\Delta t\,\Phi_\theta(\mathbf u^n;\Delta t).
\]
A shifted softplus map \(\operatorname{softplus}_{\beta}(z):=\tfrac{1}{\beta}\log\!\bigl(1+e^{\beta z}\bigr)\) is then applied to density and pressure for positivity preservation,
\[
\rho^{n+1}=\rho_{\min}+\operatorname{softplus}_{\beta}\!\bigl(\widetilde\rho^{n+1}-\rho_{\min}\bigr),
\qquad
p^{n+1}=p_{\min}+\operatorname{softplus}_{\beta}\!\bigl(\widetilde p^{n+1}-p_{\min}\bigr),
\]
while the velocities pass through unchanged, \(u^{n+1}=\widetilde u^{n+1}\) and \(v^{n+1}=\widetilde v^{n+1}\). Here \(\rho_{\min}=10^{-6}\), \(p_{\min}=10^{-8}\), and \(\beta=100\). LGNO uses an outflow boundary operator with local kernel size \(7\).

\subsubsection{Two-Dimensional Euler Equations with Periodic Boundary Conditions}

For the periodic case, the training data consist of \(10{,}000\) input-target pairs on a \(256\times256\) grid (\(\Delta t=0.01\)), generated from \(200\) trajectories evolved to \(T=0.5\) using \(50\) consecutive transitions each; reference solutions are computed on a \(512\times512\) grid and transferred to the \(256\times256\) grid. Batch size \(16\).

The initial condition is sampled from a mixture of discontinuous quadrant states and smooth periodic states, with \(80\%\) discontinuous and \(20\%\) smooth. The discontinuous component uses the same piecewise constant quadrant form as in the outflow case, now with \(\rho_j,p_j\) log-uniform on \([0.08,1.0]\), \(u_j,v_j\) uniform on \([-1,1]\), and the same Mach bound \(\sqrt{u_j^2+v_j^2}/\sqrt{\gamma p_j/\rho_j}<3.5\) (rescaling the velocity if violated).

The smooth component has
\[
\log \rho = \log \bar{\rho} + S_{\rho},
\quad
\log p = \log \bar{p} + S_{p},
\quad
u = \bar{u} + S_{u},
\quad
v = \bar{v} + S_{v},
\]
where each \(S_{\bullet}(x,y)\) is a low-frequency Fourier sum with \(2\), \(3\), or \(4\) distinct modes selected from \(\{0,1,2,3\}^2 \setminus \{(0,0)\}\). The amplitudes of \(S_{\rho}\) and \(S_p\) are sampled uniformly from \([0.08,0.22]\), while the amplitudes of \(S_u\) and \(S_v\) are sampled uniformly from \([0.05,0.35]\). If the resulting state exceeds Mach number \(3.5\), the velocity field is uniformly rescaled to satisfy the same admissibility threshold. 


\bibliography{references.bib}

\end{document}